\documentclass[11pt,twoside,english]{article}
\usepackage[T1]{fontenc}
\usepackage[utf8]{inputenc}
\usepackage{geometry}
\usepackage{color}
\usepackage{babel}
\usepackage{array}
\usepackage{units}
\usepackage{amsmath}
\usepackage{amsthm}
\usepackage{amssymb}
\usepackage{stackrel}
\usepackage{graphicx}
\usepackage[all]{xy}
\PassOptionsToPackage{normalem}{ulem}
\usepackage{ulem}
\usepackage[unicode=true,pdfusetitle,
 bookmarks=true,bookmarksnumbered=false,bookmarksopen=false,
 breaklinks=false,pdfborder={0 0 1},backref=page,colorlinks=true]
 {hyperref}
\hypersetup{
 linkcolor=blue,citecolor=blue}

\makeatletter

\providecommand{\tabularnewline}{\\}

\numberwithin{equation}{section}
\numberwithin{figure}{section}
\theoremstyle{plain}
\newtheorem{thm}{\protect\theoremname}[section]
\theoremstyle{plain}
\newtheorem*{thm*}{\protect\theoremname}
\theoremstyle{plain}
\newtheorem{conjecture}[thm]{\protect\conjecturename}
\theoremstyle{definition}
\newtheorem{defn}[thm]{\protect\definitionname}
\theoremstyle{plain}
\newtheorem{cor}[thm]{\protect\corollaryname}
\theoremstyle{plain}
\newtheorem{question}[thm]{\protect\questionname}
\theoremstyle{plain}
\newtheorem{lem}[thm]{\protect\lemmaname}
\theoremstyle{remark}
\newtheorem{rem}[thm]{\protect\remarkname}
\theoremstyle{remark}
\newtheorem{claim}[thm]{\protect\claimname}
\theoremstyle{plain}
\newtheorem{prop}[thm]{\protect\propositionname}
\theoremstyle{definition}
\newtheorem{example}[thm]{\protect\examplename}
\theoremstyle{plain}
\newtheorem{fact}[thm]{\protect\factname}

\usepackage{babel}

\usepackage{ytableau}

\usepackage{multirow}

\usepackage{appendix}

\usepackage{bbm}
\newenvironment{cellvarwidth}[1][t]
    {\begin{varwidth}[#1]{\linewidth}}
    {\@finalstrut\@arstrutbox\end{varwidth}}
\usepackage{varwidth}

\usepackage{xcolor}
\usepackage{tikz}
\usetikzlibrary{calc}
\usetikzlibrary{decorations.markings}
\usetikzlibrary{arrows.meta}
\usetikzlibrary{bending}

\tikzset{     
	mynode/.style={circle, draw=black, fill=white, inner sep=1.5pt},
	vertex/.style={circle, fill=white, inner sep=1.5pt}, 
	labelnode/.style={font=\small, text=red}, 
	edgelabel/.style={font=\small}, 
	branchlabel/.style={font=\tiny, anchor=west} 
}

\makeatother

\providecommand{\claimname}{Claim}
\providecommand{\conjecturename}{Conjecture}
\providecommand{\corollaryname}{Corollary}
\providecommand{\definitionname}{Definition}
\providecommand{\examplename}{Example}
\providecommand{\factname}{Fact}
\providecommand{\lemmaname}{Lemma}
\providecommand{\propositionname}{Proposition}
\providecommand{\questionname}{Question}
\providecommand{\remarkname}{Remark}
\providecommand{\theoremname}{Theorem}

\begin{document}
\global\long\def\F{\mathbb{\mathbb{\mathbf{F}}}}%
 
\global\long\def\rk{\mathbb{\mathrm{rk}}}%
 
\global\long\def\crit{\mathbb{\mathrm{Crit}}}%
 
\global\long\def\Hom{\mathrm{Hom}}%
 
\global\long\def\defi{\stackrel{\mathrm{def}}{=}}%
 
\global\long\def\tr{{\cal T}r }%
 
\global\long\def\id{\mathrm{id}}%
 
\global\long\def\Aut{\mathrm{Aut}}%
 
\global\long\def\wl{w_{1},\ldots,w_{\ell}}%
 
\global\long\def\alg{\le_{\mathrm{alg}}}%
 
\global\long\def\ff{\stackrel{*}{\le}}%
 
\global\long\def\mobius{M\dacute{o}bius}%
 
\global\long\def\chimax{\chi^{\mathrm{max}}}%
 
\global\long\def\fq{\mathbb{F}_{q}}%
 
\global\long\def\gl{\mathrm{GL}}%
 
\global\long\def\glm{\mathrm{GL}_{m}\left(\k\right)}%
 
\global\long\def\gln{\mathrm{GL}_{N}\left(q\right)}%
 
\global\long\def\k{K}%
 
\global\long\def\fix{\mathrm{fix}}%
 
\global\long\def\cl{\mathrm{cl}}%
 
\global\long\def\scl{\mathrm{scl}}%
 
\global\long\def\sql{\mathrm{sql}}%
 
\global\long\def\sqlh{\mathrm{sql^{*}}}%
 
\global\long\def\ssql{\mathrm{ssql}}%
 
\global\long\def\sp{\mathrm{s\pi}}%
 
\global\long\def\spm{\mathrm{s\pi^{\left(m\right)}}}%
 
\global\long\def\spq{\mathrm{s\pi_{q}}}%
 
\global\long\def\spk{\mathrm{s\pi_{K}}}%
 
\global\long\def\spf{\mathrm{s\pi^{\phi}}}%
 
\global\long\def\decomp{{\cal D}\mathrm{ecomp}_{B}}%
 
\global\long\def\decompt{{\cal D}\mathrm{ecomp}_{B}^{3}}%
 
\global\long\def\algdecomp{{\cal D}\mathrm{ecomp}_{\mathrm{alg}}}%
 
\global\long\def\algdecompt{{\cal D}\mathrm{ecomp}_{\mathrm{alg}}^{3}}%
 
\global\long\def\irr{\mathrm{Irr} }%
 
\global\long\def\wh{\mathrm{Wh} }%
 
\global\long\def\cod{\mathrm{codom} }%
 
\global\long\def\Sp{\mathrm{Sp}}%
 
\global\long\def\s{\mathrm{S}}%
 
\global\long\def\U{\mathrm{U}}%
 
\global\long\def\O{\mathrm{O}}%
 
\global\long\def\arrm{\mathrm{\overrightarrow{\mu}}}%
 
\global\long\def\arrn{\mathrm{\overrightarrow{\nu}}}%
 
\global\long\def\I{{\cal I}}%
 
\global\long\def\irr{\mathrm{Irr}}%
 
\global\long\def\stairr{\mathrm{StabIrr}}%
 
\global\long\def\triv{\mathrm{triv}}%
\global\long\def\supp{\mathrm{Supp}}%

\title{Stable Invariants of Words from Random Matrices}
\author{Doron Puder~~~~and~~~~Yotam Shomroni\\
\\
with an Appendix with Danielle Ernst-West and Matan Seidel}
\maketitle
\begin{abstract}
Let $w$ be a word in a free group. A few years ago, Magee and the
first named author discovered that the stable commutator length ($\scl$)
of $w$, a well-known topological invariant, can also be defined in
terms of certain Fourier coefficients of $w$-random unitary matrices
\cite{MPunitary}. But the random-matrix side of this equality can
be naturally tweaked by considering $w$-random permutations, $w$-random
orthogonal matrices and so on, to produce new invariants for any given
word. Are these invariants new? interesting? Do they admit an intrinsic
topological description as in the case of $w$-random unitaries and
$\scl$? 

The current paper formalizes the definition of these invariants coming
from $w$-random matrices, answers the above questions in certain
cases involving generalized symmetric groups, and poses detailed conjectures
in many others. In particular, we present a plethora of topological,
combinatorial and algebraic invariants of words which play, or are
at least conjectured to play, a similar role to the one played by
$\scl$ in the above-mentioned result. Among others, these invariants
include two invariants recently defined by Wilton \cite{wilton2024rational}:
the stable primitivity rank and a non-oriented analog of $\scl$.
\end{abstract}
\tableofcontents{}

\section{Introduction\label{sec:Introduction}}

Let $\F$ be a free group. Every word $w\in\F$ induces a measure,
called the $w$-measure, on every compact group $G$ by substituting
the letters of $w$ with independent, Haar-random elements of $G$.
For example, if $w=abab^{-2}$, a $w$-random element in $G$ is $ghgh^{-2}$
where $g$ and $h$ are independent, Haar-random elements of $G$.\footnote{Whenever we spell out a word, unless stated otherwise, we use letters
which are assumed to be distinct elements of a basis of $\F$.} 

For some very special words, the corresponding $w$-measures follow
a very structured pattern. Most notably, this is true in the case
of orientable surface words, namely, when $w_{g}=[a_{1},b_{1}]\cdots[a_{g},b_{g}]$
for some $g\ge0$: Let $\irr(G)$\marginpar{$\protect\irr(G)$} denote
the set of complex irreducible characters of $G$. Then for every
compact $G$ and every $\xi\in\irr(G)$, the expected value of $\xi$
under the $w_{g}$-measure on $G$ satisfies 
\begin{equation}
\mathbb{E}_{w_{g}}\left[\xi\right]=\frac{1}{\left(\dim\xi\right)^{2g-1}}.\label{eq:E_w=00005Bxi=00005D for orientable surface words}
\end{equation}
Similarly, for a non-orientable surface word $u_{g}=a_{1}^{2}\cdots a_{g}^{2}$,
$\mathbb{E}_{u_{g}}[\xi]=\frac{(\mathrm{FS_{\xi}})^{g}}{(\dim\xi)^{g-1}}$,
where $\mathrm{FS}_{\xi}\in\left\{ -1,0,1\right\} $ is the Frobenius-Schur
indicator of $\xi$. (These results have their seeds in the works
of Frobenius \cite{frobenius1896gruppencharaktere} and Frobenius-Schur
\cite{frobenius1906rellen}, and see, e.g., \cite[\S 1 and \S 2]{MPsurfacewords}.)
We remark that as word measures are class functions, the expected
values $(\mathbb{E}_{w}[\xi])_{\xi\in\irr(G)}$ completely determine
the $w$-measure on $G$ \cite[Fact 2.8]{magee2015word}.

The above results about surface words rely heavily on the fact that
every letter appears at most twice in the word. There is no parallel
statement for general words: the expected values of irreducible characters
certainly do not obey such clean formulas. 

However, in the context of \emph{stable }irreducible characters, some
aspects of this behavior do hold for general words. Consider a sequence
of groups $G_{\bullet}=\left\{ G_{N}\right\} _{N\ge0}$ with a stable
representation theory in the sense of \cite{church2015fi,sam_snowden_2015}\footnote{We say a few more general words about this theory in $\S$\ref{subsec:Stable-representation-theory}.
Along the paper we will mainly describe it case-by-case.} (one example to have in mind is the symmetric groups $\s_{\bullet}=\left\{ \s_{N}\right\} _{N\ge0}$).
Such a sequence admits stable irreducible characters: these are certain
sequences $\chi=\left\{ \chi_{N}\right\} _{N\ge N_{0}}$ with $\chi_{N}\in\irr(G_{N})$
(for $\s_{\bullet}$, as elaborated in $\S$\ref{subsec:Stable-primitivity-rank and S}
below, such $\chi$ is a sequence of Young diagrams with a fixed structure
outside the first row, e.g., the diagram $(n-1,1)$ corresponding
to the number of fixed points minus 1). All stable sequences of groups
mentioned in this paper, and possibly all stable sequences in general,
satisfy that for every word $1\ne w\in\F$ and every non-trivial stable
character $\chi$, there is some \marginpar{$\beta(w,\chi)$}$\beta(w,\chi)\in\mathbb{Q}_{\ge0}\cup\left\{ \infty\right\} $
such that $\mathbb{E}_{w}\left[\chi_{N}\right]$ is of order $(\dim\chi_{N})^{-\beta(w,\chi)}$.\footnote{\label{fn:E=00005Bchi=00005D is of order dim^-beta}What we mean here
is that there is a constant $c\ne0$ with $\mathbb{E}_{w}\left[\chi_{N}\right]=(\dim\chi_{N})^{-\beta(w,\chi)}\left(c+o_{N}(1)\right)$.
We stress that $\beta(w,\chi)$ may obtain non-integer values: see
Example \ref{exa:fractional beta in S_n}.} For example, in the case of $\s_{\bullet}$, this follows from \cite[Cor.~6.9 and Prop.~B.2]{hanany2020word}.
Granted, even when restricted to stable characters, this is still
much weaker than \eqref{eq:E_w=00005Bxi=00005D for orientable surface words},
in the sense that for surface words the values $\beta(w,\chi)$ have
very small dependency on $\chi$: $\beta(w_{g},\chi)=2g-1$ and $\beta(u_{g},\chi)\in\left\{ g-1,\infty\right\} $
for all $\chi$. 

A much stronger result, significantly closer to \eqref{eq:E_w=00005Bxi=00005D for orientable surface words},
was obtained in the work \cite{MPunitary} of Magee and the first
named author. They showed that if one restricts attention to the sequence
of unitary groups $U(\bullet)=\left\{ \U(N)\right\} _{N}$ (which
admits stable representation theory -- see $\S$\ref{subsec:Stable-commutator-length-and-U})
and considers \emph{all stable polynomial }irreducible characters,
then for any given $w\ne1$, all corresponding $\beta(w,\chi)$'s
are not only bounded away from zero, but their infimum coincides with
a well-known, seemingly unrelated, topological invariant of words:
the stable commutator length of $w$, or $\scl(w)$ (defined in \eqref{eq:scl, old def}
below). 
\begin{thm}
\label{thm:MP scl for U(N)}\cite{MPunitary} Let $\I_{\U,\mathrm{poly}}$
denote the set of all stable \uline{polynomial} irreducible characters
of $\U(\bullet)$. Then for every $1\ne w\in\F$,
\begin{equation}
2\cdot\scl(w)=\inf_{\triv\ne\chi\in\I_{\U,\mathrel{poly}}}\beta(w,\chi).\label{eq:MP-scl}
\end{equation}
\textup{Moreover, the infimum in the right hand side is obtained,
namely, it is, in fact, a minimum.}
\end{thm}

In $\S\S$\ref{subsec:Stable-commutator-length-and-U} and $\S$\ref{sec:scl+U}
we give the precise definitions of the notions appearing in Theorem
\ref{thm:MP scl for U(N)}. 

So Theorem \ref{thm:MP scl for U(N)} gets much closer to \eqref{eq:E_w=00005Bxi=00005D for orientable surface words}:
for an infinite, natural family of stable irreducible characters,
the corresponding $\beta(w,\chi)$'s are all uniformly bounded away
from zero, and with a meaningful tight bound. These bounds on stable
characters often show up in applications. For example, stable irreducible
characters of $\s_{\bullet}$ completely determine the asymptotic
expansion of the Witten zeta functions of $S_{N}$ as $N\to\infty$
\cite[Prop.~2.5]{LiebeckShalev}, thus showing up also in results
around surface groups as in \cite{LiebeckShalev,magee2020asymptotic,magee2025strong}.\footnote{Granted, they are used in these works only when combined with effective
bounds of all, not-necessarily-stable, irreducible representations.} 

Furthermore, Theorem \ref{thm:MP scl for U(N)} exposes a surprising
relation between topology and random matrices: it shows that using
Fourier coefficients of $w$-random unitary matrices, one can recover
an invariant of words which not only has an intrinsic definition unrelated
to $w$-measures, but is also important and useful for its own sake.
Indeed, $\scl$ gained a lot of traction in the mathematical research
-- see, e.g., the book \cite{calegari2009scl}. In particular, Bavard's
duality established its connection with quasi-morphisms and bounded
cohomology \cite{bavard1991longueur}. The $\scl$ of words in free
groups received special attention with various milestones: Culler
\cite{culler1981using} provided an algorithm to compute commutator
length and implicitly established the precise value of $\scl$ for
certain words; Duncan and Howie \cite{duncan1991genus} established
a gap: $2\cdot\scl(w)\ge1$ for every $w\ne1$; and Calegari (most
notably in \cite{calegari2009stable}) enhanced the field much further,
proving, inter alia, that $\scl(w)\in\mathbb{Q}\cup\left\{ \infty\right\} $
for every word $w$ and providing an efficient polynomial-time algorithm
to compute it. 

But the right hand side of \eqref{eq:MP-scl} can be easily tweaked
by replacing $\I_{\U,\mathrm{poly}}$ with \emph{all }stable irreducible
characters of $\U(\bullet)$, or with those of other sequences of
groups, such as the symmetric groups, the orthogonal groups, and so
on. This suggests a way to detect and define relatives of $\scl$.
Such relatives may be useful even if they only admit a fraction of
the importance of $\scl$. Moreover, as we demonstrate in $\S\S$\ref{subsec:gap of spm}
below, we can sometimes derive properties of these $\scl$-relatives
from their other capacity as random-matrix invariants.\medskip{}

The main objective of the current work is to properly single out different
``stable $w$-measure invariants'' à la right hand side of \eqref{eq:MP-scl},
and to study the following natural questions: Are these invariants
interesting -- for instance, are they new? What values can they obtain?
When does the infimum vanish? Do they coincide with some $\scl$-like
topological/combinatorial invariant as in \eqref{eq:MP-scl}? Do these
invariants exhibit similar properties to those exhibited by $\scl$? 

We answer some of these questions in some cases -- most notably for
certain natural classes of stable irreducible characters of generalized
symmetric groups -- thus showing that Theorem \ref{thm:MP scl for U(N)}
is not an isolated phenomenon. We suggest detailed conjectures in
other cases. The $\scl$-like invariants we introduce here include
two invariants recently defined by Wilton: the \emph{stable primitivity
rank} ($\sp$) \cite[Def.~10.6]{wilton2024rational}, and the \emph{stable
square length}\footnote{The term \emph{stable square length} is ours: it is not used in \cite{wilton2024rational}.}
($\ssql$) -- a non-oriented version of $\scl$ \cite[\S10.3]{wilton2024rational}.
We also introduce additional infinite families of stable invariants
of words: the mod-$m$ primitivity rank $\spm$ (Definition \ref{def:spm}),
defined for every $1\ne m\in\mathbb{Z}_{\ge0}$,\footnote{In fact, $\spm$ is a special case of the stable $\phi$-primitivity
rank $\spf$, which is defined for every irreducible character $\phi$
of any compact group -- see \cite[Def.~4.7]{PSh25}.} and the stable $K$-primitivity rank $\spk$, defined for every field
$K$ (Definition \ref{def:spq}).

To keep the length of the current paper manageable, we defer some
of the theory we develop, including some more technical proofs, to
a companion paper \cite{PSh25} and focus here on the definitions
of all invariants, their common properties, many conjectures and open
questions we find intriguing. The following theorems are some highlights
of this work, alongside the new invariants we define and Conjectures
\ref{conj:spi and stable characters of S_N}, \ref{conj:ssql and its role for U,O and Sp}
and \ref{conj:spq and stable characters} which we find very compelling.
\begin{thm*}[Theorems \ref{thm:result on wreath products and modulo-m spi}, \ref{thm:rationality of spm}
and \ref{thm:spm avoids (0,1)}]
 The new family of combinatorial invariants of words $\spm(w)$ (see
Definition \ref{def:spm}) can be defined in terms of stable Fourier
coefficients of $w$-random elements of $C_{m}\wr S_{\bullet}$. Moreover,
these invariants are rational, computable, and admit a gap in $(0,1)$.
\end{thm*}
\begin{thm*}
\cite[Thm.~1.2]{PSh25} The stable primitivity rank $\sp(w)$ can
be defined, too, in terms of stable Fourier coefficients of $w$-random
elements of groups (groups of the form $S_{m}\wr S_{\bullet}).$
\end{thm*}
\begin{thm*}
\cite[Thm.~1.4 and 7.10]{PSh25} Explicit formulas are given for the
expected value of every stable irreducible character of $S_{\bullet}$
and more generally for $G\wr S_{\bullet}$ for an arbitrary compact
group $G$. These formulas demonstrate how stable invariants like
$\sp$ play a role in these stable Fourier coefficients.
\end{thm*}
We proceed on a case by case basis, and state our main results and
conjectures. 

\subsection{Stable commutator length and word measures on $\protect\U\left(\bullet\right)$\label{subsec:Stable-commutator-length-and-U}}

The commutator length and its stable counterpart, the stable commutator
length, are two well-established invariants of words (and of elements
of general groups). The commutator length of $w\in\F$, denoted $\cl\left(w\right)$\marginpar{$\protect\cl\left(w\right)$},
is the length of the shortest product of commutators in $\F$ yielding
$w$ -- see \eqref{eq:cl}. The \emph{stable commutator length} of
$w$, denoted\marginpar{$\protect\scl\left(w\right)$} $\scl\left(w\right)$,
can be defined as
\begin{equation}
\scl\left(w\right)\defi\lim_{n\to\infty}\frac{\cl\left(w^{n}\right)}{n}=\inf_{n\in\mathbb{Z}_{\ge1}}\frac{\cl\left(w^{n}\right)}{n}.\label{eq:scl, old def}
\end{equation}
In particular, $\cl\left(w\right)=\scl\left(w\right)=\infty$ whenever
$w$ does not belong to the commutator subgroup $\left[\F,\F\right]$,
and they are both finite for $w\in\left[\F,\F\right]$. As mentioned
above, these two invariants, particularly in free groups, attracted
a lot of attention in the last decades. More recently, in \cite{MPunitary},
it was discovered that these two invariants play a role in the theory
of word measures on the unitary groups $\U\left(\bullet\right)$.

We describe the stable irreducible characters of $U(\bullet)$ in
$\S\S$\ref{subsec:stable-representation of U}. For now, let us only
mention that there is a natural correspondence between stable characters
of $\U(\bullet)$ and pairs of partitions $\lambda=\left(\lambda^{+},\lambda^{-}\right)$,
and we denote the corresponding stable character $\xi^{\lambda[\bullet]}=\left\{ \xi^{\lambda[N]}\right\} _{N}$.
For example, if $\lambda^{+}=(1)\vdash1$ and $\lambda^{-}=\emptyset\vdash0$,
then $\xi^{\lambda[N]}$ is the standard character -- the trace $\mathrm{tr}_{N}$
-- of $\U(N)$. Recall that $\mathbb{E}_{w}\left[\xi^{\lambda\left[N\right]}\right]$
denotes the expected value of $\xi^{\lambda\left[N\right]}$ under
the $w$-measure. If follows from \cite{MPunitary} that for every
$w$ and every non-trivial stable character (so $\lambda$ is not
the pair of empty partitions), $\mathbb{E}_{w}[\xi^{\lambda[N]}]$
is of order 
\[
\left(\dim\xi^{\lambda[N]}\right)^{-\beta\left(w,\xi^{\lambda[\bullet]}\right)}
\]
for some $\beta(w,\xi^{\lambda[\bullet]})\in\mathbb{Q}\cup\left\{ \infty\right\} $
(see $\S\S$\ref{subsec:stable-representation of U} for details).

For the special case of the standard stable character $\xi^{((1),\emptyset)[\bullet]}=\mathrm{tr}$,
\cite[Cor.~1.8]{MPunitary} shows a connection with the commutator
length of $w$: 
\begin{equation}
\beta\left(w,\mathrm{tr}\right)\ge2\cdot\cl(w)-1,\label{eq:bound with cl}
\end{equation}
namely, $\mathbb{E}_{w}\left[\mathrm{tr}_{N}\right]=O\left(N^{1-2\cdot\cl\left(w\right)}\right)$.
For many words (and possibly for a `generic' word), this bound is
tight, namely, $\beta(w,\mathrm{tr})=2\cdot\cl(w)-1$. Yet, there
are exceptions. For example, $\mathbb{E}_{[x,y][x,z]}[\mathrm{tr}_{N}]=0$
for all $N\ge2$ hence $\beta([x,y][x,z],\mathrm{tr})=\infty$, while
$\cl([x,y][x,z])=2$. 

This demonstrates the limitations of considering a single stable irreducible
character to derive interesting invariants. Indeed, most of the invariants
defined by words measures which are introduced in the current paper
are based on some \emph{infinite} set of stable irreducible characters,
possibly even of different sequences of groups. Usually, we consider
sets which are naturally defined. One example is the set $\I_{\U,\mathrm{poly}}$
from Theorem \ref{thm:MP scl for U(N)}, which consists of all stable
\emph{polynomial} irreducible characters of $\U(\bullet)$: the stable
character $\xi^{\lambda[\bullet]}$ is called polynomial whenever
$\lambda^{-}$ is the empty partition. In this case, the value of
the character $\xi^{\lambda[N]}$ on $A\in\U(N)$ is given by the
Schur polynomial $s_{\lambda^{+}}$ evaluated on the eigenvalues of
$A$. 

We have now defined all the notions from the statement of Theorem
\ref{thm:MP scl for U(N)}. Note that this theorem yields that the
$w$-measures on $\U\left(\bullet\right)$ \emph{determine} the value
of $\scl\left(w\right)$. In particular, if $w_{1}$ and $w_{2}$
induce the same measure on $\U\left(N\right)$ for all $N$, then
$\scl\left(w_{1}\right)=\scl\left(w_{2}\right)$. We do not know if
the parallel result about the ordinary commutator length $\cl(w)$
is true.

\subsection{Stable primitivity rank and word measures on $S_{\bullet}$\label{subsec:Stable-primitivity-rank and S}}

The first result demonstrating a deep connection between a combinatorial\textbackslash topological
invariant of words and word-measures is probably the one dealing with
the primitivity rank, denoted $\pi(w)$, and word measures on the
symmetric groups $S_{\bullet}$. First introduced in \cite{puder2014primitive},
$\pi\left(w\right)$\marginpar{$\pi(w)$} is 
\begin{equation}
\pi\left(w\right)\defi\min\left\{ \rk H\,\middle|\,w\in H\le\F~\mathrm{and}~w~\mathrm{non\textnormal{-}primitive~in}~H\right\} \label{eq:def of pi}
\end{equation}
(an element of a free group is called primitive if it is contained
in some basis).\footnote{The minimum over the empty set is $\infty$. By \cite[Lem.~4.1]{puder2014primitive},
$\pi\left(w\right)=\infty$ if and only if $w$ is primitive in the
ambient group $\F$.} In the notation of the current paper, the result, fully proved in
\cite[Thm.~1.8]{PP15}, states that
\begin{equation}
\pi(w)-1=\beta\left(w,\mathrm{std}\right),\label{eq:pp15, compact version}
\end{equation}
where $\mathrm{std}\colon\s_{N}\to\mathbb{Z}$ is the standard irreducible
character mapping $\sigma\mapsto\#\mathrm{fixed\text{-}points}(\sigma)-1$.
The primitivity rank of $w$ was later shown to be crucial for certain
properties of the one-relator group $\nicefrac{\F}{\langle\langle w\rangle\rangle}$
\cite{louder2022negative,louder2021uniform,linton2022one}. It is
also involved in further bounds on Fourier coefficients of $w$-measures
on $S_{\bullet}$ \cite{hanany2020word}, on $U(\bullet)$ \cite{brodsky2022word}
and on wreath products of the form $G\wr S_{\bullet}$ \cite{shomroni2023wreathI,shomroni2023wreathII}. 

In analogy with the roles of $\cl$ and $\scl$ in $U(\bullet)$,
it is tempting to guess that some stable version of the primitivity
rank is related to all stable irreducible characters of $S_{\bullet}$.
It is futile to naively imitate \eqref{eq:scl, old def} and define
a stable version of the primitivity rank by $\lim_{n\to\infty}\frac{\pi\left(w^{n}\right)}{n}$:
the primitivity rank of every proper power is $1$ (for $n\ge2$,
$w^{n}$ is non-primitive in $\left\langle w\right\rangle $), so
this limit vanishes for every word. A better possible definition is
one where we look for subgroups $H$ of $\F$ containing some power
$w^{n}$ of $w$, such that $w^{n}$ is not only non-primitive in
$H$ but is also the smallest power of $w$ contained in $H$. Namely,
we may consider
\begin{equation}
\inf\left\{ \frac{\rk H-1}{n}\,\middle|\,\begin{gathered}n\in\mathbb{Z}_{\ge1},w^{n}\in H\in\text{\ensuremath{\F}},\\
w^{n}~\mathrm{non\text{-}primitive~in}~H,\\
\forall1\le n'<n:~w^{n'}\notin H
\end{gathered}
\right\} .\label{eq:spi- short but unverified definition}
\end{equation}
We further explain the background to this definition in $\S\S$\ref{subsec:Stable-primitivity-rank}.
This definition is not trivial: it may give different values on words.
For example, for proper powers it gives $0$, but for $w_{g}=\left[a_{1},b_{1}\right]\cdots\left[a_{g},b_{g}\right]$
it gives $2g-1$.

Nonetheless, the ``right'' definition was introduced by Wilton in
\cite[Def.~10.6]{wilton2024rational}. While similar in spirit to
\eqref{eq:spi- short but unverified definition}, it involves not
only subgroups $H$ containing a single power of $w$, but subgroups
$H$ containing conjugates of various powers of $w$, such that $H$
does not freely splits relative to this collection of conjugates --
see Definition \ref{def:spi}. We stick to Wilton's notation and denote
his stable primitivity rank of $w$ by $\sp\left(w\right)$\marginpar{$\protect\sp\left(w\right)$}.
Wilton proves in \cite{wilton2022rationality} that $\sp\left(w\right)$
is rational and computable. We do not know whether Wilton's definition
coincides with \eqref{eq:spi- short but unverified definition}. See
$\S\S$\ref{subsec:Stable-primitivity-rank} for further discussion.

The stable irreducible characters of the symmetric groups $\s_{\bullet}$
are discussed in $\S\S$\ref{subsec:Evidence-towards-Conjecture}.
We believe that their expected values under $w$-measures are closely
related to $\sp\left(w\right)$: 
\begin{conjecture}
\label{conj:spi and stable characters of S_N}Let $\I_{\s}$ denote
the set of all stable characters of $\s_{\bullet}$. Then for every
$w\in\F$ we have 
\begin{equation}
\sp(w)=\inf_{\triv\ne\chi\in{\cal I}_{\s}}\beta\left(w,\chi\right).\label{eq: main conj about spi}
\end{equation}
Moreover, the infimum in the right hand side is attained. 
\end{conjecture}

The right hand side of \eqref{eq: main conj about spi} is trivially
at least $-1$, as $\mathbb{E}_{w}[\chi]\le\dim\chi$ for every finite
dimensional character $\chi$. It follows from \cite[\S\S1.3]{hanany2020word}
that for any $w\ne1$, this side is non-negative, and that it is zero
(and the conjecture holds) for proper powers -- see $\S\S$\ref{subsec:Evidence-towards-Conjecture}.
It is far from obvious that it is bounded away from zero for all non-powers.
However, a beautiful recent result of Cassidy\footnote{Cassidy's result appeared after the first version of the current paper.}
shows that for any $1\ne w\in\F$ which is not a proper power, $\inf_{\triv\ne\chi\in{\cal I}_{\s}}\beta(w,\chi)\ge1$
\cite[Thm.~1.5]{cassidy2024random}. Cassidy used this result to prove
an optimal asymptotic spectral gap for large random Schreier graphs
of the symmetric group. His result also yields Conjecture \ref{conj:spi and stable characters of S_N}
for words in $\F_{2}$ -- see $\S\S$\ref{subsec:Evidence-towards-Conjecture}.

As mentioned above, in the companion paper, we prove a version of
Conjecture \ref{conj:spi and stable characters of S_N} and show that
$\sp(w)$ is equal to the minimum over an infinite set of $\beta(w,\chi)$'s
for stable irreducible characters $\chi$ of the groups $S_{m}\wr S_{\bullet}$
\cite[Thm.~1.2]{PSh25}. This latter result proves that $\sp\left(w\right)$
is determined by the $w$-measures on $S_{m}\wr S_{\bullet}$, and
is thus a profinite invariant, in the following sense.
\begin{defn}
\label{def:profinite invariant}We call a map $f\colon\text{\ensuremath{\F}\ensuremath{\ensuremath{\to}\ensuremath{\mathbb{C}}}}$
\emph{profinite} if $f\left(w_{1}\right)=f\left(w_{2}\right)$ whenever
$w_{1},w_{2}\in\F$ belong to the same $\mathrm{Aut}\hat{\F}$-orbit,
where $\hat{\F}$ is the profinite completion of $\F$ (or, equivalently,
whenever $w_{1}$ and $w_{2}$ induce the same measure on every finite
group --- see \cite[Thm.~2.2]{hanany2020some}.) 
\end{defn}

We remark that it is not obvious that $\sp$ is a new invariant --
see Conjecture \ref{conj:spi=00003Dpi-1}. Conjecture \ref{conj:spi and stable characters of S_N}
is inspired by its analogous results in the cases of $U(\bullet)$
(Theorem \ref{thm:MP scl for U(N)}), of $C_{m}\wr S_{\bullet}$ (Theorem
\ref{thm:result on wreath products and modulo-m spi} below) and of
$S_{m}\wr S_{\bullet}$ \cite[Thm.~1.2]{PSh25}. Further evidence
and a more elaborated version of the conjecture are presented in $\S\S$\ref{subsec:Evidence-towards-Conjecture}.

\subsection{Stable square length and word measures on $\protect\U\left(\bullet\right)$,
$\protect\O\left(\bullet\right)$ and $\protect\Sp\left(\bullet\right)$\label{subsec:Stable-square-length and U O Sp}}

The square length $\mathrm{sql}\left(w\right)$ of $w\in\F$ is the
length of the shortest product of squares in $\F$ yielding $w$.
(It can be defined in arbitrary groups). Rather than the square length
per se, a more natural invariant to stabilize is $\sqlh\left(w\right)\defi\min\left(\mathrm{sql}\left(w\right),2\cdot\cl\left(w\right)\right)$\marginpar{$\protect\sql^{*}$}.
Indeed, as explained in $\S\S$\ref{subsec:Stable-square-length},
while the geometric definition of $\mathrm{sql}$ involves only non-orientable
surfaces, the one of $\sqlh$ involves arbitrary surfaces (orientable
and non-orientable alike) and is algorithmically more natural to work
with, as observed already in \cite{culler1981using}. Moreover, the
invariant $\sqlh$ is the one showing up in the study of word measures:
in \cite[Cor.~1.11]{MPorthsymp}, Magee and the first author show
that the expected trace of a $w$-random real orthogonal matrix in
$\O\left(N\right)$ satisfies
\begin{equation}
\mathbb{E}_{w}\left[\mathrm{tr}\right]=O\left(N^{1-\sqlh\left(w\right)}\right),\label{eq:sql hat bound for O(N)}
\end{equation}
namely, $\beta(w,\mathrm{tr})\ge\sqlh(w)-1$. As in \eqref{eq:bound with cl},
this bound is often, and probably generically, tight. 

The naive ``stable'' definition of $\mathrm{sql}$, $\lim_{n\to\infty}\frac{\mathrm{sql}\left(w^{n}\right)}{n}$,
is useless, as $\mathrm{sql}\left(w^{n}\right)=1$ for $n$ even and
$\mathrm{\mathrm{sql}}\left(w^{n}\right)\le1+\mathrm{sql}\left(w\right)$
for $n$ odd. The naive stable definition of $\sqlh$ is equally pointless.
However, Wilton suggests in \cite{wilton2024rational} a clever way
of stabilizing $\sqlh$ in free groups. In fact, Wilton defines an
invariant $\sigma_{+}\left(X\right)$ for an arbitrary 2-complex $X$
and calls it \emph{the maximal surface curvature }of $X$. When $X$
is the standard presentation complex of the one-relator group $\nicefrac{\F}{\left\langle \left\langle w\right\rangle \right\rangle }$,
the number $1-\sigma_{+}\left(X\right)$ gives rise to an invariant
of $w$ that we call the \emph{stable square length} and denote $\ssql\left(w\right)$.
We give the precise definition in $\S\S$\ref{subsec:Stable-square-length}. 

One can view $\ssql$ as a non-oriented counterpart of $\scl$. Both
invariants are rational in free groups: Calegari proves in \cite{calegari2009stable}
that $\scl\left(w\right)\in\mathbb{Q}$ for every $w\in\F$, and Wilton
proves the analogous result for $\ssql$ in \cite{wilton2022rationality}.
Yet $\ssql$ has an important feature distinguishing it from $\scl$:
while $\scl$ is interesting only inside the commutator subgroup $\left[\F,\F\right]$
and is $\infty$ outside, it follows from \cite{wilton2018essential,wilton2022rationality}
that $\ssql\left(w\right)<\infty$ for every non-primitive word. 

Recall that Theorem \ref{thm:MP scl for U(N)} connects $\scl$ with
all \emph{polynomial} stable representations of $\U\left(\bullet\right)$.
We conjecture that $\ssql\left(w\right)$ plays the same role for
\emph{arbitrary} stable representations of $\U\left(\bullet\right)$.
We believe it plays the same role for all stable representations of
the orthogonal groups $\O\left(\bullet\right)$ and the compact complex
symplectic groups $\Sp\left(\bullet\right)$.
\begin{conjecture}
\label{conj:ssql and its role for U,O and Sp} Let $w\in\F$. Let
${\cal I}_{U}$, ${\cal I}_{O}$ and ${\cal I}_{\Sp}$ denote the
sets of all stable irreducible characters of $\U(\bullet)$, $\O(\bullet)$
and $\Sp(\bullet)$, respectively. Then for ${\cal I}={\cal I}_{\U},{\cal I}_{\O},{\cal I}_{\Sp}$
we have 
\begin{equation}
\ssql(w)=\inf_{\triv\ne\chi\in{\cal I}}\beta(w,\chi).\label{eq:ssql as inf (conj)}
\end{equation}
Moreover, the infimum in the right hand side is attained. 
\end{conjecture}

A more detailed version of this conjecture, distinguishing $\U(\bullet)$
from $\O(\bullet)$ and $\Sp(\bullet)$, is summarized in Remark \ref{rem:same formulas}.
It follows from \cite{MPunitary,MPorthsymp} that for any $w\ne1$,
the right hand side of \eqref{eq:ssql as inf (conj)} is non-negative,
and that it is zero (and the conjecture holds) for proper powers --
see $\S$\ref{subsec:Stable-square-length}. It is far from obvious
that it is bounded away from zero for all non-powers, but if true,
Conjecture \ref{conj:ssql and its role for U,O and Sp}, combined
with Fact \ref{fact:basic properties of ssql}\eqref{enu:ssql on powers and gap on non-powers},
yields that the right hand side is at least $1$ for all non-powers.
We mention a remarkable recent work of Magee and de la Salle where
they show that for any non-power $1\ne w\in\F$, $\inf_{\triv\ne\chi\in{\cal I}_{\U}}\beta(w,\chi)\ge\frac{1}{6}$
\cite[Thm.~3.1(2)]{magee2024strong}.\footnote{In a forth-coming paper, Noam Ta Shma builds on \cite{magee2024strong}
and shows that $\inf_{\triv\ne\chi\in{\cal I}_{\U}}\beta(w,\chi)\ge1$
for non-powers.}

The above-mentioned result of Wilton that $\ssql\left(w\right)$ is
finite for every non-primitive word, is also important in the study
of word measures. A well-known question (e.g., \cite[Ques.~5]{jaikin2023finite})
asks for which families ${\cal C}$ of finite or, more generally,
compact groups, the fact that $w$ induces the Haar measure on every
group in ${\cal C}$ means that $w$ is necessarily primitive. This
is known for ${\cal C}$ being the symmetric groups by \cite{PP15},
and for ${\cal C}$ being $\left\{ \gl_{N}\left(q\right)\right\} _{N}$
with $q$ a fixed prime power, provided that $w\in\F_{2}$ \cite[Cor.~1.7]{West}.
It is still open whether this holds when ${\cal C}$ is the set of
unitary, orthogonal of symplectic groups. Relying on Wilton's result,
establishing Conjecture \ref{conj:ssql and its role for U,O and Sp}
would yield that only primitive words induce the Haar measure on $\U\left(\bullet\right)$
as well as on $\O\left(\bullet\right)$ or on $\Sp\left(\bullet\right)$.

\subsection{Stable mod-$m$ primitivity rank and word measures on $C_{m}\wr S_{\bullet}$\label{subsec:Stable-mod-m-primitivity and wreath}}

For every $m\in\mathbb{Z}_{\ge0}$, denote $C_{m}\defi\left\{ z\in\mathbb{S}^{1}\,\middle|\,z^{m}=1\right\} $,
so $C_{0}=\mathbb{S}^{1}$ is the full group of the unit circle, $C_{1}$
is the trivial group, and for $m\ge2$, $C_{m}$ is cyclic of order
$m$. The wreath product $C_{m}\wr S_{N}$ is the group of $N\times N$
complex matrices with entries in $C_{m}\cup\left\{ 0\right\} $ and
exactly one non-zero in every row and in every column. In particular,
$C_{1}\wr S_{N}=S_{N}$, and $C_{2}\wr S_{N}$ is the hyperoctahedral
group (also known as the signed symmetric group).

The study of word measures on $C_{m}\wr S_{\bullet}=\{C_{m}\wr S_{N}\}_{N\ge1}$
in \cite{MPsurfacewords} and on $G\wr S_{\bullet}$ for more general
groups $G$ in \cite{shomroni2023wreathI,shomroni2023wreathII} led
to interesting generalizations of the primitivity rank and of its
properties. For $1\ne m\in\mathbb{Z}_{\ge0}$, the \emph{mod-$m$
primitivity rank} of $w\in\F$, defined in \cite[Def.~1.10]{MPsurfacewords},
is
\begin{equation}
\pi^{\left(m\right)}\left(w\right)\defi\min\left\{ \rk H\,\middle|\,H\le\F_{r},w\in K_{m}\left(H\right)\right\} ,\label{eq:pi modulo m}
\end{equation}
where $K_{m}\left(H\right)$ is the kernel of the map $H\to(\nicefrac{\mathbb{Z}}{m\mathbb{Z}})^{\rk(H)}$
defined by 
\[
H\twoheadrightarrow\nicefrac{H}{\left[H,H\right]}\cong\mathbb{Z}^{\rk\left(H\right)}\twoheadrightarrow\left(\nicefrac{\mathbb{Z}}{m\mathbb{Z}}\right)^{\rk\left(H\right)}\stackrel{\mathrm{if}~m\ne0}{\cong}C_{m}^{~\rk\left(H\right)}
\]
(the last isomorphism applies only when $m\ne0$, of course). As above,
$\min\left(\emptyset\right)\defi\infty$. If one adds to the definition
\eqref{eq:pi modulo m} the condition that ``$w$ is non-primitive
in $H$'', which is redundant for $m\ne1$, then this definition
works for $m=1$ as well, in which case it coincides with $\pi\left(w\right)$.

Similarly to the $m=1$ case of $S_{N}$ described in \eqref{eq:pp15, compact version},
\cite[Thm.~1.11]{MPsurfacewords} states that for every $1\ne m\in\mathbb{Z}_{\ge0}$,
if $\mathrm{tr}$ is the stable irreducible character corresponding
to the defining representation\footnote{To be sure, the defining representation of $C_{m}\wr S_{N}$ is the
$N$-dimensional representation coming from the embedding $C_{m}\wr S_{N}\hookrightarrow\mathrm{GL}_{N}(\mathbb{C})$.} in $C_{m}\wr S_{\bullet}$, then 
\[
\pi^{\left(m\right)}(w)-1=\beta(w,\mathrm{tr}).
\]

In the current paper we define a stable version of $\pi^{\left(m\right)}$
which generalizes the stable primitivity rank. The \emph{stable mod-$m$
primitivity rank}, $\spm$, is formally defined in Definition \ref{def:spm}.
Morally, it is similar to the definition in \eqref{eq:spi- short but unverified definition}
with ``$w^{n}$ non-primitive in $H$'' replaced by ``$w^{n}\in K_{m}\left(H\right)$''. 

We review the stable representation theory of the groups $C_{m}\wr S_{\bullet}$
in $\S\S$\ref{subsec:stable characters of CmSn}. In short, every
stable irreducible character of $C_{m}\wr S_{\bullet}$ is parameterized
by a partition-valued function $\overrightarrow{\mu}\colon\irr(C_{m})\to{\cal P}$,
where ${\cal P}$ is the set of all partitions of non-negative integers
(with $\arrm$ assigning only finitely many non-empty partitions).
When the only non-empty partition is of the trivial representation,
namely, $\arrm\left(\phi\right)=\emptyset$ for every $\mathrm{triv}\ne\phi\in\irr\left(C_{m}\right)$,
one obtains the corresponding stable representation of $S_{\bullet}$
(indeed, $S_{N}$ is a quotient of $C_{m}\wr S_{N}$, so every irreducible
representation of $S_{N}$ is also one of $C_{m}\wr S_{N}$). So if
Conjecture \ref{conj:spi and stable characters of S_N} holds, then
the infimum of $\beta(w,\chi)$ over all non-trivial stable irreducible
characters $\chi$ of $C_{m}\wr S_{\bullet}$ is at most $\sp(w)$
(and is likely exactly $\sp(w)$). However, there is a natural family
of stable irreducible characters of $C_{m}\wr S_{\bullet}$ on the
other extreme: those that ``do not involve $S_{N}$ at all'', namely,
those where $\arrm\left(\mathrm{triv}\right)=\emptyset$. For those,
we can actually prove that $\beta(w,\chi)\ge\sp(w)$ \cite[Thm.~1.5]{PSh25}.
We focus here on the following, more tight result.
\begin{thm}
\label{thm:result on wreath products and modulo-m spi}Let $1\ne m\in\mathbb{Z}_{\ge0}$,
and let $\phi_{m}$ be a faithful linear character of $C_{m}$ (a
continuous  embedding of $C_{m}$ into $S^{1}\le\mathbb{C^{*}}$).
Denote by ${\cal I}_{m}$ the set of stable irreducible characters
of $C_{m}\wr S_{\bullet}$ corresponding to $\arrm\colon\irr(C_{m})\to{\cal P}$
where the only non-empty partition corresponds to $\phi_{m}$. Then
for every $w\in\F$, we have
\[
\spm(w)=\inf_{\triv\ne\chi\in{\cal I}_{m}}\beta(w,\chi).
\]
Moreover, the infimum in the right hand side is attained.
\end{thm}

The following corollary of the theorem is immediate for $m\ne0$ and
slightly less so for $m=0$.
\begin{cor}
\label{cor:sp(m) profinite}For every $1\ne m\in\mathbb{Z}_{\ge0}$,
$\spm$ is a profinite invariant. 
\end{cor}

In Theorem \ref{thm:spm avoids (0,1)} below, we rely on the random-matrix
type result from Theorem \ref{thm:result on wreath products and modulo-m spi}
in order to prove that $\spm$ admits a gap in $(0,1)$ even for words
that are proper powers.

In \cite[\S4]{PSh25} we define a more general family of stable invariants:
for every compact group $G$ and a non-trivial irreducible character
$\phi\in\irr\left(G\right)$, we introduce the stable invariant $\spf$.
If $G=C_{m}$ and $\phi=\phi_{m}$, this more general definition recovers
$\spm$.

\subsection{Stable $q$-primitivity rank and word measures on $\protect\gl_{\bullet}\left(q\right)$}

Let $q$ be a prime power. The $q$-primitivity rank $\pi_{q}$ of
a word was introduced in \cite[Def.~1.5]{West}. In Appendix \ref{sec:spq and GL_N(q)},
together with Danielle Ernst-West and Matan Seidel, we introduce a
stable version of $\pi_{q}$, denoted $\spq$. More generally, we
define the stable $K$-primitivity rank for an arbitrary field $K$.
It is analogous to the other stable invariants defined above, but
is defined in more algebraic terms and has a different flavor. The
stable $q$-primitivity rank should play the same role for word measures
on $\gl_{\bullet}\left(q\right)$, as $\sp$ is believed to play for
the symmetric groups $S_{\bullet}$. See Appendix \ref{sec:spq and GL_N(q)}
for more details.

\subsection{Common threads\label{subsec:A-unified-point-of-view}}

Table \ref{tab:all stable invariants} sums up the different invariants
we deal with in this paper, and Table \ref{tab:values on selected words}
elaborates their values on some specific words\footnote{Note that Table \ref{tab:values on selected words} compares $2\cdot\scl$,
rather than simply $\scl$, to the other invariants. Indeed, in many
respects, the number $2\cdot\scl(w)$ should have been the stable
commutator length of a word rather than $\scl(w)$ itself: it is much
more common in the major results surrounding $\scl$ (as is the case
in \eqref{eq:MP-scl}) as well as in modern definitions. In addition,
as random-matrix invariants, the other stable invariants we define
are parallel to $2\cdot\scl$.}. We now point to further similarities and relations between them.
Many of these are known for $\scl,\sp,\ssql$ and $\spm$, but not
always for $\spq$. We explain along the paper why all the invariants
are $\mathrm{Aut}\left(\F\right)$-invariant. In addition, many of
the invariants obtain rational values for all words and are computable
-- this was proved in \cite{calegari2009stable} for $\scl$, in
\cite{wilton2022rationality} for $\sp$ and $\ssql$, and in Theorem
\ref{thm:rationality of spm} below for $\spm$.

\begin{table}
\begin{tabular}{|>{\centering}m{4cm}|>{\centering}m{1cm}|>{\centering}m{1.3cm}|>{\centering}m{3cm}|>{\centering}m{2.5cm}|>{\centering}m{3cm}|}
\hline 
Stable invariant &  & Defined in & non-stable version & Related Groups & connection to $w$-measures\tabularnewline
\hline 
\hline 
stable commutator length & $\scl$ & \begin{cellvarwidth}[m]\eqref{eq:scl a la Calegari in F},

\vspace{0.05cm}
\eqref{eq:scl, old def}\end{cellvarwidth} & $\cl$ -- commutator length & \begin{cellvarwidth}[m]$\U\left(\bullet\right)$ 

(poly.~irred.~

characters)\end{cellvarwidth} & \cite[Cor.~1.11]{MPunitary} (stated here as Theorem \ref{thm:MP scl for U(N)})\tabularnewline
\hline 
stable primitivity rank & $\sp$ & Def.~\ref{def:spi} & $\pi$ -- primitivity rank & $S_{\bullet}$ & Conjecture \ref{conj:spi and stable characters of S_N}\tabularnewline
\hline 
stable square length & $\ssql$ & Def.~\ref{def:ssql} & $\sqlh=\min\left(\mathrm{sql,2\cdot\cl}\right)$ & \begin{cellvarwidth}[m]$\U\left(\bullet\right),\O\left(\bullet\right)$,\\\vspace{0.05cm}
$\Sp\left(\bullet\right)$\end{cellvarwidth} & Conjecture \ref{conj:ssql and its role for U,O and Sp}\tabularnewline
\hline 
stable mod-$m$ primitivity rank, $1\ne m\in\mathbb{Z}_{\ge0}$ & $\spm$ & Def.~\ref{def:spm} & $\pi^{\left(m\right)}$ -- mod-$m$ primitivity rank & $C_{m}\wr S_{\bullet}$ & Theorem \ref{thm:result on wreath products and modulo-m spi}\tabularnewline
\hline 
stable $q$-primitivity rank & $\spq$ & Def.~\ref{def:spq} & $\pi_{q}$ -- $q$-primitivity rank & $\gl_{\bullet}\left(q\right)$ & Conjecture \ref{conj:spq and stable characters}\tabularnewline
\hline 
\end{tabular}\caption{\label{tab:all stable invariants}This table summarizes the different
stable invariants discussed in this paper. Each stable invariant is
a stabilized version of a \textquotedblleft degree-1\textquotedblright{}
stable invariant appearing in the third column. The two right-most
columns refer to the relation of the stable invariant, conjectured
or established, with word measures, and more precisely with the expected
values of stable irreducible characters of certain families of groups.}
\end{table}

\begin{table}
\begin{tabular}{|c|c|c|c|c|c|c|c|}
\hline 
$w$ & $\spq\left(w\right)$ & $\sp\left(w\right)$ & $\sp^{\left(2\right)}\left(w\right)$ & $\sp^{\left(m\right)}\left(w\right)$, $m\ge3$ & $\sp^{\left(0\right)}\left(w\right)$ & $\ssql\left(w\right)$ & $2\cdot\scl\left(w\right)$\tabularnewline
\hline 
\hline 
$a$ & \multicolumn{7}{c|}{$\infty$}\tabularnewline
\hline 
$a^{2}$ & $0$ & $0$ & $0$ & $\infty$ & $\infty$ & $0$ & $\infty$\tabularnewline
\hline 
$a^{k}$,$k\ge3$ & $0$ & $0$ & $0$ & $\begin{cases}
0 & m\le k\\
\infty & m>k
\end{cases}$ & $\infty$ & $0$ & $\infty$\tabularnewline
\hline 
$\left[a_{1},b_{1}\right]\cdots\left[a_{g},b_{g}\right]$ & \multicolumn{7}{c|}{$2g-1$}\tabularnewline
\hline 
$a_{1}^{2}\cdots a_{g}^{2}$ & \multicolumn{3}{c|}{$g-1$} & $\infty$ & $\infty$ & $g-1$ & $\infty$\tabularnewline
\hline 
$ab^{-1}cb^{2}ac^{-1}ac$ & ? & $2$ & $\nicefrac{9}{4}$ & $\infty$ & $\infty$ & $\nicefrac{9}{4}$ & $\infty$\tabularnewline
\hline 
\end{tabular}\caption{\label{tab:values on selected words}The values of the stable invariants
on selected words. The $q$ in $\protect\spq$ is any prime power.
As we do not yet know how to compute the values of $\protect\spq\left(w\right)$
for arbitrary words, the value of $\protect\spq$ of the last word
is left as unknown. Some arguments in the computations of the values
in the table are elaborated in Examples \ref{exa:ssql of aabbaB},
\ref{exa:ssql of aBcbbaCac}, \ref{exa:spm} and \ref{exa:spm of aBcbbaCac}.}
\end{table}

\subsubsection*{Inequalities}

Let us summarize a few inequalities between the different stable invariants.
First, we have $\sp\le\ssql\le2\cdot\scl$ (this is part of the content
of Theorem A and Lemma 10.9 in \cite{wilton2024rational}). Proposition
\ref{prop:properties of spm} states that for every $1\ne m\in\mathbb{Z}_{\ge0}$,
$\sp\le\spm\le2\cdot\scl$, that $\sp^{\left(2\right)}\le\ssql$,
and that $\spm\le\sp^{\left(k\right)}$ whenever $m\mid k$.  All
these inequalities, which are summarized in Figure \ref{fig:inequalities between invariants},
are strict inequalities for certain words -- see $\S$\ref{sec:Open-questions}.
In $\S$\ref{sec:Open-questions} and Appendix \ref{sec:spq and GL_N(q)}
we make some comments regarding inequalities involving $\spq$. 
\begin{figure}

\[
\xymatrix{ & 2\cdot\boldsymbol{\scl}\ar@{-}[lddd]\ar@{-}[rd]\\
 &  & \sp^{\left(0\right)}\\
 & \mathbf{}\\
\mathbf{\boldsymbol{\ssql}}\ar@{-}[rd] & \sp^{\left(4\right)}\ar@{-}[d]\ar@{--}[uur] & \sp^{\left(6\right)}\ar@{-}[ld]\ar@{-}[d]\ar@{--}[uu] & \cdots\\
 & \sp^{\left(2\right)}\ar@{-}[rd] & \sp^{\left(3\right)}\ar@{-}[d] & \sp^{\left(5\right)}\ar@{-}[ld]\ar@{--}[luuu] & \cdots\\
 &  & \boldsymbol{\sp}
}
\]
\caption{\label{fig:inequalities between invariants}Known inequalities between
some of the stable invariants discussed in this paper. In every line
in the diagram, the upper invariant is known to be equal to or larger
than the bottom one for every word. A broken line means that there
are intermediate stable invariants bounded between the one in the
bottom and the one in the top.}

\end{figure}

In fact, we have reasons to suspect that $\sp$ is a \emph{universal}
lower bound for all these stable invariants. Combined with our conjectures
regarding the role of these invariants in the theory of word measures,
we suggest the following question:
\begin{question}
\label{que:is sp universal}Let ${\cal I}$ denote the set of \textbf{\emph{all}}
stable irreducible characters of \textbf{\emph{all }}sequences of
compact groups with stable representation theory. Is it true that
for every $w\in\F$
\begin{equation}
\sp(w)=\inf_{\triv\ne\chi\in{\cal I}}\beta(w,\chi)?\label{eq:sp universal?}
\end{equation}
\end{question}

It follows from \cite[Thm.~1.2]{PSh25} that the infimum in the right
hand side of \eqref{eq:sp universal?} is at most $\sp(w)$. An initial
attempt to state a conjecture in this spirit appears in \cite[Conj.~1.13]{hanany2020word},
where Hanany and the first author basically conjecture \eqref{eq:sp universal?}
with $\sp$ replaced by $\pi-1$. This is a combination of Question
$\ref{que:is sp universal}$ with Conjecture \ref{conj:spi=00003Dpi-1}
that $\sp=\pi-1$. While we still think that \cite[Conj.~1.13]{hanany2020word}
is very much plausible (and see Remark \ref{rem:HP conj on S}), the
current paper breaks it into a more accurate and detailed map of conjectures. 

\subsubsection*{Spectral gaps}

Unlike the \emph{non-stable }invariants in Table \ref{tab:all stable invariants}
that obtain only integer values, the \emph{stable }invariants may
obtain non-integral values. However, with the possible exception of
$\spq$, they all share a gap in the interval $\left(0,1\right)$.
Indeed, as mentioned above, $2\cdot\scl\left(w\right)\ge1$ for every
$w\ne1$ (\cite{duncan1991genus}, and see \cite{chen2018spectral}
for a lovely short proof). For $\sp$, proper powers obtain the value
$0$, and for non-powers it is known that $\sp\left(w\right)\ge1$
-- see Theorem \ref{thm:HW-LW}. As $\sp$ bounds from below all
other stable invariants, we get that $\ssql\left(w\right)$ and $\spm\left(w\right)$
are all at least $1$ for non-powers (this also gives another proof
for the gap of $\scl$). As $\ssql$ assigns zero to proper powers,
it, too, admits the $\left(0,1\right)$-gap. Finally, in $\S\S$\ref{subsec:gap of spm}
we use Theorem \ref{thm:result on wreath products and modulo-m spi}
and prove that $\spm$ avoids the interval $\left(0,1\right)$ for
proper powers as well, thus admitting a gap as well.

As $\sp(w)\ge1$ for all non-powers, a positive answer to Question
\ref{que:is sp universal} would yield that $\beta(w,\chi)\ge1$ for
any non-power $w$ and any non-trivial stable irreducible character
$\chi$. We suspect that the gap in $\left(0,1\right)$ is universal
and holds for proper powers at well.
\begin{conjecture}
\label{conj:universal spectral gap }For every $w\in\F$ and every
non-trivial stable irreducible character $\chi$ (of an arbitrary
sequence of groups with stable representation theory), $\beta(w,\chi)\notin\left(0,1\right)$.
\end{conjecture}

\subsubsection*{Do word measures separate $\mathrm{Aut}\protect\F$-orbits?}

The following question naturally arises in the study of word measures
on groups. It is an easy observation that if $w_{2}=\varphi\left(w_{1}\right)$
for some automorphism $\varphi\in\mathrm{Aut}\left(\F\right)$, then
$w_{1}$ and $w_{2}$ induce the same measure on every compact group
(e.g., \cite[Fact 2.5]{magee2015word}). But is the converse true?
The following conjecture was suggested by several mathematicians:
\begin{conjecture}
\label{que:Amit-Vishne for compact}If the words $w_{1}$ and $w_{2}$
induce the same measure on every compact group, do they necessarily
belong to the same $\mathrm{Aut}\F$-orbit?
\end{conjecture}

This conjecture appeared as \cite[Conj.~1.10]{MPunitary}. A stronger
version in which $w_{1}$ and $w_{2}$ are assumed to induce the same
measure only on every \emph{finite} group, appeared as \cite[Ques.~2.2]{amit2011characters}
and \cite[Conj.~4.2]{Shalev2013}. One possible approach to this conjecture
is to collect many invariants of words that are determined by measures
induced on groups, and then prove that two words which agree on all
these invariants must belong to the same $\Aut\F$-orbit. The current
paper suggests many such invariants. Moreover, our full results and
conjectures also suggest that for each stable invariant from Table
\ref{tab:all stable invariants}, not only is the invariant itself
determined by word measures, but also the precise parameters for which
it obtains its extremal value. See, for example, Corollary \ref{cor:same measures =00003D=00003D> same mifkad of extremal scl surfaces},
Proposition \ref{prop:detailed conjectural picture for S} and Remark
\ref{rem:same formulas}.

\subsection{Paper organization}

We start with some preliminaries in $\S$\ref{sec:Preliminaries},
where we explain some notions that are used throughout the paper.
These include core graphs, algebraic morphisms of graphs, Whitehead
graphs and Whitehead's lemma, and fatgraphs. Section \ref{sec:scl+U}
elaborates different equivalent definitions of $\scl$, gives background
to the stable characters of $\U(\bullet)$ and a more elaborated version
of Theorem \ref{thm:MP scl for U(N)}. Section \ref{sec:Stable-primitivity-rank and stable characters of S}
first introduces formally the stable primitivity rank, and then gives
evidence towards Conjecture \ref{conj:spi and stable characters of S_N}.
Then, we define $\ssql$ and explain its relation to word measures
on $\U\left(\bullet\right)$, $\O\left(\bullet\right)$ and $\Sp\left(\bullet\right)$
in $\S$\ref{sec:ssql+U + O + Sp}. Section \ref{sec:spm + G wr S}
deals with wreath products $C_{m}\wr S_{\bullet}$ and their related
stable invariants $\spm$: we prove there Theorem \ref{thm:result on wreath products and modulo-m spi}
and its Corollary \ref{cor:sp(m) profinite}, as well as many other
properties of $\spm$. Then, $\S$\ref{sec:Open-questions} gathers
many of the open questions this work gives rise to. Finally, in Appendix
\ref{sec:spq and GL_N(q)}, together with Danielle Ernst-West and
Matan Seidel, we define $\spq$ and state some of its basic properties
and its conjectural relation to word measures on $\gl_{\bullet}\left(q\right)$.

\subsection*{Notation}

Throughout the paper, $\F$ is a fixed finite-rank free group. The
bouquet $\Omega$ has wedge point $o$ and $\rk\F$ petals, where
each petal is oriented and identified with an element of some fixed
basis of $\F$, so that $\pi_{1}\left(\Omega,o\right)=\F$. The graphs
$\Gamma_{w}$, $\Gamma_{w^{\nu}}$ and $\Gamma_{w^{\sigma}}$ and
the morphisms $\eta_{w},\eta_{w^{\nu}}$ and $\eta_{w^{\sigma}}$
are introduced in Definition \ref{def:Gamma_w} (here $\nu$ is an
integer partition and $\sigma$ a permutation). Given a graph $\Gamma$,
we denote by $V(\Gamma)$ and $E(\Gamma)$ the vertex-set and edge-set
of $\Gamma$, respectively. A Serre graph is defined in $\S\S$\ref{subsec:Whitehead-graphs}.

We denote $\left[d\right]=\left\{ 1,\ldots,d\right\} $. For a finite
set $S$ and $d_{1},\ldots,d_{k}\in\mathbb{Z}_{\ge0}$ with $\sum d_{i}=|S|$,
we denote by $\binom{S}{d_{1}}$ the set of subsets of $S$ of cardinality
$d_{1}$, and by $\binom{S}{d_{1}\,\ldots\,d_{k}}$ the set of all
set-partitions of $S$ into subsets $B_{1}\sqcup\ldots\sqcup B_{k}=S$
with $\left|B_{i}\right|=d_{i}$. A partition of $n\in\mathbb{Z}_{\ge0}$
is $\mu=(\mu_{1},\ldots,\mu_{\ell})$ where $\mu_{1}\ge\mu_{2}\ge\ldots\ge\mu_{\ell}\ge1$
and $\sum\mu_{i}=n$. We denote the fact that $\mu$ is a partition
of $n$ by $\mu\vdash n$, and given $\mu$, we let $\left|\mu\right|$
denote $n=\sum\mu_{i}$ and $\ell(\mu)=\ell$ denote the number of
parts in $\mu$. When $n=0$, we denote by $\emptyset$ the empty
partition. The conjugate of a partition $\mu=(\mu_{1},\ldots,\mu_{\ell})$,
denoted $\mu'$, is obtained from $\mu$ by transposing the corresponding
Young diagram, namely, $\mu'_{1}=\ell$, $\mu'_{2}=\max\left\{ i\,|\,\mu_{i}\ge2\right\} $
and so on. Let ${\cal P}_{n}$ denote the set of all partitions of
$n$, and ${\cal P}=\bigsqcup_{n\in\mathbb{Z}_{\ge0}}{\cal P}_{n}$.

For a compact group $G$ denote by $\irr(G)$ the set of irreducible
complex characters of $G$. We use the word ``irrep'' as an abbreviation
for a complex irreducible representation.

For two functions $f\colon\mathbb{Z}_{\ge N_{0}}\to\mathbb{C}$ and
$g\colon\mathbb{Z}_{\ge N_{0}}\to\mathbb{R}$, we write $f=O(g)$
if there are a constant $c>0$ and $N\in\mathbb{Z}$ such that for
every $n\ge N$, $\left|f(n)\right|\le c\cdot g(n)$. We write $f=\Theta(g)$
if there are constants $c_{1},c_{2}>0$ and $N\in\mathbb{Z}$ such
that for every $n\ge N$, $c_{1}\cdot g(n)\le\left|f(n)\right|\le c_{2}\cdot g(n)$.

\subsection*{Acknowledgments}

We would like to thank Nir Gadish, Daniele Garzoni, Noam Kolodner,
Niv Levhari, Noam Ta Shma and Henry Wilton for beneficial discussions.
Niv Levhari also provided us with the word in Example \ref{exa:ssql of aBcbbaCac}.
We warmly thank Michael Magee for discussions and for his major role
in the several-year-old computer code that allowed us to compute some
expectations of characters in $O\left(N\right)$ -- see Example \ref{exa:measure on O by aabbaB}.
Dror Frid suggested the useful notation of $C_{0}$ for the group
$\mathbb{S}^{1}$, which allowed us to simplify some statements. This
work was supported by the European Research Council (ERC) under the
European Union’s Horizon 2020 research and innovation programme (grant
agreement No 850956), by the Israel Science Foundation, ISF grants
1140/23, as well as by the National Science Foundation under Grant
No. DMS-1926686.

\section{Preliminaries\label{sec:Preliminaries}}

\subsection{Stable representation theory\label{subsec:Stable-representation-theory}}

Stable representation theory of different sequences of groups is a
very natural concept which goes back at least to works of Murnaghan
on the symmetric groups \cite{murnaghan1938analysis}. It was more
systematically studied in recent years, see for example \cite{church2015fi,sam_snowden_2015,gadish2017categories}. 

Unfortunately, to the best of our knowledge, there is no satisfactory
abstract definition of a sequence of groups with a stable representation
theory. Yet, many classical families of groups fall into this category,
and it is usually easy to identify such sequences. In particular,
all sequences $\left\{ G_{N}\right\} _{N}$ of groups with a stable
representation theory that are analyzed in the current paper share
the following features: every stable finite-dimensional representation
has a compact description that applies simultaneously to all but finitely
many groups in the sequence, and its dimension is given by a polynomial
(often in the variable $N$); Each irreducible representation of some
$G_{N}$ is an instance of precisely one stable character of $\left\{ G_{N}\right\} _{N}$;
Every stable character behaves similarly to an element of a ring of
symmetric functions and, in particular, has a well-defined degree;
The product of two stable characters is stable; The expected values
of a stable character under a $w$-measure is a rational function
(often in the variable $N$), and so on.

Recall that we assign a number $\beta(w,\chi)\in\mathbb{Q}\cup\{\infty\}$
to any word $w\in\F$ and any non-trivial stable irreducible character
$\chi=\left\{ \chi_{N}\right\} _{N\ge N_{0}}$. In all the cases analyzed
in the paper, this number can be defined by demanding that $\mathbb{E}_{w}[\chi_{N}]$
is of order $(\dim\chi_{N})^{-\beta(w,\chi)}$ (see Footnote \ref{fn:E=00005Bchi=00005D is of order dim^-beta}
for the precise meaning of this assertion). This defines $\beta(w,\chi)$
unambiguously since a non-trivial stable irreducible character has
dimension growing with $N$. Of course, this definition of $\beta(w,\chi)$
does not require that $\mathbb{E}_{w}[\chi_{N}]$ coincide with a
rational function. More generally, one could define $\beta(w,\chi)$
as 

\begin{equation}
\beta\left(w,\chi\right)\defi\liminf_{N\to\infty}\frac{-\log\left|\mathbb{E}_{w}\left[\chi_{N}\right]\right|}{\log\left(\dim\chi_{N}\right)}.\label{eq:formal-def-of-beta}
\end{equation}

\subsection{Core graphs and $\Gamma_{w^{\nu}}$\label{subsec:Core-graphs}}

The term core graphs is used throughout this paper in the following
sense.
\begin{defn}
\label{def:core graphs}A \textbf{core graph} is a finite simplicial
graph with no leaves and no isolated vertices, namely, with all vertices
of degree $\ge2$. A core graph is not necessarily connected. 
\end{defn}

Core graphs often come equipped with a graph-morphism which is an
immersion to a fixed graph. Sometimes, this fixed graph has fundamental
group identified with $\F$, such as the bouquet \marginpar{$\Omega$}$\Omega=\bigwedge^{r}\mathbb{S}^{1}$
with $r$ oriented petals identified with the elements of a basis
of $\F$. In this case, if the core graph is connected, it is called
a Stallings core graph, after \cite{stallings1983topology}, and represents
a conjugacy class of a non-trivial f.g.~(finitely generated) subgroup
of $\F$. A not-necessarily-connected core graph equipped with an
immersion to $\Omega$ is called a \emph{multi core graph }in\emph{
}\cite{hanany2020word}, and represents a multiset of conjugacy classes
of non-trivial f.g.~subgroups. 

If $H\le\F$ is a f.g.~subgroup, there is a unique Stallings core
graph (with an immersion to $\Omega$) corresponding to the conjugacy
class $H^{\F}$: this is the connected core graph $\Gamma_{H}$ with
an immersion $\eta\colon\Gamma_{H}\to\Omega$ so that the $\eta_{*}\left(\pi_{1}\left(\Gamma_{H}\right)\right)=H^{\F}$.
Throughout the paper, we will need multi core graphs corresponding
to a multiset of powers of $w$:
\begin{defn}
\label{def:Gamma_w}Let $w\in\F$. Denote by \marginpar{$\eta_{w},\Gamma_{w}$}$\eta_{w}\colon\Gamma_{w}\to\Omega$
the Stallings core graph corresponding to $\left\langle w\right\rangle $.
For any partition $\nu=\left(\nu_{1},\ldots,\nu_{\ell}\right)\vdash d$
($d\in\mathbb{Z}_{\ge0}$), denote by \marginpar{$\eta_{w^{\nu}},\Gamma_{w^{\nu}}$}$\eta_{w^{\nu}}\colon\Gamma_{w^{\nu}}\to\Omega$
the multi core graph corresponding to the multiset $\left\{ \left\langle w^{\nu_{1}}\right\rangle ^{\F},\ldots,\left\langle w^{\nu_{\ell}}\right\rangle ^{\F}\right\} $.
If $\nu=\emptyset\vdash0$, then $\Gamma_{w^{\nu}}$ is the empty
graph. If $w=1$, then $\Gamma_{w^{\nu}}$ is the empty graph for
all $\nu$. For $\sigma\in S_{d}$, we denote by \marginpar{$\eta_{w^{\sigma}},\Gamma_{w^{\sigma}}$}$\eta_{w^{\sigma}}\colon\Gamma_{w^{\sigma}}\to\Omega$
the multi core graph corresponding to $w^{\nu}$, where $\nu$ is
the partition of $d$ corresponding to the cycle structure of $\sigma$.
\end{defn}

Topologically, $\Gamma_{w^{\nu}}$ is a disjoint union of $\ell$
cycles. Note that the number of edges in $\Gamma_{w}$ is the length
of the cyclic reduction of $w$. See \cite[Fig.~3.1]{hanany2020word}
for an illustration of $\Gamma_{w^{\left(3,2,1,1\right)}}$.

\subsection{Algebraic and free morphisms of graphs \label{subsec:Algerbaic-morphisms}}

Let $H\le\F$ be free groups. We say that $\F$ is an \emph{algebraic}
\emph{extension} of $H$ if there are no intermediate proper free
factors of $\F$. This concept goes back to Takahasi \cite{takahasi1951note}
and was coined in \cite[Def.~11.1]{KM02}. When $w\ne1$, it is a
simple observation that $\pi\left(w\right)$ from \eqref{eq:def of pi}
is the smallest rank of a proper algebraic extension of $\left\langle w\right\rangle $
in $\F$. 

In \cite[\S4]{hanany2020word} this notion was extended to multiple
subgroups and to morphisms of graphs. We usually say that a free group
$\F$ is\emph{ algebraic over} a multiset of non-trivial subgroups
$\left\{ H_{1},\ldots,H_{t}\right\} $\emph{ }if\emph{ }it is freely
indecomposable relative to $\left\{ H_{1},\ldots,H_{t}\right\} $,
namely, if it cannot be decomposed non-trivially as $\F=J_{1}*J_{2}$
so that $H_{i}$ can be conjugated into $J_{1}$ or $J_{2}$ for each
$i=1,\ldots,t$. The exception is that $\F$ is \emph{not }algebraic
over the empty multiset, even when $\rk\F=1$ (see \cite[Def.~4.6]{hanany2020word}).
Note that this property only depends on the conjugacy classes of $H_{1},\ldots,H_{t}$
in $\F$. We will mostly use this term for graph morphisms, namely,
for simplicial maps between simplicial graphs.
\begin{defn}
\label{def:algebraic morphisms}Let $f\colon\Gamma\to\Delta$ be a
graph morphism. Assume first that $\Delta$ is connected, and let
$\Gamma_{1},\ldots,\Gamma_{t}$ be the connected components of $\Gamma$.
We say that $f$ is \textbf{algebraic}\textbf{\emph{ }}if $\pi_{1}\left(\Delta\right)$
is algebraic over the multiset of conjugacy classes of subgroups $\left\{ f_{*}\left(\Gamma_{1}\right),\ldots,f_{*}\left(\Gamma_{t}\right)\right\} $.
If $\Delta$ is not necessarily connected, we say that $f$ is \textbf{algebraic}
if $f|_{f^{-1}\left(\Delta'\right)}\colon f^{-1}\left(\Delta'\right)\to\Delta'$
is algebraic for every connected component $\Delta'$ of $\Delta$. 

Likewise, if $\Delta$ is connected, we say that $f$ is \textbf{free}
if, up to a proper conjugation, the free product of $f_{*}(\Gamma_{1}),\ldots,f_{*}(\Gamma_{t})$
is a free factor of $\pi_{1}(\Delta)$ (for the precise definition,
see \cite[Def.~4.2]{hanany2020word}). For a general $\Delta$, we
say that $f$ is \textbf{free} if $f|_{f^{-1}\left(\Delta'\right)}\colon f^{-1}\left(\Delta'\right)\to\Delta'$
is free for every connected component $\Delta'$ of $\Delta$.
\end{defn}

We will usually use these notions for graph morphisms that are immersions,
although the definitions apply more generally. Algebraic morphisms
of graphs will appear in the definitions of some of the invariants
we discuss in this paper. They also play an important role in the
formula we present for the expected values of stable characters of
$S_{\bullet}$ (Theorem \ref{thm:formula for stable irreps of S_N}).
An important fact we use below is that every immersion of core graphs
$f\colon\Gamma\to\Delta$ admits a unique decomposition 
\begin{equation}
\Gamma\stackrel{f_{\text{alg}}}{\longrightarrow}\Sigma\stackrel{f_{\text{free}}}{\longrightarrow}\Delta\label{eq:algebraic-free decomposition}
\end{equation}
where $f_{\text{alg}}$ is algebraic, $f_{\text{free}}$ is free and
$f=f_{\text{free}}\circ f_{\text{alg}}$ \cite[Thm.~4.9]{hanany2020word}.

\subsection{Whitehead graphs\label{subsec:Whitehead-graphs}}

Alongside algebraic morphisms of graphs, we also use Whitehead graphs
to describe the stable invariants in this paper. To define Whitehead
graphs, it is useful to think of graphs in Serre's terminology. A
Serre graph $\Gamma$ consists of a vertex set $V\left(\Gamma\right)$
and an edge set $E\left(\Gamma\right)$, which contains two oriented
edges for every geometric edge: every $e\in E\left(\Gamma\right)$
comes with its inverse edge $\overline{e}\in E\left(\Gamma\right)$,
so $e\mapsto\overline{e}$ is a fixed-point-free involution. The map
$\iota=\iota_{\Gamma}\colon E\left(\Gamma\right)\to V\left(\Gamma\right)$
maps every oriented edge to its starting point, so the edge $e$ connects
$\iota\left(e\right)$ with $\iota\left(\overline{e}\right)$. 

The following definition is basically due to Whitehead \cite{Whi36a}.
In this generality it appears, for example, in \cite[\S2]{wilton2018essential}.
\begin{defn}
\label{def:Whitehead graphs}Let $P$ be a finite simplicial graph
which is topologically a disjoint union of cycles, and let $b\colon P\to\Delta$
be a simplicial map from $P$ to a simplicial graph $\Delta$.\footnote{Whitehead graphs can be defined more generally for arbitrary graph
morphisms (and satisfy a corresponding version of Lemma \ref{lem:Whitehead}),
but we will not need this generality here.} For every vertex $v\in V\left(\Delta\right)$ we define its \textbf{Whitehead
graph} $\wh_{b}\left(v\right)$ as follows. The vertices of $\wh_{b}\left(v\right)$
correspond to the edges $\iota_{\Delta}^{~-1}\left(v\right)$. There
are $\left|b^{-1}\left(v\right)\right|$ geometric edges in $\wh_{b}\left(v\right)$:
every $u\in b^{-1}\left(v\right)$ has precisely two incident oriented
edges $\left\{ e_{1},e_{2}\right\} =\iota_{P}^{~-1}\left(u\right)$,
and they give rise to an edge $\left(b\left(e_{1}\right),b\left(e_{2}\right)\right)$
in $\wh_{b}\left(v\right)$. This is illustrated in Figure \ref{fig:Whitehead}.
\end{defn}

\begin{center}
\begin{figure}
\includegraphics[scale=0.3]{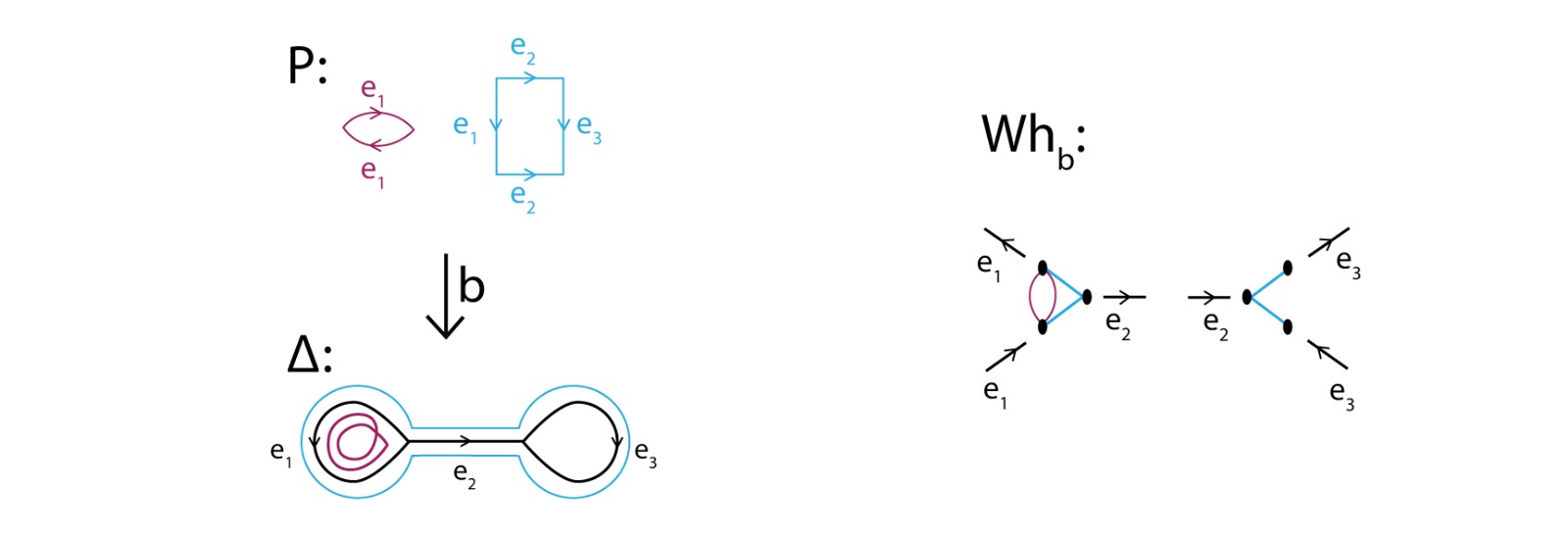}
\centering{}\caption{\label{fig:Whitehead}On the left is a morphism from $P$, a union
of two cycles, to $\Delta$, a barbell-shaped graph. The edges of
$\Delta$ are oriented and labeled in order to depict the morphism:
the image of every edge of $P$ is marked using these labels and orientations.
The two Whitehead graph of this morphism, one for each vertex of $\Delta$,
are drawn on the right. In each of the two Whitehead graphs, each
vertex is labeled by the corresponding half-edge of $\Delta$ emanating
from it.}
\end{figure}
\par\end{center}

It is illuminating to think of Whitehead graphs as in \cite[Fig.~1]{wilton2018essential}.
Whitehead graphs are fundamental tools in detecting primitive elements
in free groups and, more generally, in detecting algebraic morphisms.
\begin{lem}[Whitehead's Lemma]
\label{lem:Whitehead} Let $b\colon P\to\Delta$ be a simplicial
immersion of graphs, where $P$ is topologically a disjoint union
of cycles. If $b$ is \emph{not} algebraic, then for some $v\in V\left(\Delta\right)$,
the Whitehead graph $\wh_{b}\left(v\right)$ is either not connected
or contains a cut-vertex. 
\end{lem}

The original proof in \cite{Whi36a} of a version of Lemma \ref{lem:Whitehead}
uses 3-manifolds. However, there is a very elementary and short argument
due to Wilton \cite[Lem.~2.10]{wilton2018essential} (see also \cite{heusener2019remark}).
We sketch it here because we think it is not well-known enough. Recall
that in a Stallings' fold of a graph (introduced in \cite{stallings1983topology}),
we merge two distinct edges $e_{1}$ and $e_{2}$ emanating from the
same vertex, namely with $\iota\left(e_{1}\right)=\iota\left(e_{2}\right)$.
Such a fold is called a \emph{homotopy-equivalent fold} if it is a
homotopy equivalence, namely, if $\iota\left(\overline{e_{1}}\right)\ne\iota\left(\overline{e_{2}}\right)$.
\begin{proof}[Sketch of proof of Whitehead's Lemma \ref{lem:Whitehead}]
 Assume that $\Delta$ is connected (the general case follows immediately)
and that $b$ is not algebraic. Then there are non-trivial subgroups
$J_{1},J_{2}\le\pi_{1}\left(\Delta\right)$ so that $J_{1}*J_{2}=\pi_{1}\left(\Delta\right)$
and every connected component $P_{j}$ of $P$ satisfies that $b_{*}\left(P_{j}\right)$
can be conjugated into $J_{1}$ or $J_{2}$. For $i=1,2$, let $p_{i}\colon\tilde{\Gamma}_{i}\to\Delta$
be the covering space of $\Delta$ corresponding to (the conjugacy
class of) $J_{i}$, and let $\Gamma_{i}$ be the core of $\tilde{\Gamma}_{i}$.\footnote{Here, the core of $\tilde{\Gamma}_{i}$ is the unique subgraph which
is a core graph and homotopically equivalent to $\tilde{\Gamma}_{i}$.} The assumption about $b_{*}\left(P_{j}\right)$ yields that $b$
factors through an immersion $b'\colon P\to\Gamma_{1}\sqcup\Gamma_{2}$
(so $b=\left(p_{1}\sqcup p_{2}\right)\circ b'$). 

Find a vertex $v_{i}\in V\left(\Gamma_{i}\right)$ so that $p_{1}\left(v_{1}\right)=p_{2}\left(v_{2}\right)$
(if there is no such vertex, add to $\Gamma_{1}$ some additional
edges from $\tilde{\Gamma}_{1}$ in order to find one). Wedge $\left(\Gamma_{1},v_{1}\right)$
with $\left(\Gamma_{2},v_{2}\right)$ to obtain $\left(\Gamma,v\right)$,
with $p\colon\Gamma\to\Delta$ obtained from $p_{1}$ and $p_{2}$.
Stallings' theory of finite graphs \cite{stallings1983topology} yields
that we may obtain $\Delta$ by folding $\Gamma$ along $p$.\footnote{Folding $\Gamma$ along $p$ means that we keep folding pairs of edges
to fix situations where the morphism is not an immersion. This process
is guaranteed to terminate.} As $J_{1}*J_{2}=\pi_{1}\left(\Delta\right)$, we get 
\[
\chi\left(\Gamma\right)=\chi\left(\Gamma_{1}\right)+\chi\left(\Gamma_{2}\right)-1=\left(1-\rk J_{1}\right)+\left(1-\rk J_{2}\right)-1=1-\rk\left(\pi_{1}\left(\Delta\right)\right)=\chi\left(\Delta\right).
\]
As a general folding step cannot decrease the Euler characteristic,
every folding step in the folding process yielding $\Delta$ from
$\Gamma$ must be a homotopy-equivalent fold. Denote by $b''\colon P\to\Gamma$
the immersion so that $b=p\circ b''$. The Whitehead graph $\wh_{b''}\left(v\right)$
is clearly disconnected. If $p\colon\Gamma\to\Delta$ is already folded
(so $p$ is an isomorphism), then $\wh_{b}(p(v))=\wh_{b''}(v)$ and
we are done. Otherwise, consider the last folding step of $\Gamma$
along $p$. Say we folded two edges $e_{1},e_{2}$ with $\iota\left(e_{1}\right)=\iota\left(e_{2}\right)$,
and that $u\in V\left(\Delta\right)$ is the image of $\iota\left(\overline{e_{1}}\right)$
and of $\iota\left(\overline{e_{2}}\right)$, then the image of $\overline{e_{1}}$
in $\Delta$ (which is also the image of $\overline{e_{2}}$) is a
cut-vertex in $\wh_{b}\left(u\right)$. (This folding step can be
pictured as moving from the right side of Figure \ref{fig:unfolding}
to its left side.) 
\end{proof}
\begin{lem}
\label{lem:unconnected Whitehead graph means non-algebraic}Let $b\colon P\to\Delta$
be a simplicial immersion of core graphs, where $P$ is topologically
a disjoint union of cycles. If $\wh_{b}\left(v\right)$ is disconnected
for some $v\in V\left(\Delta\right)$, then $b$ is not algebraic. 
\end{lem}

\begin{proof}
Without loss of generality, assume that $\Delta$ is connected. Let
$E_{1},E_{2}\subseteq\iota_{\Delta}^{-1}\left(v\right)$ be two non-empty,
disjoint complementing subsets with no edges between them in $\wh_{b}\left(v\right)$.
Define a new core graph $\Delta'$ by splitting $v$ into two distinct
vertices $v_{1}$ and $v_{2}$, with $\iota_{\Delta'}\left(e\right)=v_{i}$
for every $e\in E_{i}$. Define $f\colon\Delta'\to\Delta$ naturally
by mapping both $v_{1}$ and $v_{2}$ to $v$. If $\Delta'$ is connected,
then clearly $f_{*}\left(\pi_{1}\left(\Delta',v_{1}\right)\right)$
is a proper free factor (of co-rank one) in $\pi_{1}\left(\Delta,v\right)$.
If $\Delta'$ is disconnected with connected components $\Delta'_{1}$
and $\Delta'_{2}$, then $\pi_{1}\left(\Delta,v\right)=f_{*}\left(\pi_{1}\left(\Delta'_{1},v_{1}\right)\right)*f_{*}\left(\pi_{1}\left(\Delta'_{2},v_{2}\right)\right)$.
In either case, $b$ factors through $b'\colon P\to\Delta'$, which
shows that $b$ is not algebraic.
\end{proof}
Another useful notion developed in \cite{wilton2018essential} is
that of \emph{unfolding}. 
\begin{lem}
\label{lem:unfolding}\cite[Lem.~2.8]{wilton2018essential} Let $b\colon P\to\Delta$
be a simplicial immersion of finite graphs, where $P$ is topologically
a disjoint union of cycles, and assume that $\wh_{b}\left(v\right)$
has a cut vertex for some $v\in V\left(\Delta\right)$. Then we may
\textbf{unfold} $b$ and obtain an immersion $b'\colon\text{P\ensuremath{\to\Delta'}}$\textup{
so that $f\colon\Delta'\to\Delta$ is a single homotopy-equivalent
fold and $b=f\circ b'$. }
\end{lem}

If $e_{0}\in\iota_{\Delta}^{-1}\left(v\right)$ is a cut-vertex, the
unfolding is obtained by splitting $v$ into two distinct vertices
$v_{1}$ and $v_{2}$ and splitting $e_{0}$ (and $\overline{e_{0}}$)
into two distinct edges $e_{1}$ and $e_{2}$ ($\overline{e_{1}}$
and $\overline{e_{2}}$, respectively). The full description appears
in the proof of \cite[Lem.~2.8]{wilton2018essential}, and see also
Figure \ref{fig:unfolding}. The following immediate Corollary also
appears in \cite[\S2]{wilton2018essential}, in a somewhat different
language. Together with Lemma \ref{lem:Whitehead}, it yields an algorithm
to determine if $b\colon P\to\Delta$ is algebraic.
\begin{center}
\begin{figure}
\includegraphics[scale=0.8]{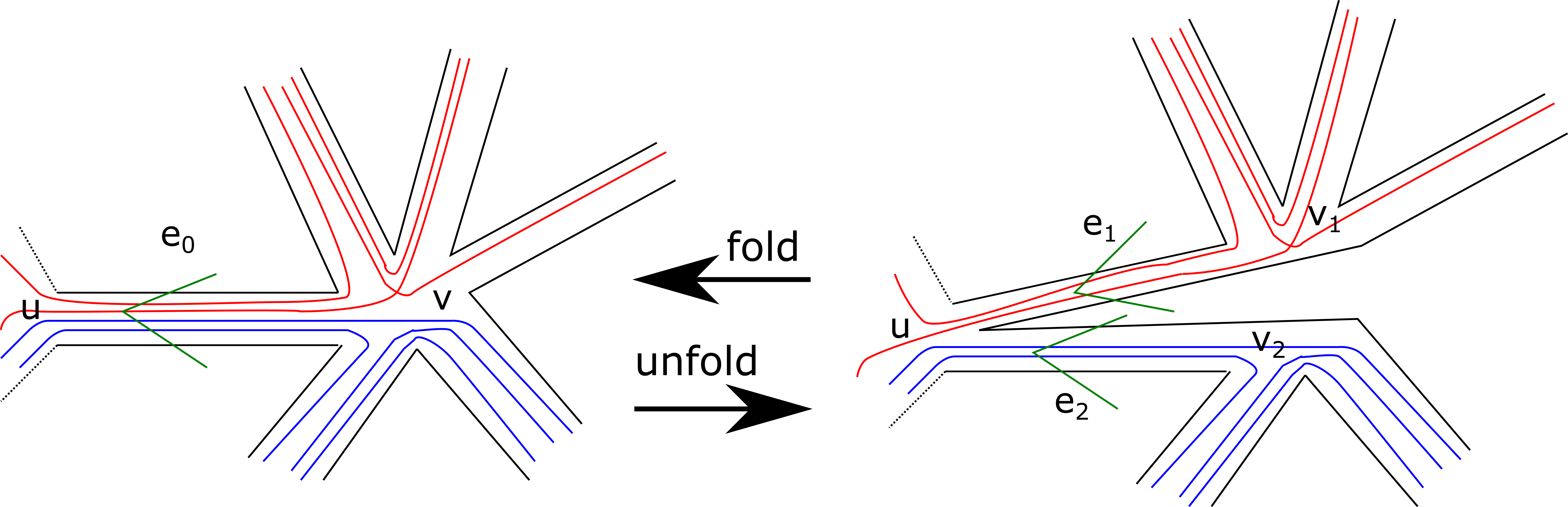}
\centering{}\caption{\label{fig:unfolding}Folding and unfolding: consider an immersion
of finite graphs $b\colon P\to\Delta$, where $P$ is a union of cycles.
Let $v\in V\left(\Delta\right)$ and assume that $\protect\wh_{b}\left(v\right)$
contains a cut vertex $e_{0}\in\iota_{\Delta}^{-1}\left(v\right)$.
Denote $u=\iota_{\Delta}\left(\overline{e_{0}}\right)$. On the left,
we draw the vertex $v$ together with its incident half-edges, with
the edges of $\Delta$ drawn as thick ribbons. The red and blue lines
represent pieces of the $b$-image of $P$ that traverse $v$. Unfolding
at $e_{0}$ results in the picture on the right. Conversely, a homotopy-equivalent
folding of $\overline{e_{1}}$ and $\overline{e_{2}}$ on the right
picture results in the left picture where the merged vertex admits
a cut vertex.}
\end{figure}
\par\end{center}
\begin{cor}
\label{cor:checking algebraicity by unfolding}Let $b\colon P\to\Delta$
be an \textbf{algebraic} simplicial immersion of core graphs, where
$P$ is topologically a disjoint union of cycles. Then after finitely
many unfolding steps, one obtains an immersion $\overline{b}\colon P\to\overline{\Delta}$
where all Whitehead graphs are connected and without cut vertices.
\end{cor}

\begin{proof}
First, every unfolding step is a homotopy equivalence, so the resulting
map of an unfolding step is, too, algebraic. By Lemma \ref{lem:unconnected Whitehead graph means non-algebraic},
$b$, as well as every unfolding of it, has all its Whitehead graphs
connected. As long as there are cut-vertices, we may unfold by Lemma
\ref{lem:unfolding}. It remains to show that this process ends after
finitely many unfolding steps. 

Let $b'\colon P\to\Delta'$ the result of some sequence of unfolds.
By a simple induction, as $\Delta$ is a core graph (all vertices
have degree $\ge2$), so is $\Delta'$. Thus, for every $v\in V\left(\Delta'\right)$,
$\wh_{b'}\left(v\right)$ contains at least two vertices and, being
connected, at least one edge. The total number of edges in $\bigsqcup_{v\in V\left(\Delta'\right)}\wh_{b'}\left(v\right)$
is fixed (and equal to $|E\left(P\right)|/2$), so $\left|V\left(\Delta'\right)\right|\le\left|E\left(P\right)\right|/2$.
Because every unfolding step increases the number of vertices in $\Delta'$
by one, there can be only finitely many such steps. 
\end{proof}

\subsection{Fatgraphs\label{subsec:Fatgraphs}}

Below, when we discuss $\scl$ and $\ssql$, we use \emph{fatgraphs},
possibly with twists. A fatgraph $\mathbb{G}$ (also known as a ribbon
graph) is a graph $G$ together with a cyclic ordering of the edges
incident to each vertex. Every fatgraph $\mathbb{G}$ gives rise to
an orientable surface with boundary by fattening every vertex to a
disc and every edge to a ribbon. We may allow twists -- at most one
twist in every edge -- to allow for non-orientable surfaces.

Let $\Omega$ be a bouquet with $\pi_{1}\left(\Omega,o\right)=\F$
as above. A \emph{fatgraph over $\F$} is a fatgraph $\mathbb{G}$
together with a simplicial map from the underlying graph $G$ to $\Omega$,
namely, so that every edge of $G$ is oriented and labeled by an element
of the fixed basis of $\F$. Twists may also be allowed in fatgraphs
over $\F$. Every fatgraph $\mathbb{G}$ over $\F$ gives rise to
a surface equipped with a map to the bouquet $\Omega$ (the map is
the composition of the retraction of the surface to the underlying
graph with the map from the graph to $\Omega$). For example, Figure
\ref{fig:ssql of aabbaB} depicts a fatgraph over $\F$ with two twisted
edges, where the map to $\Omega$ is described by an orientation and
a label at every edge.

Given a finite fatgraph $\mathbb{G}$ over $\F$ (possibly with twists),
its boundary $\partial\mathbb{G}$ corresponds to a collection of
closed paths in the underlying graph $G$. If $G$ is a core graph
(all vertices have degree $\ge2$), these paths are immersed, i.e.,
non-backtracking, in $G$. Let $P$ be a disjoint union of cycle-graphs
and $b\colon P\to G$ a graph immersion so that $b\left(P\right)$
coincides with $\partial\mathbb{G}$. The Whitehead graphs of $b$
are then all cycles corresponding to the cyclic order of the incident
edges at every vertex. In a sense, the converse is also true. Indeed,
consider diagrams of simplicial maps of finite graphs
\begin{equation}
\begin{gathered}\xymatrix{P\ar@{->>}[dr]^{g~~}\ar[r]^{b} & G\ar[d]^{f}\\
 & \Omega
}
\end{gathered}
\label{eq:diagram of fatgraph}
\end{equation}
where $P$ is topologically a disjoint union of cycles and $\Omega$
is the bouquet.
\begin{lem}
\label{lem:fatgraphs with twists and whitehead} 
\begin{enumerate}
\item A diagram as in \eqref{eq:diagram of fatgraph} corresponds to a fatgraph
over $\F$, possibly with twists, if and only if for every $v\in V\left(G\right)$,
the Whitehead graph $\wh_{b}\left(v\right)$ is a cycle. 
\item Every fatgraph, possibly with twists, can be defined via a diagram
as in \eqref{eq:diagram of fatgraph}.
\item \label{enu:orientability of fatgraphs by WH graphs}A diagram as in
\eqref{eq:diagram of fatgraph} describes an orientable fatgraph (a
fatgraph that can be constructed with no twists), if and only if one
can pick an orientation of the cycles of $P$ and of the cycle in
$\wh_{b}\left(v\right)$ for every $v\in V\left(G\right)$, such that
$b$ becomes orientation preserving.
\end{enumerate}
\end{lem}

To explain Item \ref{enu:orientability of fatgraphs by WH graphs}
in the lemma, let $P_{i}$ be a connected component of $P$. When
we track the cycle $b\left(P_{i}\right)$ in $G$ along the orientation
of $P_{i}$, every visit of a vertex $v$ of $G$, corresponds to
traversing an edge of $\wh_{b}\left(v\right)$. We think of $b$ as
orientation preserving if we always traverse this edge of $\wh_{b}\left(v\right)$
along its own orientation coming from the orientation on the cycle
which is $\wh_{b}\left(v\right)$. We omit the simple proof of Lemma
\ref{lem:fatgraphs with twists and whitehead}.
\begin{rem}
In Lemma \ref{lem:fatgraphs with twists and whitehead}, $\wh_{b}\left(v\right)$
may be a cycle even when $\deg\left(v\right)<3$: if $\deg\left(v\right)=2$,
a cycle is a bigon; if $\deg\left(v\right)=1$, a cycle is a loop
at the single vertex of $\wh_{b}\left(v\right)$, and if $v$ is isolated,
the empty graph is considered a cycle. In practice, though, we only
care about core graphs $G$, so every vertex has degree $\ge2$.
\end{rem}

\section{Stable commutator length and polynomial stable characters of $\protect\U\left(\bullet\right)$\label{sec:scl+U}}

\subsection{Stable commutator length\label{subsec:Stable-commutator-length}}

We recall in some details some basic facts about stable commutator
length, as they inspire many of the definitions of the other invariants
in the following sections. Recall that in any group $G$, the commutator
length of $h\in G$ is 
\begin{equation}
\cl\left(h\right)\defi\min\left\{ g\in\mathbb{Z}_{\ge0}\,\middle|\,\exists a_{1},b_{1},\ldots,a_{g},b_{g}\in G~\text{s.t.}~h=\left[a_{1},b_{1}\right]\cdots\left[a_{g},b_{g}\right]\right\} .\label{eq:cl}
\end{equation}
Let $\Sigma_{g,1}$ be a genus-$g$ orientable surface with one boundary
component and a marked point $v$ at the boundary. It is an easy observation
(and see \cite[\S1]{culler1981using}), that if $\left(X,x_{0}\right)$
is a pointed topological space with $\pi_{1}\left(X,x_{0}\right)\cong G$,
then 
\[
\cl\left(h\right)=\min\left\{ g\in\mathbb{Z}_{\ge0}\,\middle|\,\exists f\colon\left(\Sigma_{g,1},v\right)\to\left(X,x_{0}\right)~\mathrm{s.t.}~f_{*}\left(\left[\partial\Sigma_{g,1}\right]\right)=\left[h\right]\right\} ,
\]
where $\partial\Sigma_{g,1}$ is the loop at $v$ following the boundary
of $\Sigma_{g,1}$. 

\subsubsection*{$\protect\scl$ in arbitrary groups}

Recall from \eqref{eq:scl, old def} that $\scl\left(h\right)=\inf_{n\in\mathbb{Z}_{\ge1}}\frac{\cl\left(h^{n}\right)}{n}$,
so, equivalently,
\begin{equation}
\scl\left(h\right)=\inf\left\{ \frac{g}{n}\,\middle|\,f\colon\left(\Sigma_{g,1},v\right)\to\left(X,x_{0}\right)~\mathrm{s.t.}~f_{*}\left(\left[\partial\Sigma_{g,1}\right]\right)=\left[h^{n}\right]\right\} .\label{eq:scl old geometric definition}
\end{equation}
The equivalent definitions \eqref{eq:scl, old def} and \eqref{eq:scl old geometric definition}
were the ones used for decades (see, e.g., \cite{bavard1991longueur}),
until Calegari noticed that one could also allow an arbitrary number
of boundary components in the definition, and that, in fact, this
flexibility is extremely useful in analyzing $\scl$\footnote{Calegari's definition is also the one that works nicely and naturally
with word measures on $\U\left(\bullet\right)$ -- see $\S\S$\ref{subsec:stable-representation of U}.}. By \cite[Lem.~2.6 and 2.7]{calegari2009stable},
\begin{equation}
2\cdot\scl\left(h\right)=\inf\left\{ \frac{-\chi^{-}\left(\Sigma\right)}{\deg\left(\rho\right)}\,\middle|\,\begin{gathered}\xymatrix{P\ar@{->>}[d]^{\rho~~}\ar[r]^{b} & \Sigma\ar[d]^{f}\\
\left(\mathbb{S}^{1},s_{0}\right)\ar[r]^{h} & \left(X,x_{0}\right)
}
\end{gathered}
\right\} ,\label{eq:scl a la Calegari}
\end{equation}
where $\Sigma$ is a compact oriented surface with boundary in every
connected component, $\chi^{-}\left(\Sigma\right)=\sum_{\Sigma'}\min\left(0,\chi\left(\Sigma'\right)\right)$
the sum going over all the connected components $\Sigma'$ of $\Sigma$,
$\rho$ is a covering map of degree $\deg\left(\rho\right)$, $b$
is an isomorphism between $P$ and the boundary $\partial\Sigma$,
and the diagram commutes. Moreover, the cycle $\mathbb{S}^{1}$ is
oriented according to $h$, this orientation is pulled back to the
cycles in $P$, and then $b$ is orientation preserving where $\partial\Sigma$
gets an induced orientation from an orientation of $\Sigma$. Diagrams
as in \eqref{eq:scl a la Calegari} satisfying all these conditions
are called \emph{admissible monotone surfaces} in \cite{calegari2009stable}.

\subsubsection*{$\protect\scl$ in free groups}

We now focus on $\scl$ in free groups, where much more is known.
Our topological space $X$ is now the bouquet $\Omega$ with $\pi_{1}\left(\Omega,o\right)=\F$
as above. For every $1\ne w\in\F$, there cannot be admissible surfaces
with connected components of non-negative Euler characteristic, so
in \eqref{eq:scl a la Calegari} we may replace $\chi^{-}\left(\Sigma\right)$
with $\chi\left(\Sigma\right)$, and get that for $w\ne1$
\begin{equation}
2\cdot\scl\left(w\right)=\inf\left\{ \frac{-\chi\left(\Sigma\right)}{\deg\left(\rho\right)}\,\middle|\,\begin{gathered}\xymatrix{P\ar@{->>}[d]^{\rho~~}\ar[r]^{b} & \Sigma\ar[d]^{f}\\
\left(\mathbb{S}^{1},s_{0}\right)\ar[r]^{w} & \left(\Omega,o\right)
}
\end{gathered}
\right\} .\label{eq:scl a la Calegari in F}
\end{equation}

Calegari uses \eqref{eq:scl a la Calegari in F} to prove the main
result of \cite{calegari2009stable}: $\scl$ is rational and computable
in free groups. The rationality is obtained by showing the existence,
for every $1\ne w\in\left[\F,\F\right]$, of \emph{extremal surfaces
-- }admissible surfaces where the infimum in \eqref{eq:scl a la Calegari in F}
is attained. Using the LERF property of free groups, Calegari further
shows that these extremal surfaces are $\pi_{1}$-injective {[}Ibid,
Lem.~2.9{]}. The same argument can be used to show that these surfaces
are, moreover, ``efficient'', in the following sense. 
\begin{defn}
\label{def:efficient map with surfaces}A diagram as in \eqref{eq:scl a la Calegari in F}
is called \textbf{efficient} if it satisfies the following two properties:
\begin{enumerate}
\item $f$ is $\pi_{1}$-injective, and
\item for every two distinct points $s_{1},s_{2}\in\rho^{-1}\left(s_{0}\right)\subseteq P$
and every path $q$ in $\Sigma$ from $b\left(s_{1}\right)$ to $b\left(s_{2}\right)$,
the image of $q$ in $\Omega$, which is a loop at $o$, is \emph{not}
null-homotopic. In particular, $b$ is injective on $\rho^{-1}\left(s_{0}\right)$. 
\end{enumerate}
\end{defn}

If the diagram is not efficient, the argument in the proof of {[}Ibid.{]}
shows that after possibly passing to a cover of $\Sigma$ and cutting
the cover along an embedded path or loop, we may obtain a diagram
with a strictly smaller ratio $\frac{-\chi\left(\Sigma'\right)}{\deg\left(\rho'\right)}$.
So the efficiency of the diagram is a necessary condition for it to
be extremal. Thus, \cite{calegari2009stable} yields the following
corollary:
\begin{cor}
\label{cor:scl in free groups - restricted definition}Let $1\ne w\in\F$.
Then 
\begin{equation}
2\cdot\scl\left(w\right)=\inf\left\{ \frac{-\chi\left(\Sigma\right)}{\deg\left(\rho\right)}\,\middle|\,\begin{gathered}\xymatrix{P\ar@{->>}[d]^{\rho~~}\ar[r]^{b} & \Sigma\ar[d]^{f}\\
\left(\mathbb{S}^{1},s_{0}\right)\ar[r]^{w} & \left(\Omega,o\right)
}
\end{gathered}
~\mathrm{s.t.~the~diagram~is~efficient}\right\} ,\label{eq:scl in free groups restricting definition}
\end{equation}
where, as above, the diagram represents an admissible monotone surface:
$\Sigma$ is a compact oriented surface with boundary in every connected
component, $\rho$ is a covering map of degree $\deg\left(\rho\right)$,
$b$ is an orientation-preserving isomorphism between $P$ and the
boundary $\partial\Sigma$, and the diagram commutes. 
\end{cor}

In the case of $\scl$, the restriction to efficient maps is not only
unnecessary, but in a sense also counter-productive. Indeed, Calegari's
proof of the rationality (and efficient computability -- see \cite[\S\S\S 4.1.7-4.1.8]{calegari2009scl})
of $\scl$ in free groups heavily relies on the flexibility given
by the less restrictive definition \eqref{eq:scl a la Calegari in F}:
in the proof, he constructs admissible surfaces out of small pieces,
and these surfaces are not guaranteed to be efficient -- only the
extremal ones are. 

However, the more restrictive definition \eqref{eq:scl in free groups restricting definition}
is crucial if one cares to define other, similar invariants in the
same spirit, and to compare them with $\scl$. This is especially
true for defining the non-oriented analog of $\scl$ that we call
stable square length. See $\S\S$\ref{subsec:Stable-square-length}.\\

Finally, in the set of admissible surfaces from Corollary \ref{cor:scl in free groups - restricted definition}
we may further restrict to orientable fatgraphs (namely, fatgraphs
without twists), as defined in $\S\S$\ref{subsec:Fatgraphs}. A surface
$\Sigma$ with a map $f$ to the bouquet $\Omega$ is called incompressible
if no non-nullhomotopic embedded loop is mapped to a null-homotopic
loop in $\Omega$. Every incompressible admissible surface $f\colon\Sigma\to\Omega$
can be constructed, up to homotopy, from a fatgraph over $\F$, possibly
with twists. This is essentially due to Culler \cite{culler1981using},
and see \cite[Prop.~2.7]{calegari2015surface}. By the aforementioned
\cite[Lem.~2.9]{calegari2009stable}, every extremal surface is incompressible,
(this also follows from $f$ being $\pi_{1}$-injective) so we may
further restrict the surfaces in \eqref{eq:scl in free groups restricting definition}
to surfaces defined by fatgraphs. By Lemma \ref{lem:fatgraphs with twists and whitehead}
we obtain the following (recall Definition \ref{def:Gamma_w} of $\Gamma_{w}$).
\begin{cor}
\label{cor:scl in free groups - def with WH graphs}Let $1\ne w\in\F$.
Then 
\begin{equation}
2\cdot\scl\left(w\right)=\inf\left\{ \frac{-\chi\left(\Gamma\right)}{\deg\left(\rho\right)}\,\middle|\,\begin{gathered}\xymatrix{P\ar@{->>}[d]^{\rho~~}\ar[r]^{b} & \Gamma\ar[d]^{f}\\
\Gamma_{w}\ar[r]^{\eta_{w}} & \Omega
}
\end{gathered}
\begin{gathered}\mathrm{s.t.~\rho~\mathrm{is~a~covering~map,}}\\
b~\mathrm{is~an~immersion~of~core~graphs,}\\
\mathrm{the~diagram~is~efficient,}\\
\wh_{b}\left(v\right)~\mathrm{is~a~cycle}~\forall v\in V\left(\Gamma\right),\\
\mathrm{and}~b~\mathrm{is~orientation~preserving}
\end{gathered}
\right\} ,\label{eq:scl in free groups with WH}
\end{equation}
($b$ is orientation preserving in the sense of Lemma \ref{lem:fatgraphs with twists and whitehead}\eqref{enu:orientability of fatgraphs by WH graphs}).
As in Definition \ref{def:efficient map with surfaces}, \emph{efficient}
here means that $f$ is $\pi_{1}$-injective and that for every $v\in V\left(\Gamma_{w}\right)$,
every two distinct points $s_{1},s_{2}\in\rho^{-1}\left(v\right)$
and every path $q$ from $b\left(s_{1}\right)$ to $b\left(s_{2}\right)$,
the image of $q$ in $\Omega$ is not null-homotopic.
\end{cor}

\subsection{Stable characters of $\protect\U\left(\bullet\right)$\label{subsec:stable-representation of U}}

Complex irreducible representations (irreps for short) of the unitary
group $\U\left(N\right)$ are parameterized by non-increasing integer
sequences $a_{1}\ge a_{2}\ge\ldots\ge a_{N}$ of length $N$.\footnote{For this and other basic facts about the irreducible representations
of $U(N)$ see, e.g., \cite[\S36]{bump2013lie}.} A stable irreducible character of $U(\bullet)$ is given by fixing
such a sequence of length $N_{0}$ and padding it with additional
$N-N_{0}$ zeros (in suitable positions so that it remains non-increasing),
to get a non-increasing integer sequence of length $N$, for every
$N\ge N_{0}$ (e.g., \cite[\S3]{sam_snowden_2015}). We formalize
it by considering a pair of partitions $\lambda=(\lambda^{+},\lambda^{-})$,
with $\lambda^{+}=(\lambda_{1}^{+},\ldots,\lambda_{\ell^{+}}^{+})$
and $\lambda^{-}=(\lambda_{1}^{-},\ldots,\lambda_{\ell^{-}}^{-})$,
and for every $N\ge\ell^{+}+\ell^{-}$ letting 
\[
\lambda[N]\defi\left(\lambda_{1}^{+},\ldots,\lambda_{\ell^{+}}^{+},\underbrace{0,\ldots,0}_{N-\ell^{+}-\ell^{-}},-\lambda_{\ell^{-}}^{-},\ldots,-\lambda_{1}^{-}\right)
\]
be a non-increasing sequence of length $N$. We denote the corresponding
character of $\U(N)$ by $\xi^{\lambda[N]}$\marginpar{$\xi^{\lambda[N]}$},
and the corresponding stable character by \marginpar{$\xi^{\lambda[\bullet]}$}$\xi^{\lambda[\bullet]}=\{\xi^{\lambda\left[N\right]}\}_{N\ge\ell^{+}+\ell^{-}}$.
We also let $\left|\lambda\right|\defi\left|\lambda^{+}\right|+\left|\lambda^{-}\right|$
and $\ell(\lambda)\defi\ell(\lambda^{+})+\ell(\lambda^{-})=\ell^{+}+\ell^{-}$.
If $\lambda^{-}$ is the empty partition $\emptyset$, the corresponding
irreps are called polynomial, and for every $A\in\U(N)$, $\xi^{\lambda[N]}(A)$
is equal to the Schur polynomial $s_{\lambda^{+}}$ evaluated on the
eigenvalues of $A$.\footnote{In contrast, it follows from \eqref{eq:Koike} that a non-polynomial
character is given by a polynomial in $\alpha_{1},\ldots,\alpha_{N},\overline{\alpha_{1}},\ldots,\overline{\alpha_{N}}$,
where $\alpha_{1},\ldots,\alpha_{N}$ are the eigenvalues of $A$.
This polynomial is symmetric on $\alpha_{1},\ldots,\alpha_{N}$, and
separately on $\overline{\alpha_{1}},\ldots,\overline{\alpha_{N}}$.} In this case, we abuse notation and denote by $\xi^{\lambda^{+}[\bullet]}=\left\{ \xi^{\lambda^{+}\left[N\right]}\right\} _{N\ge\ell^{+}}$
the stable character corresponding to $(\lambda^{+},\emptyset)$.

For example, if $\lambda=\left((1),\emptyset\right)$, then $\lambda\left[N\right]=\left(1,0,\ldots,0\right)$
for all $N\ge1$, and the corresponding irrep is the \emph{standard
}representation $\U\left(N\right)\hookrightarrow\gl_{N}\left(\mathbb{C}\right)$,
which is polynomial. We have $\xi^{\lambda\left[N\right]}\left(A\right)=\mathrm{tr}_{N}\left(A\right)$
for every $A\in\U\left(N\right)$. As another example, $\lambda=\left(\lambda^{+},\lambda^{-}\right)$
with $\lambda^{+}=\lambda^{-}=(1)\vdash1$ gives rise to the stable
non-polynomial character $\xi^{\lambda[\bullet]}=\{\xi^{\lambda\left[N\right]}\}_{N\ge2}$
with $\lambda\left[N\right]=\left(1,0,\ldots0,-1\right)$. In this
case, the corresponding irrep of $\U\left(N\right)$ is $\left(N^{2}-1\right)$-dimensional
and $\xi^{\lambda[N]}(A)=\mathrm{tr}(A)\cdot\mathrm{tr}(\overline{A})-1$. 

For every $k\in\mathbb{Z}$, define $\zeta_{k}\colon\U\left(N\right)\to\mathbb{C}$
by $A\mapsto\mathrm{tr}(A^{k})$. For every (non-increasing, say)
integer sequence $\mu=\left(\mu_{1},\ldots,\mu_{\ell}\right)$, define
$\zeta_{\mu}\defi\zeta_{\mu_{1}}\cdots\zeta_{\mu_{\ell}}$. Note that
$\zeta_{\mu}$ is defined for all $N$. We have the following formula
(e.g., \cite[Cor.~7.17.5]{Stanley1999enumerative2}), connecting these
trace maps to stable polynomial irreducible characters. Let $\lambda=(\nu,\emptyset)$.
Then for every $N\ge\ell(\nu)$,
\begin{equation}
\xi^{\lambda\left[N\right]}=\xi^{\nu\left[N\right]}=\sum_{\mu\vdash\left|\nu\right|}p_{\mu}\cdot\chi^{\nu}\left(\mu\right)\cdot\zeta_{\mu}=\frac{1}{\left|\nu\right|!}\sum_{\sigma\in S_{\left|\nu\right|}}\chi^{\nu}\left(\sigma\right)\cdot\zeta_{\sigma},\label{eq:formula for stable poly irreps of U}
\end{equation}
where $p_{\mu}$ is the probability that a random permutation in $S_{\left|\mu\right|}$
has cycle structure $\mu$, where $\chi^{\nu}\left(\mu\right)$ is
the value of the character of $S_{\left|\nu\right|}$ corresponding
to $\nu$ on a permutation of cycle structure $\mu$, and where $\zeta_{\sigma}\defi\zeta_{\mu}$
if $\mu$ is the cycle structure of $\sigma$. For example, if $\nu=(1)$
then $\xi^{\nu\left[N\right]}=\zeta_{1}$ and if $\nu=(2,1)$ then
$\xi^{\nu\left[N\right]}=\frac{1}{3}\zeta_{\left(1,1,1\right)}-\frac{1}{3}\zeta_{\left(3\right)}$.
In particular, it follows from \eqref{eq:formula for stable poly irreps of U}
that 
\begin{equation}
\dim\left(\xi^{\nu\left[N\right]}\right)=\xi^{\nu\left[N\right]}\left(1\right)=N^{\left|\nu\right|}\left(\frac{\dim\left(\chi^{\nu}\right)}{\left|\nu\right|!}+O\left(\frac{1}{N}\right)\right)=\Theta\left(N^{\left|\nu\right|}\right).\label{eq:dim of stable poly irrep of U}
\end{equation}
In \eqref{eq:dim of arbitrary stable irrep of U} we give the same
estimate for arbitrary stable characters of $U(\bullet)$, where $\lambda^{-}$
is not necessarily empty. Indeed, Koike \cite[Eq.~(0.3)]{koike1989decomposition}
gives a formula for the value of $\xi^{\lambda\left[N\right]}$ for
every $\lambda=(\lambda^{+},\lambda^{-})$ in terms of polynomial
characters. For every $N\ge\ell(\lambda)$ and $A\in\U\left(N\right)$,
\begin{equation}
\xi^{\lambda\left[N\right]}\left(A\right)=\sum_{\tau,\nu,\omega}\left(-1\right)^{\left|\tau\right|}\mathrm{c}_{\tau,\nu}^{\lambda^{+}}c_{\tau',\omega}^{\lambda^{-}}\xi^{\nu\left[N\right]}\left(A\right)\xi^{\omega\left[N\right]}\left(A^{-1}\right),\label{eq:Koike}
\end{equation}
where the sum is over partitions $\tau,\nu,\omega$, the numbers $c_{\tau,\nu}^{\lambda^{+}}$
and $c_{\tau^{t},\omega}^{\lambda^{-}}$ are the Littlewood-Richardson
coefficients (e.g., \cite[\S I.9]{macdonald1998symmetric}), and $\tau'$
is the conjugate of $\tau$. Note that the summation is finite as,
for example, we must have $\left|\tau\right|+\left|\nu\right|=\left|\lambda^{+}\right|$
and $\left|\tau\right|+\left|\omega\right|=\left|\lambda^{-}\right|$
for the Littlewood-Richardson coefficients to not vanish. Relying
on \eqref{eq:dim of stable poly irrep of U}, one can observe that
in the expression obtained for $\dim(\xi^{\lambda[N]})=\xi^{\lambda[N]}(1)$
from \eqref{eq:Koike}, the summands of highest order are those where
$\tau=\emptyset$. In this case, the Littlewood-Richardson coefficients
vanish unless $\nu=\lambda^{+}$ and $\omega=\lambda^{-}$. As $c_{\emptyset,\nu}^{\nu}=1$,
we get

\begin{equation}
\dim\left(\xi^{\lambda\left[N\right]}\right)=\dim\left(\xi^{\lambda^{+}[N]}\right)\dim\left(\xi^{\lambda^{-}[N]}\right)\left(1+O\left(\frac{1}{N}\right)\right)=\Theta\left(N^{\left|\lambda\right|}\right).\label{eq:dim of arbitrary stable irrep of U}
\end{equation}
Moreover, $\dim(\xi^{\lambda[N]})$ is given by a polynomial in $N$
of degree $\left|\lambda\right|$. 

Combining \eqref{eq:formula for stable poly irreps of U} and \eqref{eq:Koike},
we get a formula for $\xi^{\lambda\left[N\right]}$ in terms of the
$\zeta_{\mu}$'s. For example, for $\lambda^{+}=(2)$ and $\lambda^{-}=(1,1)$
we get 
\begin{eqnarray*}
\xi^{\lambda\left[N\right]}\left(A\right) & = & \xi^{\left(2\right)\left[N\right]}\left(A\right)\xi^{\left(1,1\right)\left[N\right]}\left(A^{-1}\right)-\xi^{\left(1\right)\left[N\right]}\left(A\right)\xi^{\left(1\right)\left[N\right]}\left(A^{-1}\right)+\xi^{\left(\right)\left[N\right]}\left(A\right)\xi^{\left(\right)\left[N\right]}\left(A^{-1}\right)\\
 & = & \frac{\zeta_{\left(1,1\right)}\left(A\right)+\zeta_{\left(2\right)}\left(A\right)}{2}\cdot\frac{\zeta_{\left(1,1\right)}\left(A^{-1}\right)-\zeta_{\left(2\right)}\left(A^{-1}\right)}{2}-\zeta_{\left(1\right)}\left(A\right)\zeta_{\left(1\right)}\left(A^{-1}\right)+1\cdot1,
\end{eqnarray*}
hence
\begin{equation}
\xi^{\lambda\left[N\right]}=\frac{\zeta_{\left(1,1,-1,-1\right)}+\zeta_{\left(2,-1,-1\right)}-\zeta_{\left(1,1,-2\right)}-\zeta_{\left(2,-2\right)}}{4}-\zeta_{\left(1,-1\right)}+1.\label{eq:xi^lambda as a sum over zeta_mu - example}
\end{equation}
For every integer sequence $\mu$, the expectations of $\mathbb{E}_{w}\left[\zeta_{\mu}\right]$
in $\left\{ \U\left(N\right)\right\} _{N}$ coincide with a rational
function in $\mathbb{Q}\left(N\right)$ for all sufficiently large
$N$ \cite[Prop.~1.1]{MPunitary}. If $w\ne1$ and $\mu$ contains
no zeros, then, moreover, the degree of this rational function is
non-positive\footnote{Here, the degree of $f\in\mathbb{Q}(N)$ is the degree of the numerator
minus the degree of the denominator, or $-\infty$ if $f=0$.} (e.g., {[}Ibid., Cor.~1.13{]}). We thus have the same result for
$\mathbb{E}_{w}\left[\xi^{\lambda\left[N\right]}\right]$: it coincides
with a rational function, which we denote by $f(w,\lambda)$, for
all sufficiently large $N$, and if $w\ne1$, this function has non-positive
degree. For example, if $\lambda=((1),\emptyset)$ and $w=\left[x,y\right]^{2}$,
we have $f(w,\lambda)=\mathbb{E}_{\left[x,y\right]^{2}}\left[\xi^{\left(1\right)\left[N\right]}\right]=\mathbb{E}_{\left[x,y\right]^{2}}\left[\mathrm{tr}\right]=\frac{-4}{N^{3}-N}$
for all $N\ge2$. If $\left|\lambda\right|\ne0$, then $f(w,\lambda)$
is of order 
\[
N^{\deg\left(f(w,\lambda)\right)}=\left(\dim\xi^{\lambda[N]}\right)^{\deg\left(f(w,\lambda)\right)/\left|\lambda\right|}.
\]
Recall that, by definition, $\mathbb{E}_{w}[\xi^{\lambda[N]}]$ is
of order $(\dim\xi^{\lambda[N]})^{-\beta(w,\xi^{\lambda[\bullet]})}$,
hence
\[
\beta\left(w,\xi^{\lambda[\bullet]}\right)=\frac{-\deg\left(f(w,\lambda)\right)}{\left|\lambda\right|}\in\mathbb{Q}_{\ge0}\cup\left\{ \infty\right\} .
\]
For example, $\beta([x,y]^{2},\xi^{(1)[\bullet]})=3$. 

\subsection{Polynomial stable characters of $\protect\U\left(\bullet\right)$
and Theorem \ref{thm:MP scl for U(N)}\label{subsec:Polynomial-stable-representation of U}}

We now focus on polynomial stable characters $\xi^{\nu[\bullet]}$
where $\nu$ is a partition and explain how Theorem \ref{thm:MP scl for U(N)}
follows from \cite{MPunitary}. The main result of {[}Ibid.{]} is
a formula for the expected value of the $\zeta_{\mu}$-s under word
measures. By {[}Ibid., Thm.~1.7{]}, $\mathbb{E}_{w}\left[\zeta_{\mu}\right]$
is given in terms of maps from orientable surfaces to the bouquet
$\Omega$, where each surface has $\ell$ boundary components mapping
to loops corresponding to $w^{\mu_{1}},\ldots,w^{\mu_{\ell}}$ in
the bouquet. More precisely, 
\begin{equation}
\mathbb{E}_{w}\left[\zeta_{\mu}\right]=\sum_{\left[\Sigma,f\right]}c_{\Sigma,f}\cdot N^{\chi\left(\Sigma\right)},\label{eq:MP19}
\end{equation}
where the sum is over all admissible surfaces with boundary as above
(up to some natural equivalence relation -- see {[}Ibid., Def.~1.3{]}),
and $c_{\Sigma,f}$ is an integer which is the $L^{2}$-Euler characteristic
of some stabilizer of $\left[\Sigma,f\right]$. 

When $\mu$ is a partition (namely, when all the $\mu_{i}$'s are
positive), these maps and surfaces are precisely the admissible monotone
surfaces from \eqref{eq:scl a la Calegari in F}. Together with \eqref{eq:formula for stable poly irreps of U},
we get a formula for $\mathbb{E}_{w}\left[\xi^{\nu\left[N\right]}\right]$
involving admissible monotone surfaces for $w$, all of which have
degree $\left|\nu\right|$. Together with \eqref{eq:dim of stable poly irrep of U},
this implies that 

\[
\mathbb{E}_{w}\left[\xi^{\nu\left[N\right]}\right]=O\left(N^{-2\scl\left(w\right)\cdot\left|\nu\right|}\right)\stackrel{\eqref{eq:dim of stable poly irrep of U}}{=}O\left(\left(\dim\xi^{\nu\left[N\right]}\right)^{-2\scl\left(w\right)}\right),
\]
namely, $\beta(w,\xi^{\nu[\bullet]})\ge2\scl(w)$.

On the other hand, let $1\ne w\in\left[\F,\F\right]$. Calegari's
result guarantees the existence of extremal surfaces for $w$. Say
that $\left[\Sigma,f\right]$ is such an extremal surface with boundary
components mapping to $w^{\mu_{1}},\ldots,w^{\mu_{\ell}}$, with $\mu=\left(\mu_{1},\ldots,\mu_{\ell}\right)$
a partition, and fix $\nu=\left(\left|\mu\right|\right)\vdash\left|\mu\right|$.
Then $\chi^{\nu}$ is the trivial character of $S_{\left|\mu\right|}$,
and in this case \eqref{eq:formula for stable poly irreps of U} becomes
\[
\mathbb{E}_{w}\left[\xi^{\left(\left|\mu\right|\right)\left[N\right]}\right]=\sum_{\omega\vdash\left|\mu\right|}p_{\omega}\mathbb{E}_{w}\left[\zeta_{\omega}\right],
\]
where all the coefficients $p_{\omega}$ are positive. We already
explained above why $\mathbb{E}_{w}\left[\zeta_{\omega}\right]$ is
in\linebreak{}
$O(\left(\dim\xi^{\left(\left|\mu\right|\right)\left[N\right]}\right)^{-2\scl\left(w\right)})$.
Because $\left[\Sigma,f\right]$ is extremal, we have $\chi\left(\Sigma\right)=-2\scl\left(w\right)\cdot\left|\mu\right|$
by the definition of an extremal surface. Moreover, every admissible
surface of degree $\left|\mu\right|$ and Euler characteristic $-2\scl\left(w\right)\cdot\left|\mu\right|$
is extremal. Finally, because extremal surfaces are $\pi_{1}$-injective,
it follows that their coefficient $c_{\Sigma,f}$ is equal to 1 {[}Ibid.,
Lem.~5.1{]}. We obtain that $\mathbb{E}_{w}\left[\zeta_{\mu}\right]=\Theta\left(N^{-2\scl\left(w\right)\cdot\left|\mu\right|}\right)$,
so
\[
\beta\left(w,\xi^{\left(\left|\mu\right|\right)[\bullet]}\right)=2\scl(w).
\]
This explains the statement of Theorem \ref{thm:MP scl for U(N)}.

Let us point out one further corollary from the above discussion,
which is relevant for separating $\mathrm{Aut}\F$-orbits by measures
on groups.
\begin{cor}
\label{cor:same measures =00003D=00003D> same mifkad of extremal scl surfaces}If
$w_{1}$ and $w_{2}$ induce the same measure on $\U\left(N\right)$
for all $N$, then not only do we have the equality $\scl\left(w_{1}\right)=\scl\left(w_{2}\right)$,
but both words also have the precise same ``census'' of extremal
surfaces. Namely, for every partition $\mu=\left(\mu_{1},\ldots,\mu_{\ell}\right)$,
the number of non-equivalent extremal surfaces\footnote{Here, extremal surfaces are equivalent in the sense of \cite[Def.~1.3]{MPunitary}.}
of $w_{1}$ with boundary components corresponding to the powers $\mu_{1},\ldots,\mu_{\ell}$,
is identical to the number of such surfaces of $w_{2}$.
\end{cor}

\begin{proof}
The discussion above shows that 
\begin{equation}
\mathbb{E}_{w}\left[\zeta_{\mu}\right]=N^{-2\scl\left(w\right)\cdot\left|\mu\right|}\left(d_{w,\mu}+O\left(N^{-1}\right)\right),\label{eq:leading term of E_w=00005BXi=00005D in U}
\end{equation}
where $d_{w,\mu}$ is the number of non-equivalent extremal surfaces
of $w$ with boundary corresponding to $\mu$. If $w_{1}$ and $w_{2}$
induce the same measure on $\U\left(N\right)$ for all $N$, then
$\mathbb{E}_{w_{1}}\left[\zeta_{\mu}\right]=\mathbb{E}_{w_{2}}\left[\zeta_{\mu}\right]$,
so $d_{w_{1},\mu}=d_{w_{2},\mu}$.
\end{proof}
\begin{rem}
\label{rem:infinitely many extremal lambdas for U}If $\Sigma,f$
is an extremal surface for $w$, then any topological finite-degree
cover is also an extremal surfaces. We obtain that there are infinitely
many partitions $\lambda$ for which $\beta(w,\xi^{\lambda[\bullet]})=2\scl(w)$.
\end{rem}

\begin{rem}
\label{rem:scl(1)}Consider the trivial element $1_{\F}\in\F$. The
classical definition \eqref{eq:scl, old def} of $\scl$ gives $\scl\left(1_{\F}\right)=0$.
It may make more sense to define $\scl\left(1_{\F}\right)=-\frac{1}{2}$.
This would agree with \eqref{eq:scl in free groups restricting definition},
with the formula of $\scl$ for surface words (see Table \ref{tab:values on selected words}),
and, most convincingly, would make Theorem \ref{thm:MP scl for U(N)}
apply also for $w=1$.
\end{rem}

\begin{rem}
The entire discussion about word measures in the current paper revolves
around the asymptotics of stable characters. It is natural to wonder
about bounds on ``large-degree'' characters, or even on all characters
of a given group, under word measures. Recent achievements in this
direction are \cite{avni2022fourier,avni2024fourier}.
\end{rem}

\section{Stable primitivity rank and stable characters of $S_{\bullet}$\label{sec:Stable-primitivity-rank and stable characters of S}}

Stable primitivity rank was introduced by Wilton in \cite[Def.~10.6]{wilton2024rational}
as a special case of his \emph{maximal irreducible curvature}, denoted
$\rho_{+}\left(X\right)$, which he defines for every 2-complex $X$
\cite[Def.~6.3]{wilton2024rational}. Wilton defines $\sp\left(w\right)\defi1-\rho_{+}\left(X\right)$,
where $X$ is the standard presentation complex of $\nicefrac{\F}{\left\langle \left\langle w\right\rangle \right\rangle }$.
He also proves in \cite[Thm.~A]{wilton2022rationality} that $\rho_{+}\left(X\right)$,
and therefore $\sp\left(w\right)$, is rational and algorithmically
computable. Here we give our own perspective on $\sp$, and explain
its connection to word measures on $S_{\bullet}$. 

\subsection{Stable primitivity rank\label{subsec:Stable-primitivity-rank}}

The notion of primitivity rank \eqref{eq:def of pi} was conceived
in \cite{puder2014primitive} as part of a work on word measures on
$S_{\bullet}$. Subsequently, it was discovered that $\pi\left(w\right)$
played an important role in the structure of the one-relator group
defined by $w$ \cite{louder2022negative,linton2022one}. Although
stable primitivity rank was originally introduced in \cite{wilton2024rational}
without direct connection to word measures, it also arises very naturally
from word measures on $S_{\bullet}$, as we explain in $\S\S$\ref{subsec:Evidence-towards-Conjecture}. 

Recall that $\pi\left(w\right)$ is defined as the smallest rank of
a subgroup $H\le\F$ in which $w$ is a non-primitive element. Equivalently,
it is the smallest rank of a proper algebraic extension of $\left\langle w\right\rangle $
inside $\F$. We already mentioned in $\S\S$\ref{subsec:Stable-primitivity-rank and S}
that a naive attempt to imitate the classical definition \eqref{eq:scl, old def}
of $\scl$ and ``stabilize'' the primitivity rank by $\inf_{n}\frac{\pi\left(w^{n}\right)}{n}$
is inutile: for $n\ge2$, $w^{n}$ is a non-primitive element in $\left\langle w\right\rangle $,
so $\pi\left(w^{n}\right)=1$. A possible solution is to look for
algebraic extensions of $\left\langle w^{n}\right\rangle $ that give
us 'genuine' new algebraic extensions, which were not relevant for
smaller powers of $w$. This is equivalent to looking for algebraic
extensions of $\left\langle w^{n}\right\rangle $ which do not contain
smaller powers of $w$ (if $w\in H\le\F$ and $n\ge1$, then $H$
is algebraic over $\left\langle w\right\rangle $ if and only if it
is algebraic over $\left\langle w^{n}\right\rangle $). We are led
to the following definition, introduced above in \eqref{eq:spi- short but unverified definition}:
\begin{equation}
\tilde{\sp}\left(w\right)\defi\inf\left\{ \frac{\rk H-1}{n}\,\middle|\,\begin{gathered}n\in\mathbb{Z}_{\ge1},w^{n}\in H\in\text{\ensuremath{\F}},\\
w^{n}~\mathrm{non\text{-}primitive~in}~H,\\
\forall1\le n'<n:~w^{n'}\notin H
\end{gathered}
\right\} .\label{eq:spi tilde}
\end{equation}
(The ``$-1$'' in the numerator makes sense for several reasons
and is meant to make this definition inline with that of $\sp$ below).
Note that $\tilde{\sp}\left(1_{\F}\right)=-1$ (obtained with $H=\left\{ 1_{\F}\right\} $).
If $w\ne1$ is a proper power, then $\tilde{\sp}\left(w\right)=0$.
If $w\ne1$ is not a proper power, then 
\[
\tilde{\sp}\left(w\right)=\inf\left\{ \frac{\rk H-1}{n}\,\middle|\,\begin{gathered}n\in\mathbb{Z}_{\ge1},w^{n}\in H\in\text{\ensuremath{\F}},\\
w^{n}~\mathrm{non\text{-}primitive~and~non\text{-}power~in}~H
\end{gathered}
\right\} .
\]
The definition in \eqref{eq:spi tilde} can be seen as an analog of
the classical definition \eqref{eq:scl, old def} of $\scl$. However,
a more natural definition is one in the spirit of Calegari's alternative
definition \eqref{eq:scl a la Calegari} of $\scl$. Such a definition
considers subgroups containing not only a single power of $w$, but
a set of different conjugates of powers. Indeed, this is the version
of the definition that arises naturally from analyzing word measures
in $S_{\bullet}$ (see $\S\S$\ref{subsec:Evidence-towards-Conjecture}).
It is also the one introduced by Wilton in \cite{wilton2024rational},
which allows him to generalize Calegari's techniques and prove rationality
and computability of stable primitivity rank in \cite{wilton2022rationality}.
We present the definition in terms of core graphs. (Recall the notation
$\Gamma_{w}$ from Definition \ref{def:Gamma_w}.)
\begin{defn}
\label{def:spi}Let $1\ne w\in\F$. The \textbf{stable primitivity
rank} of $w$ is 
\begin{equation}
\sp\left(w\right)\defi\inf\left\{ \frac{-\chi\left(\Gamma\right)}{\deg\left(\rho\right)}\,\middle|\,\begin{gathered}\xymatrix{P\ar@{->>}[d]^{\rho~~}\ar[r]^{b} & \Gamma\ar[d]^{f}\\
\Gamma_{w}\ar[r]^{\eta_{w}} & \Omega
}
\end{gathered}
~\begin{gathered}\mathrm{s.t.}~\Gamma~\mathrm{is~a~core~graph~and}~f~\mathrm{an~immersion,}\\
\rho~\mathrm{is~a~finite\text{-}degree~covering~map,}\\
b~\mathrm{is~algebraic,the~diagram~efficient,and}\\
b~\mathrm{not~an~isomorphism~on~any~component}
\end{gathered}
\right\} .\label{eq:def of spi with core graphs}
\end{equation}
The last condition means that $b|_{b^{-1}\left(\Gamma_{0}\right)}\colon b^{-1}\left(\Gamma_{0}\right)\to\Gamma_{0}$
is not an isomorphism for any connected component $\Gamma_{0}$ of
$\Gamma$. We also define $\sp\left(1_{\F}\right)\defi-1$.
\end{defn}

As above, $\Omega$ is a bouquet with $\pi_{1}\left(\Omega,o\right)=\F$.
The map $\rho$ is a covering map of degree $\deg\left(\rho\right)$,
so $P$ is a disjoint union of cycle-graphs corresponding to different
powers of $w$ with total exponent $\deg\left(\rho\right)$.
\begin{rem}[Comments around Definition \ref{def:spi}]
\label{rem:about the definition of spi}
\begin{enumerate}
\item The requirement that $b$ be not an isomorphism on any connected component
of $\Gamma$ is needed because for every $w$ and every connected
cover $P$, the trivial diagram with $b\colon P\stackrel{\cong}{\to}\Gamma$
an isomorphism satisfies all the other requirements but represents
a trivial algebraic extension giving $\frac{-\chi\left(\Gamma\right)}{\deg\left(\rho\right)}=\frac{0}{\deg\left(\rho\right)}=0$.
\item As $\Gamma$ is a core graph with $f\colon\Gamma\to\Omega$ an immersion,
it is a multi core graph in the sense of \cite{hanany2020word} (see
$\S\S$\ref{subsec:Core-graphs}). 
\item The efficiency is required to guarantee that $b$ does not factor
as $P\stackrel{h}{\to}P'\stackrel{b'}{\to}\Gamma$ where $P'$ is
a cover of $\mathbb{S}^{1}$ of degree $<\deg\left(\rho\right)$.
This generalizes the requirement in \eqref{eq:spi tilde} that $w^{n'}\notin H$.
\item \label{enu:pullbacks}The fact that $f$ is an immersion means it
is $\pi_{1}$-injective, and that, moreover, every non-backtracking
path in $\Gamma$ maps to a non-backtracking path in $\Omega$. So
the efficiency of the diagram in \eqref{eq:def of spi with core graphs}
is equivalent to that $b$ is injective on some (or, equivalently,
every) fiber of $\rho$. Efficiency here is also equivalent to that
$P$ is the pullback of $w$ and $f$ in the sense of \cite[P.~9248]{hanany2020word}
(obtained after omitting contractible components).
\item \label{enu:may demand connectivity }We may demand in Definition \ref{def:spi}
that $\Gamma$ be connected: indeed, for not-necessarily-connected
$\Gamma$ assume that $\Gamma=\Gamma_{1}\sqcup\ldots\sqcup\Gamma_{m}$
are its connected components and let $\rho_{i}\colon P_{i}\to\mathbb{S}^{1}$
is the restriction of $\rho$ to $P_{i}=b^{-1}\left(\Gamma_{i}\right)$.
Then there exists some $i\in\left\{ 1,\ldots,m\right\} $ so that
\[
\frac{-\chi\left(\Gamma_{i}\right)}{\deg\left(\rho_{i}\right)}\le\frac{-\chi\left(\Gamma_{1}\right)-\ldots-\left(\Gamma_{m}\right)}{\deg\left(\rho_{1}\right)+\ldots+\deg\left(\rho_{m}\right)}=\frac{-\chi\left(\Gamma\right)}{\deg}.
\]
So we may restrict to connected core graphs $\Gamma$, and then the
requirement on $b$ is that it is algebraic and not an isomorphism. 
\item \label{enu:replace algebraicity with non-free}We may replace the
requirement that $b$ be algebraic, with $b$ being non-free (see
Definition \ref{def:algebraic morphisms}) on every component of $\Gamma$.
By the algebraic-free decomposition \eqref{eq:algebraic-free decomposition}
and \cite[Prop.~4.3(3)]{hanany2020word}, this does not change the
infimum. The condition that $b$ is not an isomorphism on any component
is then redundant: an isomorphism is always free. Combined with Item
\ref{enu:may demand connectivity }, we may simply demand that $\Gamma$
be connected and $b$ non-free.
\end{enumerate}
\end{rem}

\begin{claim}
\label{claim:sp is Aut invariant} Stable primitivity rank is $\mathrm{Aut}\F$-invariant.\footnote{This also follows from $\sp$ being profinite \cite[Thm.~1.2]{PSh25}.}
\end{claim}

\begin{proof}
This follows from the fact that Definition \ref{def:spi} can be given
in completely group-theoretic terms. Indeed, by Remark \ref{rem:about the definition of spi}\eqref{enu:may demand connectivity },
we may restrict to diagram as in \eqref{eq:def of spi with core graphs}
with $\Gamma$ connected. For any such diagram, pick an arbitrary
vertex $v\in V\left(\Gamma\right)$ and let $H=f_{*}\left(\pi_{1}\left(\Gamma,v\right)\right)\le\F$.
Assume that $P=P_{1}\sqcup\ldots\sqcup P_{m}$ are the connected components
of $P$ with $\rho|_{P_{i}}\colon P_{i}\to\Gamma_{w}$ a degree-$n_{i}$
cover. For every $i=1,\ldots,m$, the path $b\left(P_{i}\right)$
corresponds to the conjugacy class in $H$ of $u_{i}w^{n_{i}}u_{i}^{-1}$
for some $u_{i}\in\F$. The assumption that $b$ is algebraic is equivalent
to that $H$ is algebraic over the multiset of subgroups $\left\{ \left\langle u_{1}w^{n_{1}}u_{1}^{-1}\right\rangle ,\ldots,\left\langle u_{m}w^{n_{m}}u_{m}^{-1}\right\rangle \right\} $
(see $\S\S$\ref{subsec:Algerbaic-morphisms}). The assumption that
the diagram is efficient, namely, that $b$ is injective on every
fiber of $\rho$, is equivalent to that for any $i,j\in\left\{ 1,\ldots,m\right\} $
and $k\in\mathbb{Z}$, if $u_{i}w^{k}u_{j}^{-1}\in H$, then $i=j$
and $n_{i}\mid k$. 
\end{proof}
It is an easy observation that a map such as $b\colon P\to\Gamma$
from a collection of cycles to a core graph which is algebraic and
not an isomorphism on any connected component, must cover every edge
at least twice (e.g., \cite[Lem.~4.1]{puder2015expansion}). Helfer
and Wise \cite[Thm.~4.1]{helfer2016counting}, and independently Louder
and Wilton \cite[Thm.~1.2]{louder2017stackings}, both proved that
if one considers \eqref{eq:def of spi with core graphs} with the
requirement that $b$ is algebraic and non-isomorphism on any component
replaced by the requirement that $\left|b^{-1}\left(e\right)\right|\ge2$
for every $e\in E\left(\Gamma\right)$, then the infimum is $\ge1$
for non-powers. The following immediate corollary also appears as
\cite[Thm.~11.15]{wilton2024rational}.
\begin{thm}
\cite{helfer2016counting,louder2017stackings}\label{thm:HW-LW} If
$1\ne w\in\F$ is not a proper power then $\sp\left(w\right)\ge1$.
\end{thm}

\begin{cor}
\label{cor:ineq of spi, spi tilde and pi}Let $1\ne w\in\F$ be a
non-power. Then
\begin{equation}
1\le\sp\left(w\right)\le\tilde{\sp}\left(w\right)\le\pi\left(w\right)-1.\label{eq:inequalities of sp,spi tilde and pi}
\end{equation}
If $1\ne w$ is a proper power then $\sp\left(w\right),\tilde{\sp}\left(w\right)$
and $\pi\left(w\right)-1$ are all zero, if $w$ is primitive they
are all $\infty$, and if $w=1$ they are all $-1$. 
\end{cor}

\begin{proof}
The first inequality in \eqref{eq:inequalities of sp,spi tilde and pi}
is Theorem \ref{thm:HW-LW}. The remaining inequalities hold because
$\tilde{\sp}\left(w\right)$ is equal to the infimum over all diagrams
in Definition \ref{def:spi} where the cover $P$ is \emph{connected},
and $\pi\left(w\right)-1$ the infimum over all diagrams with a \emph{trivial
}(degree-1) cover. The other claims are immediate.
\end{proof}
For non-primitive words, we trivially have $\pi\left(w\right)\le\rk\F$,
so for a non-primitive non-power $w$ we have $\sp\left(w\right)\in\left[1,\rk\F-1\right]$.
In particular, in $\F_{2}$, the rank-2 free groups, every $1\ne w\in\F_{2}$
which is not primitive and not a proper power satisfies 
\[
\sp\left(w\right)=\tilde{\sp}\left(w\right)=\pi\left(w\right)-1=1.
\]
As hinted above, we suspect that as in the case of $\scl$, the restriction
to connected covers does not change the outcome.
\begin{conjecture}
\label{conj:spi and spi tilde}For all $w\in\F$ we have $\sp\left(w\right)=\tilde{\sp}\left(w\right)$.
\end{conjecture}

In fact, the following stronger conjecture was presented by Wilton
as a ``stability conjecture'' in an online talk \cite[Minute 57]{wilton2021youtbue}.
It also follows indirectly from \cite[Conj.~1.8]{hanany2020word}
-- see Remark \ref{rem:HP conj on S}.
\begin{conjecture}[Wilton]
\label{conj:spi=00003Dpi-1} For all $w\in\F$ we have $\sp\left(w\right)=\pi\left(w\right)-1$.
\end{conjecture}

\begin{rem}
\label{rem:spi and covers}Another advantage of $\sp$ over $\tilde{\sp}$
is that it works better with covers. Recall from Remark \ref{rem:about the definition of spi}\eqref{enu:pullbacks}
that instead of efficiency of the diagrams in \eqref{eq:def of spi with core graphs},
we may define $P$ to be the pullback of $f$ and $\eta_{w}$. It
is easy to see that we may replace $\Gamma$ by any degree-$k$ cover
of $\Gamma$ (and replace $P$ with the appropriate pullback) to obtain
another valid diagram with the same quotient $\frac{-\chi\left(\Gamma\right)}{\deg\left(\rho\right)}$.
Thus, every valid, extremal diagram for $w$, where the quotient gives
the precise value of $\sp\left(w\right)$, gives rise to infinitely
many such valid, extremal diagrams by covers. 
\end{rem}

Definition \ref{def:spi} looks similar to the definition of $\scl$
in free groups from Corollary \ref{cor:scl in free groups - def with WH graphs},
but the following version makes the connection even more apparent.
\begin{prop}
\label{prop:spi in terms of WH}For every $1\ne w\in\F$, 
\begin{equation}
\sp\left(w\right)=\inf\left\{ \frac{-\chi\left(\Gamma\right)}{\deg\left(\rho\right)}\,\middle|\,\begin{gathered}\xymatrix{P\ar@{->>}[d]^{\rho~~}\ar[r]^{b} & \Gamma\ar[d]^{f}\\
\Gamma_{w}\ar[r]^{\eta_{w}} & \Omega
}
\end{gathered}
~\begin{gathered}\mathrm{s.t.}~\Gamma~\mathrm{is~a~core~graph,}\\
\rho~\mathrm{is~a~finite\text{-}degree~covering~map},\\
\forall v\in V\left(\Gamma\right),~\wh_{b}\left(v\right)~\mathrm{is~connected,without}\\
\mathrm{cut\text{-}vertices~and}~\wh_{b}\left(v\right)\ne\xymatrix{\bullet\ar@{-}[r] & \bullet}
\\
\mathrm{and~the~diagram~is~efficient}
\end{gathered}
\right\} .\label{eq:def of spi with WH}
\end{equation}
\end{prop}

Note that in \eqref{eq:def of spi with WH} we omitted the requirements
from \eqref{eq:def of spi with core graphs} that b is algebraic and
not an isomorphism on connected components, and that $f$ is an immersion
(so $\Gamma$ is no longer necessarily a multi core graph in the sense
of \cite{hanany2020word}). The efficiency of the diagram in \eqref{eq:def of spi with WH}
is much more subtle than in \eqref{eq:def of spi with core graphs}
and has the same meaning as in Corollary \ref{cor:scl in free groups - def with WH graphs}.
\begin{proof}[Proof of Proposition \ref{prop:spi in terms of WH}]
 We want to show that every diagram in \eqref{eq:def of spi with WH}
gives rise to a diagram as in \eqref{eq:def of spi with core graphs}
and vice versa. Assume first we are given a diagram as in \eqref{eq:def of spi with WH}.
There is a unique folded version of $f$ à la Stallings, namely, $f$
factors uniquely as $\Gamma\stackrel{f_{1}}{\to}\Gamma'\stackrel{f_{2}}{\to}\Omega$,
where $f_{1}$ is obtained by a sequence of folds, and $f_{2}$ is
an immersion \cite[\S3.3]{stallings1983topology}. We claim that the
diagram
\begin{equation}
\begin{gathered}\xymatrix{P\ar@{->>}[d]^{\rho~~}\ar[r]^{f_{1}\circ b} & \Gamma'\ar[d]^{f_{2}}\\
\Gamma_{w}\ar[r]^{\eta_{w}} & \Omega
}
\end{gathered}
\label{eq:after folding}
\end{equation}
satisfies the assumptions of \eqref{eq:def of spi with core graphs}.
Indeed, the efficiency of the original diagram means that $f$ is
$\pi_{1}$-injective, therefore so is $f_{1}$. This means that all
the folds constituting $f_{1}$ are homotopy-equivalent folds and
$f_{1}$ is an isomorphism in the level of $\pi_{1}$. Thus $f_{1}\circ b$
is algebraic. The conditions on $\wh_{b}\left(v\right)$ in \eqref{eq:def of spi with WH}
guarantee that $\wh_{b}\left(v\right)$ has no vertices of degree
$\le1$, so $\left|b^{-1}\left(e\right)\right|\ge2$ for every $e\in E\left(\Gamma\right)$.
Folding can only increase the size of the preimage of an edge, so
$\left|(f_{1}\circ b)^{-1}\left(e\right)\right|\ge2$ for every $e\in E\left(\Gamma'\right)$.
In particular, $f_{1}\circ b$ is not an isomorphism on any connected
component of $\Gamma'$. The efficiency of the initial diagram guarantees
the efficiency of \eqref{eq:after folding}: $f_{1}$ merges the vertices
$u,v\in V\left(\Gamma\right)$ if and only if there is a path from
$u$ to $v$ that maps to a null-homotopic path in $\Omega$. Finally,
because $f_{1}$ is a homotopy equivalence, we have $\chi\left(\Gamma'\right)=\chi\left(\Gamma\right)$,
so $\frac{-\chi\left(\Gamma'\right)}{\deg\left(\rho\right)}=\frac{-\chi\left(\Gamma\right)}{\deg\left(\rho\right)}$.

Conversely, assume we are given a diagram as in \eqref{eq:def of spi with core graphs}.
By Corollary \ref{cor:checking algebraicity by unfolding}, $b$ factors
as $P\stackrel{\overline{b}}{\to}\overline{\Gamma}\stackrel{b_{2}}{\to}\Gamma$
where $b_{2}$ is obtained by a sequence of homotopy-equivalent folds,
and all Whitehead graphs of $\overline{b}$ are connected and without
cut-vertices. We claim that 
\begin{equation}
\begin{gathered}\xymatrix{P\ar@{->>}[d]^{\rho~~}\ar[r]^{\overline{b}} & \overline{\Gamma}\ar[d]^{f\circ b_{2}}\\
\Gamma_{w}\ar[r]^{\eta_{w}} & \left(\Omega,o\right)
}
\end{gathered}
\label{eq:after unfolding}
\end{equation}
satisfies the assumptions in \eqref{eq:def of spi with WH}. Again,
because $b_{2}$ induces an isomorphism of fundamental groups, the
algebraicity of $b$ is equivalent to that $\overline{b}$. Because
$\eta_{w}\circ\rho$ is an immersion, so is $\overline{b}$, and together
with the conditions on the Whitehead graphs of $\overline{b}$ we
get that $\overline{\Gamma}$ is a core graph. If for some $v\in V\left(\overline{\Gamma}\right)$
we have $\wh_{\overline{b}}\left(v\right)=\xymatrix{\bullet\ar@{-}[r] & \bullet}
$, then $\overline{b}$ must be an isomorphism of cycles, which means
that so must $b$, contradicting the assumption in \eqref{eq:def of spi with core graphs}.
The efficiency of the original diagram clearly yields the efficiency
of \eqref{eq:after unfolding}. Finally, because $b_{2}$ is an homotopy
equivalence, $\frac{-\chi\left(\overline{\Gamma}\right)}{\deg\left(\rho\right)}=\frac{-\chi\left(\Gamma\right)}{\deg\left(\rho\right)}$. 
\end{proof}
The following corollary is also recorded in \cite[Lem.~10.9]{wilton2024rational}.
\begin{cor}
\label{cor:spi le 2scl}For every $w\in\F$ we have $\sp\left(w\right)\le2\scl\left(w\right)$.
\end{cor}

\begin{proof}
Notice that if a (Whitehead) graph is a cycle, then, in particular,
it is connected, admits no cut-vertices, and is not equal to $\xymatrix{\bullet\ar@{-}[r] & \bullet}
$. So every diagram in Corollary \ref{cor:scl in free groups - def with WH graphs}
satisfies all the conditions in \eqref{eq:def of spi with WH}, and
the infimum in Corollary \ref{cor:scl in free groups - def with WH graphs}
is over a subset of the set in Proposition \ref{prop:spi in terms of WH}.
\end{proof}
In fact, there is generically a large difference between $\sp$ and
$\scl$. While $\sp$ is bounded from above by $\rk\F-1$ for all
non-primitive words, $\scl$ is $\infty$ whenever $w\notin\left[\F,\F\right]$,
and within $\left[\F,\F\right]$ its value for random words of length
$n$ is of order $\frac{n}{\log n}$ \cite{calegari2013random}. By
Corollary \ref{cor:spi le 2scl}, Conjecture \ref{conj:spi=00003Dpi-1}
yields a conjecture of Heuer that $\pi\left(w\right)-1\le2\scl\left(w\right)$
\cite[Conj.~6.3.2]{heuer2019constructions}.
\begin{rem}
The ordinary primitivity rank can be defined for subgroups of $\F$
as well as for words \cite[Def.~1.7]{PP15}. The definition of stable
primitivity rank can be extended to subgroups too - see \cite{shomroni2023probabilisticHNC}.
\end{rem}

\subsection{Evidence towards Conjecture \ref{conj:spi and stable characters of S_N}\label{subsec:Evidence-towards-Conjecture}}

Recall that $\irr\left(S_{N}\right)$, the complex irreducible characters
of $S_{N}$, are parameterized by partitions of $N$ or, equivalently,
by Young diagrams with $N$ squares. We let $\chi^{\mu}$\marginpar{$\chi^{\mu}$}
denote the irreducible character of $S_{\left|\mu\right|}$ corresponding
to the partition $\mu=(\mu_{1},\ldots,\mu_{\ell(\mu)})$. A stable
irreducible representation of $S_{\bullet}=\{S_{N}\}_{N}$ is given
by padding a given partition $\mu$ with a large first number. Namely,
for every $N\ge\left|\mu\right|+\mu_{1}$, we denote by\footnote{We use the same notation here as in stable sequences of highest weight
vectors in the representation theory of $\U\left(\bullet\right)$.
The correct meaning of the notation $\lambda\left[N\right]$ (or $\mu\left[N\right]$)
each time should be inferred from the context.} $\mu\left[N\right]$ the partition of $N$ given by $\left(N-\left|\mu\right|,\mu_{1},\mu_{2},\ldots,\mu_{\ell}\right)$,
and the corresponding stable irreducible character by \marginpar{$\chi^{\mu[\bullet]}$}$\chi^{\mu[\bullet]}=\left\{ \chi^{\mu\left[N\right]}\right\} _{N\ge\left|\mu\right|+\mu_{1}}$.
For example, if $\mu=\left(1\right)\vdash1$, then for every $N\ge2$,
$\mu\left[N\right]=\left(N-1,1\right)$ and $\chi^{\mu[\bullet]}$
is the standard character from \eqref{eq:pp15, compact version}. 

As mentioned above, it is shown in \cite{hanany2020word} that for
every word $w$ and partition $\mu$ the values of $\mathbb{E}_{w}\left[\chi^{\mu\left[N\right]}\right]$
coincide with a rational function $f(w,\mu)\in\mathbb{Q}(N)$ for
every large enough $N$. It follows immediately from the hook length
formula that $\dim\chi^{\mu\left[N\right]}=\Theta\left(N^{\left|\mu\right|}\right)$,
so 
\begin{equation}
\beta(w,\chi^{\mu[\bullet]})=-\frac{\deg(f(w,\mu)}{\left|\mu\right|}\in\mathbb{Q}_{\ge-1}\cup\left\{ \infty\right\} \label{eq:beta for S_bullet}
\end{equation}
is well defined whenever $\mu\ne\emptyset$ (namely, whenever the
stable irreducible character is not trivial). 
\begin{example}
\label{exa:fractional beta in S_n}For $w=a^{2}b^{2}a^{-1}b$ and
$\mu=(2,1)\vdash3$, 
\begin{equation}
\mathbb{E}_{a^{2}b^{2}a^{-1}b}\left[\chi^{(2,1)[N]}\right]=\frac{12N^{2}-91N+126}{N^{7}-21N^{6}+175N^{5}-735N^{4}+1624N^{3}-1764N^{2}+720N}\label{eq:example of non-integral beta}
\end{equation}
for $N\ge9$, so $\beta(a^{2}b^{2}a^{-1}b,\chi^{(2,1)[\bullet]})=\frac{5}{3}$.
We thank Noam Ta Shma for computing \eqref{eq:example of non-integral beta}.
\end{example}

Recall from $\S\S$\ref{subsec:Stable-primitivity-rank and S} that
${\cal I}_{S}$ denotes the set of all stable irreducible characters
of $S_{\bullet}$. Conjecture \ref{conj:spi and stable characters of S_N}
states that $\sp(w)=\inf_{\triv\ne\chi\in{\cal I}_{S}}\beta(w,\chi)$
and that the infimum is attained. There are many cases where this
conjecture is known:
\begin{itemize}
\item \textbf{The identity element}: Trivially, $\sp(1_{\F})=-1$ and $\beta(1_{\F},\chi)=-1$
for all $\triv\ne\chi\in{\cal I}_{S}$.
\item \textbf{Surface words:} By Frobenius and Schur, as the Frobenius-Schur
indicator of all irreducible characters of the symmetric group is
$1$, surface words satisfy $\beta(w_{g},\chi)=2g-1$ and $\beta(u_{g},\chi)=g-1$
for all $\triv\ne\chi\in{\cal I}_{S}$, and it is not hard to show
that these are precisely the corresponding values of $\sp$ of these
words (e.g., \cite[Example 7.5]{wilton2024rational}).
\item \textbf{Proper powers: }Let $w\ne1$ be a proper power. As mentioned
above, $\sp(w)=0$. Moreover, $\beta(w,\chi^{(1)[\bullet]})=0$ (this
is a special case of \cite[Thm.~1.8]{PP15} mentioned in \eqref{eq:pp15, compact version}
above, which was already proven in \cite{nica1994number}). Finally,
$\beta(w,\chi)\ge0$ holds for any $\triv\ne\chi\in{\cal I}_{S}$
by combining the formula in \cite[Prop.~B.2]{hanany2020word} and
the results in \cite[\S4]{linial2010word} (see also \cite[Thm.~1.15]{puder2024local}).
\item \textbf{Primitive words: }Let $w\in\F$ be primitive. As mentioned
in $\S\S$\ref{subsec:Stable-primitivity-rank}, $\sp(w)=\infty$.
Because $w$ induces the uniform distribution on every finite group
(e.g., \cite[Obs.~1.2]{PP15}), we have $\mathbb{E}_{w}[\chi_{N}]=0$
for every non-trivial character $\chi_{N}$ of $S_{N}$, so $\beta(w,\chi)=\infty$
for all $\triv\ne\chi\in{\cal I}_{S}$.
\item \textbf{Words with $\pi(w)=2$: }In this case, $1\le\sp(w)\le\pi(w)-1=1$
so $\sp(w)=1$. By \eqref{eq:pp15, compact version}, $\beta(w,\chi^{(1)[\bullet]})=1$,
and by Cassidy's result \cite[Thm.~1.5]{cassidy2024random} $\beta(w,\chi)\ge1$
for any $\triv\ne\chi\in{\cal I}_{S}$.
\item \textbf{Rank-2 free group:} Every element in the rank-2 free group
$\F_{2}$ falls into one of the previous cases.
\end{itemize}
The above cases are all, in a sense, degenerate cases, where we do
not see the full definition of $\sp$ at play. However, we now present
a formula which shows how the diagrams in the definition of $\sp$
show up in the expected value of stable irreducible characters, thus
giving more evidence towards Conjecture \ref{conj:spi and stable characters of S_N}.
The formula uses some additional notation. First, let $\eta\colon\Gamma\to\Delta$
be an immersion of core graphs. Denote by $\algdecompt\left(\eta\right)$
the set of decompositions of $\eta$ into a triple of immersions $\left(\eta_{1},\eta_{2},\eta_{3}\right)$
such that $\eta=\eta_{3}\circ\eta_{2}\circ\eta_{1}$, and $\eta_{1}$
and $\eta_{2}$ are algebraic. As algebraic morphisms of core graphs
are surjective \cite[Thm.~4.7]{hanany2020word}, $|\algdecompt\left(\eta\right)|<\infty$.
Second, if $\eta$ as above is algebraic, there is a map $C_{\eta}^{\mathrm{alg}}\colon\mathbb{Z}_{\ge1}\to\mathbb{Q}$,
introduced in \cite[\S\S6.2]{hanany2020word} (following \cite[\S5]{PP15}),
which is obtained as a Möbius inversion of functions generalizing
$\mathbb{E}_{w}[\mathrm{tr}_{N}]$. By \cite[Cor.~6.16]{hanany2020word},
\begin{equation}
C_{\eta}^{\mathrm{alg}}\left(N\right)=\begin{cases}
N^{\chi\left(\Gamma\right)}, & \mathrm{if}~\eta=\id~\mathrm{is~an~isomorphism}\\
O\left(N^{\min\left(\chi\left(\Gamma\right),\chi\left(\Delta\right)\right)-1}\right), & \mathrm{otherwise.}
\end{cases}\label{eq:bounds on C^alg}
\end{equation}

\begin{thm}
\label{thm:formula for stable irreps of S_N}Let $d\in\mathbb{Z}_{\ge0}$,
let $\mu\vdash d$ be a partition and let $1\ne w\in\F$ be a non-power.
For every $N\ge2d$ we have 
\begin{equation}
\mathbb{E}_{w}\left[\chi^{\mu\left[N\right]}\right]=\frac{1}{d!}\sum_{\sigma\in S_{d}}\chi^{\mu}\left(\sigma\right)\sum_{\substack{\left(\eta_{1},\eta_{2},\eta_{3}\right)\in\algdecompt\left(w^{\sigma}\to\F\right)\colon\\
\eta_{1}~\mathrm{is~efficient,~and}\\
\mathrm{codomain}\left(\eta_{2}\right)~\mathrm{has~no~cycles}
}
}C_{\eta_{2}}^{\mathrm{alg}}\left(N\right).\label{eq:formula for stable irreps of S_N}
\end{equation}
\end{thm}

We defer the proof of Theorem \ref{thm:formula for stable irreps of S_N}
to the companion paper \cite{PSh25} (where it appears as Theorem
1.4). We explain here how it relates to Conjecture \ref{conj:spi and stable characters of S_N}.
In \eqref{eq:formula for stable irreps of S_N}, $w^{\sigma}\to\F$
is the morphism $\eta_{w^{\sigma}}\colon\Gamma_{w^{\sigma}}\to\Omega$
from Definition \ref{def:Gamma_w}. Let $\Gamma_{w^{\sigma}}\stackrel{\eta_{1}}{\to}\Sigma\stackrel{\eta_{2}}{\to}\Sigma'\stackrel{\eta_{3}}{\to}\Omega$
be an element of $\algdecompt\left(\eta_{w^{\sigma}}\right)$. By
``$\eta_{1}$ is efficient'' we mean that the diagram
\begin{equation}
\xymatrix{\Gamma_{w^{\sigma}}\ar[d]_{\rho}\ar[r]^{\eta_{1}} & \Sigma\ar[d]^{\eta_{3}\circ\eta_{2}}\\
\Gamma_{w}\ar[r]_{\eta_{w}} & \Omega
}
\label{eq:diagram in formula for irrep of S}
\end{equation}
is efficient, namely, that $\eta_{1}$ is injective on each fiber
of $\rho$ ($\eta_{3}\circ\eta_{2}$ is an immersion and so automatically
$\pi_{1}$-injective). By ``$\mathrm{codomain}\left(\eta_{2}\right)$
has no cycles'', we mean that $\Sigma'$ has no connected component
which is a cycle-graph. Whenever $\eta_{1}$ is \emph{not} an isomorphism
on any connected component of $\Sigma$, the diagram \eqref{eq:diagram in formula for irrep of S}
satisfies all the conditions of Definition of \ref{def:spi}, so $\sp\left(w\right)\le\frac{-\chi\left(\Sigma\right)}{\deg\left(\rho\right)}=\frac{-\chi\left(\Sigma\right)}{d}$.
By \eqref{eq:bounds on C^alg}, $C^{\mathrm{alg}}\left(\eta_{2}\right)=O(N^{\chi\left(\Sigma\right)})=O(N^{-d\sp\left(w\right)})=O((\dim\chi^{\mu[N]})^{-\sp(w)})$,
which agrees with the conjecture. In fact, unless $\eta_{2}$ is an
isomorphism, we get $O(N^{-d\sp\left(w\right)-1})$.

We are left with the case that $\eta_{1}$ is an isomorphism on some
connected component of its codomain $\Sigma$. So this component must
be a cycle isomorphic to one of the cycles in $\Gamma_{w^{\sigma}}$.
Recall that $\Sigma'$ contains no cycles, so $\eta_{2}$ cannot be
an isomorphism in this case. If $\eta_{2}\circ\eta_{1}$ is efficient,
then 
\[
\xymatrix{\Gamma_{w^{\sigma}}\ar[r]^{\eta_{2}\circ\eta_{1}}\ar[d]_{\rho} & \Sigma'\ar[d]^{\eta_{3}}\\
\Gamma_{w}\ar[r]_{\eta_{w}} & \Omega
}
\]
satisfies the conditions in the definition of $\sp\left(w\right)$,
and $C^{\mathrm{alg}}\left(\eta_{2}\right)=O(N^{\chi\left(\Sigma'\right)-1})=O(N^{-d\sp\left(w\right)-1})$.
The only case where we do not have this bound is when $\Sigma$ contains
connected components that are cycles, $\Sigma'$ does not, and $\eta_{2}\circ\eta_{1}$
is not efficient. We conjecture these elements, too, contribute $O(N^{-d\sp\left(w\right)-1})$.
\begin{conjecture}
\label{conj:C^alg from cycle and non-efficient}Let $1\ne w\in\F$
be a non-power, $\sigma\in S_{d}$, and let $\Gamma_{w^{\sigma}}\stackrel{\eta_{1}}{\to}\Sigma\stackrel{\eta_{2}}{\to}\Sigma'\stackrel{\eta_{3}}{\to}\Omega$
be in $\algdecompt\left(w^{\sigma}\to\F\right)$ such that $\eta_{1}$
is efficient, $\Sigma$ contains a cycle (as a connected component),
$\Sigma'$ does not, and $\eta_{2}\circ\eta_{1}$ is not efficient.
Then
\[
C_{\eta_{2}}^{\mathrm{alg}}\left(N\right)=O\left(N^{-d\sp\left(w\right)-1}\right).
\]
\end{conjecture}

For example, consider the simple commutator $w=\left[a,b\right]\in\F_{2}$
with $\sigma=\left(12\right)\in S_{2}$, $\eta_{1}=\id$ and $\Sigma'=\Omega$,
so $\eta_{2}=\eta_{w^{2}}$. In this case $\sp\left(w\right)=\pi\left(w\right)-1=1$.
While the bound \eqref{eq:bounds on C^alg} only gives $\mathbb{E}_{w}\left[\eta_{2}\right]=O\left(N^{-2}\right)$,
the actual value we know from direct computation is $\mathbb{E}_{w}\left[\eta_{2}\right]=\frac{4}{N\left(N-1\right)\left(N-2\right)\left(N-3\right)}$,
which is of order $N^{-4}$. (By Conjecture \ref{conj:C^alg from cycle and non-efficient},
it should be $O\left(N^{-3}\right)$.)
\begin{prop}
\label{prop:detailed conjectural picture for S}Conjecture \ref{conj:spi and stable characters of S_N}
follows from Conjecture \ref{conj:C^alg from cycle and non-efficient}.
Moreover, Conjecture \ref{conj:C^alg from cycle and non-efficient}
yields that for any $1\ne w\in\F$ a non-power and partition $\mu\vdash d$,
\begin{equation}
\mathbb{E}_{w}\left[\chi^{\mu\left[N\right]}\right]=\left[\frac{1}{d!}\sum_{\sigma\in S_{d}}\chi^{\mu}\left(\sigma\right)\cdot\mathrm{crit}_{w,\sigma}\right]N^{-d\sp\left(w\right)}+O\left(N^{-d\sp\left(w\right)-1}\right),\label{eq:precise leading term for E_w in S-1}
\end{equation}
where $\mathrm{crit}_{w,\sigma}$ is the number of diagrams
\[
\xymatrix{\Gamma_{w^{\sigma}}\ar@{->>}[d]^{\rho~~}\ar[r]^{b} & \Gamma\ar[d]^{f}\\
\Gamma_{w}\ar[r]^{\eta_{w}} & \Omega
}
\]
as in Definition \ref{def:spi} of $\sp\left(w\right)$ which achieve
the lower bound of $\sp\left(w\right)$, namely, such that $\frac{-\chi\left(\Gamma\right)}{d}=\sp\left(w\right)$.
\end{prop}

\begin{proof}
We explained above why all the summands in the formula of Theorem
\ref{thm:formula for stable irreps of S_N} for $\mathbb{E}_{w}[\chi^{\mu\left[N\right]}]$,
except for possibly those in Conjecture \ref{conj:C^alg from cycle and non-efficient}
are $O(N^{-d\sp\left(w\right)})$, so with Conjecture \ref{conj:C^alg from cycle and non-efficient}
we get $\mathbb{E}_{w}[\chi^{\mu\left[N\right]}]=O(N^{-d\sp\left(w\right)})$.
Moreover, assuming Conjecture \ref{conj:C^alg from cycle and non-efficient},
the only summands of order $\Theta(N^{-d\sp\left(w\right)})$ are
those where $\Sigma=\Sigma'$ and $\Sigma$ contains no cycles. The
number of these elements in $\algdecompt\left(w^{\sigma}\to\F\right)$
is precisely $\mathrm{crit}_{w,\sigma}$, and each one of them satisfies
$C_{\eta_{2}}^{\mathrm{alg}}\left(N\right)=1$. This gives \eqref{eq:precise leading term for E_w in S-1}. 

It is left to prove that there exist non-trivial partitions $\mu$
for which $\mathbb{E}_{w}[\chi^{\mu\left[N\right]}]=\Theta(N^{-d\sp\left(w\right)})$.
By Wilton's rationality theorem \cite{wilton2022rationality}, every
non-primitive word $w$ admits extremal diagrams for $\sp$. If there
is an extremal diagram of degree $d$, then for $\mu=\left(d\right)$,
the coefficient of $N^{-d\sp\left(w\right)}$ in \eqref{eq:precise leading term for E_w in S-1}
is $\frac{1}{d!}\sum_{\sigma\in S_{d}}\mathrm{crit}_{w,\sigma}$,
which is strictly positive, so 
\[
\mathbb{E}_{w}\left[\chi^{\left(d\right)\left[N\right]}\right]=\Theta\left(N^{-d\sp\left(w\right)}\right).
\]
\end{proof}
If Conjecture \ref{conj:C^alg from cycle and non-efficient} (or the
weaker Conjecture \ref{conj:spi and stable characters of S_N}) holds,
then $\sp\left(w\right)$ is determined by the $w$-measures on $S_{\bullet}$
and can be defined in terms of word measures on $S_{\bullet}$. In
\cite[Thm.~1.2]{PSh25} we show that $\sp$ is determined by the $w$-measures
on $S_{m}\wr S_{N}$ for all $m$ and $N$.
\begin{rem}
\label{rem:HP conj on S}Before Wilton introduced the definition of
$\sp$, Hanany and the first author conjectured in \cite[Conj.~1.8]{hanany2020word}
that $\beta(w,\chi)\ge\pi(w)-1$ for all $\triv\ne\chi\in{\cal I}_{\s}$.
The right way to think about this conjecture is as a combination of
Conjecture \ref{conj:spi and stable characters of S_N} that $\sp(w)=\inf_{\triv\ne\chi\in{\cal I}_{\s}}\beta(w,\chi)$
and Conjecture \ref{conj:spi=00003Dpi-1} that $\sp\left(w\right)=\pi\left(w\right)-1$.
In fact, Conjecture \ref{conj:spi=00003Dpi-1} follows from Conjecture
\ref{conj:spi and stable characters of S_N} and {[}Ibid., Conj.~1.8{]},
recalling from Corollary \ref{cor:ineq of spi, spi tilde and pi}
that $\sp\left(w\right)\le\pi\left(w\right)-1$.

A more general conjecture in the same paper {[}Ibid., Conj.~1.13{]}
states that for an arbitrary non-trivial stable character of any sequence
of groups with stable representation theory, $\beta(w,\chi)\ge\pi(w)-1$.
For $\U\left(\bullet\right)$, $\O\left(\bullet\right)$ and $\Sp\left(\bullet\right)$,
this conjecture would follow from Conjecture \ref{conj:ssql and its role for U,O and Sp},
Conjecture \ref{conj:spi=00003Dpi-1} and the inequality $\sp\le\ssql$
(Fact \ref{fact:basic properties of ssql}\eqref{enu:sp < ssql < 2scl}).
For the stable irreducible characters of $C_{m}\wr S_{\bullet}$ covered
in Theorem \ref{thm:result on wreath products and modulo-m spi},
it would follow from Conjecture \ref{conj:spi=00003Dpi-1} together
with the inequality $\sp\le\spm$ (Proposition \ref{prop:properties of spm}).
\end{rem}

\section{Stable square length and stable characters of $\protect\U\left(\bullet\right),\protect\O\left(\bullet\right),\protect\Sp\left(\bullet\right)$\label{sec:ssql+U + O + Sp}}

\subsection{Stable square length\label{subsec:Stable-square-length}}

In any group $G$, the square length of $h\in G$ is 
\[
\mathrm{sql}\left(h\right)=min\left\{ g\in\mathbb{Z}_{\ge0}\,\middle|\,\exists a_{1},\ldots,a_{g}\in G~\mathrm{so~that}~h=a_{1}^{~2}\cdots a_{g}^{~2}\right\} .
\]
As with commutator length, it is an easy observation \cite[\S1.1]{culler1981using}
that if $\left(X,x_{0}\right)$ is a pointed topological space with
$\pi_{1}\left(X,x_{0}\right)\cong G$, then
\[
\mathrm{sql}\left(h\right)=\min\left\{ g\in\mathbb{Z}_{\ge0}\,\middle|\,\exists f\colon\left(\Lambda_{g,1},v\right)\to\left(X,x_{0}\right)~\mathrm{so~that}~f_{*}\left(\left[\partial\Lambda_{g,1}\right]\right)=\left[h\right]\right\} ,
\]
where $\Lambda_{g,1}$ is the genus-$g$ \emph{non-}orientable surface
with one boundary component and a marked point $v$ at the boundary.
In \cite{culler1981using}, Culler presents an algorithm to compute
$\mathrm{sql}$ in free groups. In fact, this algorithm first computes
$\sqlh=\min\left(\mathrm{sql},2\cdot\cl\right)$, using the observation
that 
\[
1-\sqlh\left(h\right)=\max\left\{ \chi\left(\Sigma\right)\,\middle|\,\exists f\colon\left(\Sigma,v\right)\to\left(X,x_{0}\right)~\mathrm{\mathrm{so~that}}~f_{*}\left(\left[\partial\Sigma\right]\right)=\left[h\right]\right\} 
\]
where $\Sigma$ is a compact surface (either orientable or not) with
one boundary component and a marked point $v\in\partial\Sigma$. Given
$\sqlh\left(w\right)$ of a word, the algorithm then proceeds to computing
$\mathrm{sql}\left(w\right)$ relying on that $\sqlh$$\left(w\right)\le\mathrm{sql}\left(w\right)\le\sqlh\left(w\right)+1$
{[}Ibid., Thm.~2.2{]}. Another instance where $\mathrm{\sqlh}$ appears
rather than $\mathrm{sql}$ is in \cite[Cor.~1.11]{MPorthsymp}, which
states that in $\O\left(N\right)$ (as well as in $\Sp\left(N\right)$),
the expected value of the trace of a $w$-random matrix is $O(N^{1-\sqlh\left(w\right)})$. 

As explained in $\S\S$\ref{subsec:Stable-square-length and U O Sp},
one cannot use the naive formulas $\inf_{n}\frac{\mathrm{sql}\left(w^{n}\right)}{n}$
or $\inf_{n}\frac{\sqlh\left(w^{n}\right)}{n}$ to stabilize these
invariants. Using admissible surfaces as in \eqref{eq:scl a la Calegari}
or \eqref{eq:scl a la Calegari in F}, when now the surfaces are not
necessarily orientable and not necessarily monotone, is equally useless,
as the exact same problems arise: there are trivial solutions for
certain degrees. For example, for every $w$, there is a Möbius band
with boundary mapping to $w^{2}$ (this gives $\frac{-\chi\left(\Sigma\right)}{\deg}=\frac{0}{2}=0$),
or an annulus with two boundary components mapping each to $w$ (this,
too, gives $\frac{0}{2}=0$). However, it is possible to get a meaningful
definition if one uses the more restricted definition of $\scl$ in
free groups as in Corollary \ref{cor:scl in free groups - restricted definition}.
This definition of a stable version of $\sqlh$ is the special case
of Wilton's maximal surface curvature $\sigma_{+}$, when applied
to the presentation complex of $\nicefrac{\F}{\left\langle \left\langle w\right\rangle \right\rangle }$
\cite[Def.~6.4 and \S10.3]{wilton2024rational}. Wilton does not use
the term \emph{stable square length} in his paper.
\begin{defn}
\label{def:ssql}\cite{wilton2024rational} Let $1\ne w\in\F$. The
\emph{stable square length} of $w$ is 
\begin{equation}
\ssql\left(w\right)\defi\inf\left\{ \frac{-\chi\left(\Sigma\right)}{\deg\left(\rho\right)}\,\middle|\,\begin{gathered}\xymatrix{P\ar@{->>}[d]^{\rho~~}\ar[r]^{b} & \Sigma\ar[d]^{f}\\
\left(\mathbb{S}^{1},s_{0}\right)\ar[r]^{w} & \left(\Omega,o\right)
}
\end{gathered}
~\begin{gathered}\mathrm{s.t.~\Sigma~\mathrm{is~a~compact~surface~with}}\\
\mathrm{boundary~in~every~connected~component,}\\
\rho~\mathrm{is~a~covering~map~of~degree}~\deg\left(\rho\right),\\
b~\mathrm{is~an~isomorphism}~P\cong\partial\Sigma,\\
\mathrm{and~the~diagram~is~efficient}
\end{gathered}
\right\} .\label{eq:ssql in free groups with efficient surfaces}
\end{equation}
Diagrams obtaining the infimum for a given $w$ are called \emph{$\ssql$-extremal
surfaces }(of $w$). We also define $\ssql\left(1_{\F}\right)=-1$
(thinking of the disk as an admissible surface of degree one for $w=1$).
\end{defn}

As in Definition \ref{def:efficient map with surfaces}, the diagram
being efficient means here that $f$ is $\pi_{1}$-injective and that
for every distinct $o_{1},o_{2}\in\rho^{-1}\left(s_{0}\right)$ and
every path $p$ in $\Sigma$ from $b\left(o_{1}\right)$ to $b\left(o_{2}\right)$,
the image of $p$ in $\Omega$ is not null-homotopic. Note that this
rules out the trivial solutions of the Möbius band and annulus mentioned
above. There are two differences between the definition of $\scl$
in Corollary \ref{cor:scl in free groups - restricted definition}
and the definition of $\ssql$: the surface $\Sigma$ is not necessarily
orientable, and, more importantly, the map $b$ is not necessarily
orientation preserving. (Of course, if $\Sigma$ is non-orientable,
it does not even make sense to demand that $b$ be orientation preserving.)
\begin{rem}
\label{rem:def of ssql}We make the following remarks concerning Definition
\ref{def:ssql}.
\begin{enumerate}
\item \label{enu:ssql with orientable surfaces}For every diagram as in
\eqref{eq:ssql in free groups with efficient surfaces}, an arbitrary
finite degree cover of $\Sigma$ (with a suitable cover of $P$) yields
another valid diagram with the same quotient $\frac{-\chi\left(\Sigma'\right)}{\deg\left(\rho'\right)}$.
In particular, every non-orientable $\Sigma$ admits a $2$-cover
which is orientable, so one can restrict to diagrams with orientable
surfaces. The important point, therefore, is that $b$ need not be
orientation preserving. 
\item Recall that in the case of $\scl$, one can restrict the definition
to surfaces with one boundary component (this was, in fact, the original
definition), and that we suspect the analogous restriction is possible
in the definition of $\sp$ (Conjecture \ref{conj:spi and spi tilde}).
It is plausible that the same is true for $\ssql.$ Notice, though,
that if $w\notin\left[\F,\F\right]$, then surfaces admissible for
$w$ with one boundary component must be non-orientable.
\item \label{enu:ssql without pi_1-inj issue}We already noted above some
trivial surfaces -- the Möbius band and annulus -- which are ruled
out by the condition on the images through $b$ of the fiber $\rho^{-1}\left(s_{0}\right)$
(this is part of the definition of efficiency). We do not know what
happens if we omit the other part of efficiency, namely, the constraint
that $f$ be $\pi_{1}$-injective: whether the infimum gets different
values. However, the proofs of some of the properties of $\ssql$
we state do depend on this additional constraint.\\
In this regard, especially troubling is the fact that in the definitions
of $\mathrm{sql}$ and $\sqlh$, the map $f$ need not be $\pi_{1}$-injective
(it trivially satisfy the other efficiency condition), and we are
not even sure that $\ssql\le\sqlh-1$. If this is not the case, the
name \emph{stable square }length\emph{ }may not be justified. Moreover,
in light of \eqref{eq:sql hat bound for O(N)} (which should be tight
for generic words), Conjecture \ref{conj:ssql and its role for U,O and Sp}
with ${\cal I}_{O}$ strongly suggests that this inequality should
hold.
\item Definition \ref{def:ssql} can be extended to elements of every group:
one just needs to replace $\Omega$ with a suitable space, as in \eqref{eq:scl a la Calegari}. 
\end{enumerate}
\end{rem}

Recall Corollary \ref{cor:scl in free groups - def with WH graphs}
and the paragraph preceding it. As with $\scl$, because every surface
$f\colon\Sigma\to\Omega$ in \eqref{eq:ssql in free groups with efficient surfaces}
is, in particular, incompressible, it can be constructed as a fatgraph
over $\F$, possibly with twists. By Lemma \ref{lem:fatgraphs with twists and whitehead},
we get the following. (Recall the notation $\Gamma_{w}$ from Definition
\ref{def:Gamma_w}.)
\begin{cor}
\label{cor:ssql - def with WH graphs}Let $1\ne w\in\F$. Then 
\begin{equation}
\ssql\left(w\right)=\inf\left\{ \frac{-\chi\left(\Gamma\right)}{\deg\left(\rho\right)}\,\middle|\,\begin{gathered}\xymatrix{P\ar@{->>}[d]^{\rho~~}\ar[r]^{b} & \Gamma\ar[d]^{f}\\
\Gamma_{w}\ar[r]^{\eta_{w}} & \Omega
}
\end{gathered}
\begin{gathered}\mathrm{s.t.~\rho~\mathrm{is~a~covering~map,}}\\
b~\mathrm{is~an~immersion~of~core~graphs,}\\
\mathrm{the~diagram~is~efficient,~and}\\
\wh_{b}\left(v\right)~\mathrm{is~a~cycle}~\forall v\in V\left(\Gamma\right)
\end{gathered}
\right\} .\label{eq:ssql with WH}
\end{equation}
As in Definition \ref{def:efficient map with surfaces}, \emph{efficient}
here means that $f$ is $\pi_{1}$-injective and that for every $v\in V\left(\Gamma_{w}\right)$,
every two distinct points $s_{1},s_{2}\in\rho^{-1}\left(v\right)$
and every path $p$ from $b\left(s_{1}\right)$ to $b\left(s_{2}\right)$,
the image of $p$ in $\Omega$ is not null-homotopic.
\end{cor}

In \eqref{eq:ssql with WH}, unlike in \eqref{eq:scl in free groups with WH},
the map $b$ need not be orientation preserving -- this is the only
difference between $\ssql$ and $\scl$ when defined using fatgraphs
and Whitehead graphs. The following are some facts about $\ssql$,
mostly due to Wilton.
\begin{fact}
\label{fact:basic properties of ssql}
\begin{enumerate}
\item \label{enu:sp < ssql < 2scl}For every $w\in\F$, $\sp\left(w\right)\le\ssql\left(w\right)\le2\scl\left(w\right)$.
\item \label{enu:ssql rational}For every $w\in\F$, $\ssql\left(w\right)\in\mathbb{Q}\cup\left\{ \infty\right\} $.
\item \label{enu:ssql finite for imprimitive}$\ssql\left(w\right)=\infty$
if and only if $w$ is primitive.
\item \label{enu:ssql on powers and gap on non-powers}If $w\ne1$ is a
proper power, then $\ssql\left(w\right)=0$. If $w\ne1$ is not a
proper power, then $\ssql\left(w\right)\ge1$.
\item \label{enu:ssql AutF-invariant}$\ssql$ is $\mathrm{Aut}\F$-invariant.
\end{enumerate}
\end{fact}

\begin{proof}
Item \ref{enu:sp < ssql < 2scl} follows immediately from \cite[Thm.~A and Lem.~10.9]{wilton2024rational},
but also from comparing Corollaries \ref{cor:scl in free groups - def with WH graphs}
and \ref{cor:ssql - def with WH graphs} and Proposition \ref{prop:spi in terms of WH}.
Item \ref{enu:ssql rational} is a special case of \cite[Thm.~C]{wilton2022rationality}.
Item \ref{enu:ssql finite for imprimitive} is a deep theorem that
follows from a combination of \cite{wilton2018essential} and \cite{wilton2022rationality}.
If $w=u^{n}$ for some $1\ne u\in\F$ and $n\ge2$, it is easy to
construct a valid diagram for $w$ as in \eqref{eq:ssql in free groups with efficient surfaces}
with the surface $\Sigma$ being an annulus with each boundary component
mapping to $w$. For $\Sigma$ an annulus, $\chi\left(\Sigma\right)=0$,
which shows that $\ssql\left(w\right)=0$. If $w\ne1$ is not a proper
power then $\ssql\left(w\right)\ge\sp\left(w\right)\ge1$ by Item
\ref{enu:sp < ssql < 2scl} and Theorem \ref{thm:HW-LW}, proving
Item \ref{enu:ssql on powers and gap on non-powers}. For any homotopy
equivalence $h\colon\Omega\to\Omega$ and a diagram as in \eqref{eq:ssql in free groups with efficient surfaces},
one can obtain another valid diagram by replacing $f$ by $h\circ f$
and $w$ by $h\circ w$. Because every automorphism of $\F$ is realizable
by a homotopy equivalence of $\Omega$, Item \ref{enu:ssql AutF-invariant}
follows. (Item \ref{enu:ssql AutF-invariant} follows also from \cite[Prop.~8.3]{wilton2024rational}.)
\end{proof}
Wilton describes in \cite{wilton2022rationality} a concrete algorithm
to compute $\ssql$. Unlike the algorithm for computing $\scl$ (implemented
as ``scallop'' \cite[P.~106]{calegari2009scl}), Wilton's algorithm
is not efficient and can probably run in practice only for relatively
short words. If a word $w$ contains every letter at most three times,
then in every diagram as in \eqref{eq:ssql with WH}, the map $f$
must be an immersion, and $\ssql$ can be computed more easily, in
an algorithm similar to Calegari's. We give two concrete examples.
\begin{example}
\label{exa:ssql of aabbaB}Consider the word $w=a^{2}b^{2}ab^{-1}\in\F_{2}$.
As $w$ is not a power, we have $\ssql\left(w\right)\ge1$. Figure
\ref{fig:ssql of aabbaB} shows a surface $\Sigma$ with Euler characteristic
$-2$ and one boundary component reading $w^{2}$, and which satisfies
the assumptions of Definition \ref{def:ssql}. As $\frac{-\chi\left(\Sigma\right)}{\deg\left(\rho\right)}=\frac{2}{2}=1$,
we conclude that $\ssql\left(w\right)=1$. Moreover, in this case
it can be shown that every $\ssql$-extremal surface of $w$ is a
cover of this degree-2 surface. \vspace{-22bp}
\end{example}

\begin{center}
\begin{figure}
\includegraphics{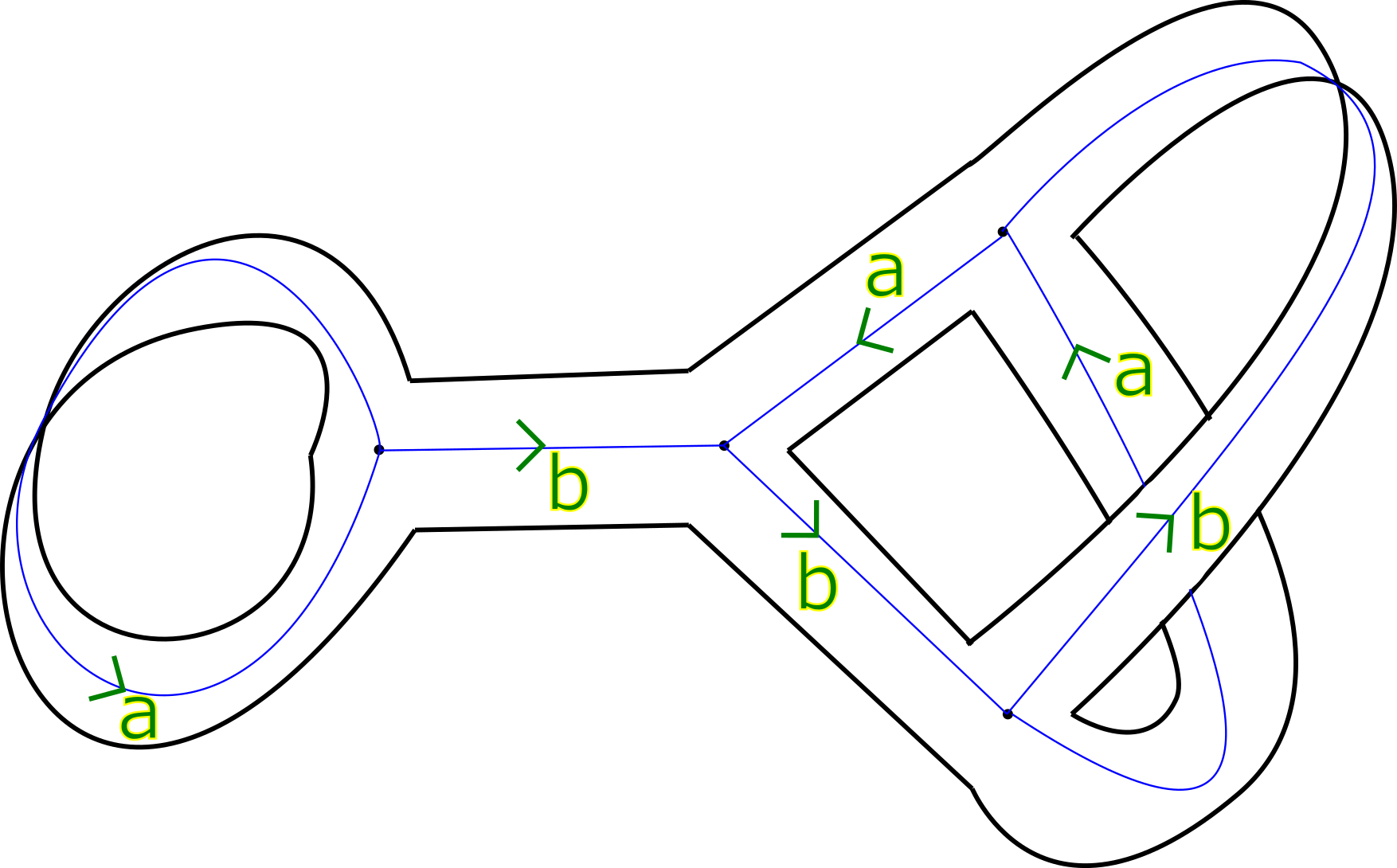}
\centering{}\caption{\label{fig:ssql of aabbaB}A degree-2 $\protect\ssql$-extremal surface
for $w=a^{2}b^{2}ab^{-1},$ showing that $\protect\ssql\left(w\right)=1$.}
\end{figure}
\par\end{center}
\begin{example}
\label{exa:ssql of aBcbbaCac}As another example, consider the word
$w=ab^{-1}cb^{2}ac^{-1}ac$ from Table \ref{tab:values on selected words}.
As mentioned in that table, $\ssql\left(w\right)=2.25$. To see this,
note that the ordinary Whitehead graph of $w$ is bipartite without
double edges. This means that in any diagram as in \eqref{eq:ssql with WH},
the local Whitehead graphs $\wh_{b}\left(v\right)$ must be cycles
of even length $\ge4$, namely, that every vertex in $\Gamma$ has
even degree at least $4$. Now, $-\chi\left(\Gamma\right)=\frac{1}{2}\sum_{v\in V\left(\Gamma\right)}\left(\deg\left(v\right)-2\right)$,
and $\deg\left(\rho\right)=\frac{\sum_{v\in V\left(\Gamma\right)}\deg\left(v\right)}{\left|w\right|}$.
So up to a constant, $\frac{-\chi\left(\Gamma\right)}{\deg\left(\rho\right)}=\frac{\left|w\right|}{2}\cdot\frac{\sum_{v\in V\left(\Gamma\right)}\left(\deg\left(v\right)-2\right)}{\sum_{v\in V\left(\Gamma\right)}\deg\left(v\right)}$
is a weighted average of the numbers $\frac{\deg(v)-2}{\deg(v)}$,
each of which is at least $\frac{1}{2}$, a value obtained only when
the degree is 4. As there exist 4-regular diagrams, where each of
the nine possible 4-cycles appears the same number of times, we conclude
the the infimum is $\frac{|w|}{2}\cdot\frac{1}{2}=\frac{9}{4}$. We
thank Niv Levhari for producing this example. We compute other stable
invariants of this word in Example \ref{exa:spm of aBcbbaCac}.
\end{example}

\subsection{Conjecture \ref{conj:ssql and its role for U,O and Sp} for $\protect\U\left(\bullet\right)$\label{subsec:conj for Arbitrary-stable-irreducible of U}}

As explained in $\S\S$\ref{subsec:stable-representation of U}, the
main result of \cite{MPunitary} is a formula for $\mathbb{E}_{w}\left[\zeta_{\mu}\right]$
in terms of orientable compact surfaces with boundary components corresponding
to $w^{\mu_{1}},\ldots,w^{\mu_{k}}$. These are more general surfaces
than those appearing in Definition \ref{def:ssql} (although they
must be orientable -- as in Remark \ref{rem:def of ssql}\eqref{enu:ssql with orientable surfaces}):
the surfaces considered in \cite{MPunitary} do not necessarily give
rise to efficient diagrams in \eqref{eq:ssql in free groups with efficient surfaces}.
However, surfaces which do give rise to an efficient diagram have
their coefficient $c_{\Sigma,f}$ from \eqref{eq:MP19} equal to $1$
(this follows from \cite[Lem.~5.1]{MPunitary}). For this reason,
for the similarities with Theorem \ref{thm:MP scl for U(N)} (and
more concretely Corollary \ref{cor:same measures =00003D=00003D> same mifkad of extremal scl surfaces})
and with Theorem \ref{thm:result on wreath products and modulo-m spi},
and for other considerations we omit here, we conjecture the following
elaborated version of Conjecture \ref{conj:ssql and its role for U,O and Sp}
regarding $\U\left(\bullet\right)$.
\begin{conjecture}
\label{conj:arbitrary stable chars of U}For $1\ne w\in\F$ and an
integer sequence $\mu=\left(\mu_{1},\ldots,\mu_{k}\right)$ (say,
with no zeros), let $d_{w,\mu}$ denote the number of \textbf{\emph{orientable}}
$\ssql$-extremal surfaces from Definition \ref{def:ssql} with oriented
boundary components corresponding to $w^{\mu_{1}},\ldots,w^{\mu_{k}}$,
up to the equivalence defined in \cite[Def.~1.3]{MPunitary}. Then
for every pair of partitions $\lambda=(\lambda^{+},\lambda^{-})$,
\[
\mathbb{E}_{w}\left[\xi^{\lambda\left[N\right]}\right]=N^{-\left|\lambda\right|\ssql\left(w\right)}\left(\left(\sum_{\mu\colon\left|\mu\right|=\left|\lambda\right|}a_{\lambda,\mu}d_{w,\mu}\right)+O\left(N^{-1}\right)\right),
\]
where $a_{\lambda,\mu}$ is the coefficient of $\zeta_{\mu}$ in the
formula for $\xi^{\lambda\left[N\right]}$ as in \eqref{eq:xi^lambda as a sum over zeta_mu - example}. 
\end{conjecture}

For example, for $\lambda=\left(2,-1,-1\right)$, the conjecture is
that $\mathbb{E}_{w}\left[\xi^{\lambda\left[N\right]}\right]=O\left(N^{-4\ssql\left(w\right)}\right)$,
and the coefficient of $N^{-4\ssql\left(w\right)}$ is 
\[
\frac{1}{4}\left(d_{w,\left(1,1,-1,-1\right)}+d_{w,\left(2,-1,-1\right)}-d_{w,\left(1,1,-2\right)}-d_{w,\left(2,-2\right)}\right).
\]
Recall that $\dim\xi^{\lambda[N]}=\Theta(N^{|\lambda|})$, so $N^{-|\lambda|\ssql(w)}=\Theta((\dim\xi^{\lambda[N]})^{-\ssql(w)}$.
It can be easily shown that when $\lambda^{+}$ and $\lambda^{-}$
both have at most one row, then the coefficients $a_{\lambda,\mu}$
are positive for every $\mu$ with $\left|\mu\right|=\left|\lambda\right|$.
These irreducible characters, therefore, show how Conjecture \ref{conj:arbitrary stable chars of U}
yields the case of ${\cal I}_{U}$ in Conjecture \ref{conj:ssql and its role for U,O and Sp}.
We remark that Conjecture \eqref{conj:arbitrary stable chars of U}
is true for words which are proper powers, in which case $\ssql(w)=0$:
this follows from {[}Ibid., Cor.~1.13{]}.
\begin{rem}
A significant, highly non-trivial step towards Conjecture \ref{conj:arbitrary stable chars of U}
was obtained by Brodsky in \cite{brodsky2022word}. He shows {[}Ibid.,~Cor.~1.10{]}
that for any $\lambda\ne\left((1),\emptyset\right),\left(\emptyset,(1)\right),\left(\emptyset,\emptyset\right)$
and every non-power $w$, $\mathbb{E}_{w}\left[\xi^{\lambda\left[N\right]}\right]=O\left(N^{-\pi\left(w\right)}\right)$.
For $\lambda=\left((1),(1)\right)$ he gets even closer to Conjecture
\ref{conj:arbitrary stable chars of U} and proves that $\mathbb{E}_{w}\left[\xi^{\lambda\left[N\right]}\right]=O((\dim\xi^{\lambda\left[N\right]})^{1-\pi\left(w\right)})$
{[}Ibid., Cor.~1.12{]}. Recall that the conjecture is trivial for
$\lambda=\left(\emptyset,\emptyset\right)$ and known for $\lambda=\left((1),\emptyset\right),(\emptyset,(1))$
by \eqref{eq:bound with cl}.
\end{rem}

\subsection{Stable irreducible characters of $\protect\O\left(\bullet\right)$
and $\protect\Sp\left(\bullet\right)$\label{subsec:Stable-irreducible-characters O and Sp}}

There is a duality between $\O$ and $\Sp$ in their (stable) representation
theory (e.g., \cite[\S4.3.4]{sam_snowden_2015}), and in word measures
(e.g., \cite[Thm.~1.2]{MPorthsymp}). In particular, it follows that
our conjectures for both families are equivalent, and for simplicity
we restrict our discussion here to the orthogonal group.

The irreducible complex representations of $\O\left(N\right)$ are
in bijection with partitions $\nu$ such that the sum of lengths of
the first two columns is at most $N$ \cite[\S19.5]{fulton2013representation}.
Stable irreducible characters are in bijection with partitions, where
the partition $\nu$ gives rise to the stable irreducible character
$\psi^{\nu[\bullet]}=\{\psi^{\nu\left[N\right]}\}_{N\ge\nu'_{1}+\nu'_{2}}$
(recall that $\nu'$ denotes the transpose of $\nu$). See \cite[\S4]{sam_snowden_2015}.

We do not know of a reference with a concrete formula, but it should
be possible to write every stable irreducible character of $\O\left(\bullet\right)$
as a linear combination of power sums of the form $\zeta_{\mu}$ with
$\mu$ a partition. For example, it is possible to obtain the following
formulas for $\left|\nu\right|\le3$. Of course, $\psi^{\emptyset\left[N\right]}\equiv\zeta_{\emptyset}=1$
is the trivial character. 
\begin{eqnarray}
\psi^{\left(1\right)\left[N\right]} & = & \zeta_{\left(1\right)}\nonumber \\
\psi^{\left(2\right)\left[N\right]} & = & \frac{1}{2}\zeta_{\left(1,1\right)}+\frac{1}{2}\zeta_{\left(2\right)}-1\nonumber \\
\psi^{\left(1,1\right)\left[N\right]} & = & \frac{1}{2}\zeta_{\left(1,1\right)}-\frac{1}{2}\zeta_{\left(2\right)}\nonumber \\
\psi^{\left(3\right)\left[N\right]} & = & \frac{1}{6}\zeta_{\left(1,1,1\right)}+\frac{1}{2}\zeta_{\left(2,1\right)}+\frac{1}{3}\zeta_{\left(3\right)}-\zeta_{\left(1\right)}\nonumber \\
\psi^{\left(2,1\right)\left[N\right]} & = & \frac{1}{3}\zeta_{\left(1,1,1\right)}-\frac{1}{3}\zeta_{\left(3\right)}-\zeta_{\left(1\right)}\nonumber \\
\psi^{\left(1,1,1\right)\left[N\right]} & = & \frac{1}{6}\zeta_{\left(1,1,1\right)}-\frac{1}{2}\zeta_{\left(2,1\right)}+\frac{1}{3}\zeta_{\left(3\right)}\label{eq:examples for formulas in O}
\end{eqnarray}

Let $\mu=\left(\mu_{1},\ldots,\mu_{k}\right)$ be a partition. The
main result of \cite{MPorthsymp} is a formula for $\mathbb{E}_{w}\left[\zeta_{\mu}\right]$
in the group $\O\left(N\right)$, which is given by equivalence classes
of $\left(\Sigma,f\right)$ where $\Sigma$ is a compact surface (not
necessarily orientable), $f\colon\Sigma\to\Omega$ is a map to the
bouquet, and the boundary components of $\Sigma$ are mapped to $w^{\mu_{1}},\ldots,w^{\mu_{k}}$
{[}Ibid, $\S\S$1.1{]}. Similar in spirit to the formula \eqref{eq:MP19}
for $\U\left(N\right)$, the formula in the current case is 
\[
\mathbb{E}_{w}\left[\zeta_{\mu}\right]=\sum_{\left[\left(\Sigma,f\right)\right]}c_{\Sigma,f}\left(N-1\right)^{\chi\left(\Sigma\right)},
\]
where $c_{\Sigma,f}$ is an integer which is equal to the $L^{2}$-Euler
characteristic of some natural stabilizer of $f$. These surfaces
do not necessarily give rise to efficient diagrams as in \eqref{eq:ssql in free groups with efficient surfaces},
but surfaces that do give rise to efficient diagrams satisfy that
$c_{\Sigma,f}=1$.
\begin{conjecture}
\label{conj:stable chars of O}For $1\ne w\in\F$ and a partition
$\mu=\left(\mu_{1},\ldots,\mu_{k}\right)$, let $e_{w,\mu}$ denote
the number of $\ssql$-extremal surfaces (orientable or otherwise)
from Definition \ref{def:ssql} with (non-oriented) boundary components
corresponding to $w^{\mu_{1}},\ldots,w^{\mu_{k}}$, up to the equivalence
defined in \cite[Def.~1.3]{MPorthsymp}. Then for every partition
$\nu$,
\[
\mathbb{E}_{w}\left[\psi^{\nu\left[N\right]}\right]=N^{-\left|\nu\right|\ssql\left(w\right)}\left(\left(\sum_{\mu\vdash\left|\nu\right|}b_{\nu,\mu}e_{w,\mu}\right)+O\left(N^{-1}\right)\right),
\]
where $b_{\nu,\mu}$ is the coefficient of $\zeta_{\mu}$ in the formula
for $\psi^{\nu\left[N\right]}$. 
\end{conjecture}

We expect that Conjecture \ref{conj:stable chars of O} should yield
Conjecture \ref{conj:ssql and its role for U,O and Sp} in the case
of $\O(\bullet)$ (and of $\Sp(\bullet)$). This depends on having
formulas for $\psi^{\nu[\bullet]}$ in terms of the power sums of
the form $\zeta_{\mu}$ which behave nicely as in the cases of $U(n)$
and $\s_{n}$. 
\begin{example}
\label{exa:measure on O by aabbaB}Consider the word $w=a^{2}b^{2}ab^{-1}$
from Example \ref{exa:ssql of aabbaB} and Figure \ref{fig:ssql of aabbaB}.
Using the notation from Conjecture \ref{conj:stable chars of O},
the fact that the $\ssql$-extremal surface in Figure \ref{fig:ssql of aabbaB}
is the only degree-2 $\ssql$-extremal surface for $w$ and every
other one is a cover of it, yields that $e_{w,\left(2\right)}=2$
(there are two possible starting points for $w^{2}$ in that extremal
surface -- see \cite[Def.~1.3]{MPorthsymp}), but $e_{w,\left(1,1\right)}=0$.
With the formulas in \eqref{eq:examples for formulas in O}, Conjecture
\ref{conj:stable chars of O} says that $\mathbb{E}_{w}\left[\psi^{\left(2\right)\left[N\right]}\right]=N^{-2}\left(1+O\left(N^{-1}\right)\right)$
and $\mathbb{E}_{w}\left[\psi^{\left(1,1\right)\left[N\right]}\right]=N^{-2}\left(-1+O\left(N^{-1}\right)\right)$.
And indeed, we verified these special cases by a computer\footnote{Using the formula \cite[Thm.~3.4]{MPorthsymp} and \cite[\S3]{CS},
Magee and the first author wrote a computer code to compute the rational
function which is equal to $\mathbb{E}_{w}\left[\zeta_{\mu}\right]$
in orthogonal groups, for short words and small partitions $\mu$.}. In fact, the computer yields that $\mathbb{E}_{w}[\zeta_{\left(2\right)}]=1+\frac{2}{N^{2}}+\frac{2}{N^{3}}-\frac{18}{N^{4}}+\ldots$
and $\mathbb{E}_{w}[\zeta_{\left(1,1\right)}]=1-\frac{20}{N^{4}}+\ldots$.
The leading term of $1$ in the former corresponds to the trivial
Möbius band and in the latter to the trivial annulus.
\end{example}

\begin{rem}
\label{rem:same formulas}There seems to be a close connection between
the formulas we have, or expect, for the largest possible term in
the expected value of stable polynomial irreducible characters in
different groups. Let $\mu\vdash d$ be a partition. From \eqref{eq:formula for stable poly irreps of U}
and \eqref{eq:leading term of E_w=00005BXi=00005D in U} it follows
that in $\U\left(\bullet\right)$,
\[
\mathbb{E}_{w}\left[\xi^{\mu\left[N\right]}\right]=N^{-d\cdot2\scl\left(w\right)}\left(\frac{1}{d!}\sum_{\sigma\in S_{d}}\chi^{\mu}\left(\sigma\right)d_{w,\sigma}+O\left(N^{-1}\right)\right).
\]
Proposition \ref{prop:detailed conjectural picture for S} shows how
Conjecture \ref{conj:C^alg from cycle and non-efficient} yields that
in $S_{\bullet}$
\[
\mathbb{E}_{w}\left[\chi^{\mu\left[N\right]}\right]=N^{-d\cdot\sp\left(w\right)}\left(\frac{1}{d!}\sum_{\sigma\in S_{d}}\chi^{\mu}\left(\sigma\right)\mathrm{crit}_{w,\sigma}+O\left(N^{-1}\right)\right).
\]
Although we did not find a precise formula in the literature for $\mathbb{E}_{w}\left[\psi^{\mu\left[N\right]}\right]$
in terms of the $\zeta_{\nu}$'s, it seems from the examples in \eqref{eq:examples for formulas in O}
that Conjecture \ref{conj:stable chars of O} translates to a parallel
formula in $\O\left(\bullet\right)$, namely, that
\[
\mathbb{E}_{w}\left[\psi^{\mu\left[N\right]}\right]=N^{-d\cdot\ssql\left(w\right)}\left(\frac{1}{d!}\sum_{\sigma\in S_{d}}\chi^{\mu}\left(\sigma\right)e_{w,\sigma}+O\left(N^{-1}\right)\right).
\]
\end{rem}

\section{Stable mod-$m$ primitivity rank and stable characters of $C_{m}\wr S_{\bullet}$\label{sec:spm + G wr S}}

\subsection{Stable mod-$m$ primitivity rank \label{subsec:spm}}

Recall that the mod-$m$ primitivity rank of a word $w\in\F$, denoted
$\pi^{\left(m\right)}\left(w\right)$, is the smallest rank of a subgroup
$H\le\F$ with $w\in K_{m}\left(H\right)$ -- see \eqref{eq:pi modulo m}.
Here we define a stable version of this invariant, in analogy with
the stable primitivity rank and the other stable invariants mentioned
above. We will see that when $m\ne1$, the definition of $\spm$ can
be relatively ``light'', with fewer conditions than in $\sp^{\left(1\right)}=\sp$
(similarly to the way the definition of $\scl$ in \eqref{eq:scl a la Calegari in F}
is lighter than the definition of $\ssql$).

A key feature of $\pi^{\left(m\right)}\left(w\right)$ is the following
observation. Assume that $w\in H\le\F$ and let $p$ be the oriented
closed path corresponding to $w$ in the Stallings core graph corresponding
to $H$. Then $w\in K_{m}\left(H\right)$ if and only if the number
of times $p$ traverses every edge (counted with signs corresponding
to the orientation of each traverse) is in $m\mathbb{Z}$ \cite[Lem.~3.2]{MPsurfacewords}.
This naturally gives rise to the following stable counterpart of $\pi^{\left(m\right)}$. 
\begin{defn}
\label{def:spm}Let $1\ne w\in\F$ and $1\ne m\in\mathbb{Z}_{\ge0}$.
The \textbf{mod-$m$ stable primitivity rank} of $w$ is 
\begin{equation}
\spm\left(w\right)\defi\inf\left\{ \frac{-\chi\left(\Gamma\right)}{\deg\left(\rho\right)}\,\middle|\,\begin{gathered}\xymatrix{P\ar@{->>}[d]^{\rho~~}\ar[r]^{b} & \Gamma\ar[d]^{f}\\
\Gamma_{w}\ar[r]^{\eta_{w}} & \Omega
}
\end{gathered}
~\begin{gathered}\mathrm{s.t.}~\Gamma~\mathrm{is~a~core~graph~and}~f~\mathrm{an~immersion,}\\
\rho~\mathrm{is~a~finite\text{-}degree~covering~map,}\\
\text{the~diagram~is~efficient,}\\
\text{and}\,\forall e\in E\left(\Gamma\right),~n_{b}\left(e\right)\in m\mathbb{Z}.
\end{gathered}
\right\} ,\label{eq:spm with core graphs}
\end{equation}
where $P$ inherits an orientation from $\Gamma_{w}$ and $n_{b}\left(e\right)$
is the signed sum of times $e\in E\left(\Gamma\right)$ is covered
by $b$: a positive sign if the orientation of the edge in $P$ agrees
with that of $e$, and a negative sign otherwise. We also define $\spm\left(1_{\F}\right)\defi-1$. 
\end{defn}

As in Definition \ref{def:spi}, efficiency here is equivalent to
that $b$ is injective on some (or every) fiber of $\rho$. In practice,
we may demand that $b$ be surjective: otherwise, restricting to $b(P)$
instead of the whole $\Gamma$ does not increase $\frac{-\chi(\Gamma)}{\deg\rho}$
(recall that every connected component of a core graph has non-positive
Euler characteristic). We will see below (following Proposition \ref{prop:spf with algebraic morphisms})
that there is also a unified definition encompassing the $m=1$ case
as well. 
\begin{example}
\label{exa:spm}Let $w=u^{d}$ be a proper power, with $1\ne u\in\F$
a non-power and $d\ge2$. If $2\le m\le d$ then $\spm\left(w\right)=0$.
Indeed, consider $\Gamma=\Gamma_{w}$. There are $d$ vertices in
$\Gamma_{w}$ from which $w$ can be read as a cycle. Choose $m$
of them. Then one can map $m$ distinct copies of $\Gamma_{w}$ into
$\Gamma$, where in each of them the starting point of $w$ is mapped
to a distinct vertex out of the $m$ chosen ones. It is easy to see
that this diagram, of degree $m$, satisfies the conditions in \eqref{eq:spm with core graphs},
so $\spm\left(w\right)\le\frac{0}{m}=0$ and thus $\spm\left(w\right)=0$.

As another example, recall that $\pi^{\left(m\right)}\left(w\right)=\min\left\{ -\chi\left(\Gamma\right)+1\right\} $
where $\Gamma$ is as in \eqref{eq:spm with core graphs}, but the
diagram is of degree one. So $\spm\left(w\right)\le\pi^{\left(m\right)}\left(w\right)-1$
for every $w\in\F$. In particular, if $\pi^{\left(m\right)}\left(w\right)<\infty$,
which happens if and only if $w\in K_{m}\left(\F\right)$, then $\spm\left(w\right)\le\pi^{\left(m\right)}\left(w\right)-1\le\rk\F-1$.
\end{example}

\begin{rem}
\label{rem:connected definition does not work}Recall from Conjecture
\ref{conj:spi and spi tilde} that we suspect that in the definition
of the stable primitivity rank we may restrict to diagrams where $P$
is connected. This is not true in general for $\spm$. For example,
consider the simple power word $w=a^{3}$. As explained in Example
\ref{exa:spm}, $\sp^{\left(2\right)}\left(a^{3}\right)=0$. However,
it is easy to see that if we restrict to diagrams with a connected
cover $P$, then there are no valid diagrams and the value of the
invariant becomes $\infty$. In particular, the analog of Wilton's
Conjecture \ref{conj:spi=00003Dpi-1} for $\spm$ and $\pi^{\left(m\right)}$
is false, so $\spm$ is a genuine new invariant\footnote{This can also be seen from the fact that $\spm$ may obtain non-integer
values -- see Example \ref{exa:spm of aBcbbaCac}.}.
\end{rem}

Before moving on to establishing Theorem \ref{thm:result on wreath products and modulo-m spi}
in the subsequent subsections, we record some properties of the invariants
$\spm$, many of which are parallel to known properties of $\scl$,
$\sp$ and $\ssql$. We start with rationality. For $\sp^{(1)}=\sp$,
this is the content of \cite[Thm.~A]{wilton2022rationality}. The
case $m\ne1$ is significantly simpler and follows the lines of Calegari's
proof in \cite{calegari2009stable} of the rationality of $\scl$.
\begin{thm}
\label{thm:rationality of spm} Let $1\ne m\in\mathbb{Z}_{\ge0}$
and $w\in\F$. Then $\spm(w)\in\mathbb{Q}\cup\{\infty\}$, the value
is computable and if finite, it is attained.
\end{thm}

The last assertion means that the infimum in Definition \ref{def:spm}
is, in fact, a minimum: as long as the set of valid diagrams is not
empty, there are extremal diagrams attaining the infimum.
\begin{proof}
The strategy is to construct diagrams via linear programming. First,
consider a valid diagram as in \eqref{eq:spm with core graphs}. Decorate
the Whitehead graphs of the morphism $b\colon P\to\Gamma$ as follows.
For every $v\in V\left(\Gamma\right)$ and every edge $e$ of $\wh_{b}\left(v\right)$,
let $u\in b^{-1}\left(v\right)$ be the vertex of $P$ to which $e$
corresponds. Orient $e$ according to the orientation of $P$, and
decorate it by the vertex $\rho\left(u\right)\in V\left(\Gamma_{w}\right)$.
The efficiency of the diagram is equivalent to that for every $v\in V(\Gamma)$,
all the edges of $\wh_{b}\left(v\right)$ have distinct decorations.
This is illustrated in Figure \ref{fig:pieces - example with m=00003D3}.
It follows that (the decorated) $\wh_{v}\left(b\right)$ is a subgraph
of (the decorated) $\wh_{\eta_{w}}\left(o\right)$, where $\eta_{w}\colon\Gamma_{w}\to\Omega$
is as in Definition \ref{def:Gamma_w} and $o$ is the sole vertex
of $\Omega$. The fact that $\Gamma$ is a core graph means that $v$
has degree $\ge2$, namely, $\wh_{b}(v)$ has at least two vertices.
The condition $n_{b}(e)\in m\mathbb{Z}$ for every $e\in E(\Gamma)$
translates to that for every $v\in V(\Gamma)$, in every vertex of
the Whitehead graph $\wh_{b}(v)$, the difference between the number
of incoming edges and that of outgoing edges is divisible by $m$
(in particular, if $m=0$, this difference must be zero).\footnote{\label{fn:only connected pieces}We may assume further, without loss
of generality, that each Whitehead graph $\wh_{b}(v)$ is connected:
we already explained why we may assume $b$ is surjective, meaning
that $\wh_{b}(v)$ has no isolated vertices, and if $\wh_{b}(v)$
is not connected, we may split the corresponding vertex of $\Gamma$,
remain with a core graph and decrease $\frac{-\chi(\Gamma)}{\deg(\rho)}$.
However, this observation is not needed in the proof.} The crucial observation here is that it follows that there is only
a finite number of possible decorated Whitehead graphs which may appear
as a Whitehead graph of a diagram for $w$ as in \eqref{eq:spm with core graphs}.

Following Calegari, we say that such a subgraph of the decorated $\wh_{\eta_{w}}\left(o\right)$
satisfying the latter conditions of divisibility by $m$ and having
at least two vertices is a \emph{piece}. Another important property
of the pieces involved in a given diagram as in \eqref{eq:spm with core graphs}
is the following. Assume that one of these pieces has a vertex $s$
corresponding to an outgoing $a$-edge of $\Omega$, namely, that
the corresponding oriented edge $e$ of $\Gamma$ projects to an outgoing
$a$-edge of $\Omega$. Assume further that $s$ is incident to $k$
Whitehead edges labeled by $x_{1},\ldots,x_{k}\in V(\Gamma_{w})$.
We call such a vertex of a piece an $\{x_{1},\ldots,x_{k}\}$-incoming-$a$-vertex.
Then the vertex $t$ in the Whitehead graph of the other endpoint
of $e$ which corresponds to $\overline{e}$ is necessarily a $\{a.x_{1},\ldots,a.x_{k}\}$-outgoing-$a$-vertex.
(Here, for a vertex $x\in V(\Gamma_{w})$ with an outgoing $a$-edge,
$a.x$ denotes the vertex on the other end of that $a$-edge.)

\begin{figure}
\includegraphics[scale=0.6]{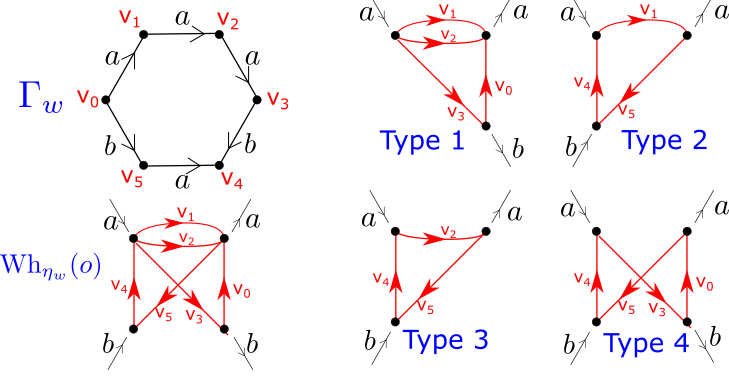}
\centering{}\caption{\label{fig:pieces - example with m=00003D3}Consider $w=a^{3}ba^{-1}b^{-1}$.
The upper left graph is $\Gamma_{w}$, and the bottom left one is
(the decorated) $\protect\wh_{\eta_{w}}(o)$ (the black oriented edge
labeled by $a$ or $b$ at every vertex only marks the corresponding
half-edge of $\Omega$, and is not a genuine part of the Whitehead
graph). The other four graphs depict the four types of pieces which
are valid for this word when $m=3$. Here we restrict to pieces representing
connected Whitehead graphs only (see Footnote \ref{fn:only connected pieces}).}
\end{figure}

Now assume we are given a finite collection of pieces (possibly with
repetitions). If we can glue all the vertices of the pieces in pairs,
we obtain a valid diagram which satisfies all the conditions in \eqref{eq:spm with core graphs}.
Minus the Euler characteristic of the resulting diagram is 
\begin{equation}
\sum_{\mathfrak{p}:~\mathrm{piece~in~the~collection}}\left(\frac{1}{2}\left|V(\mathfrak{p})\right|-1\right),\label{eq:expression for Euler char}
\end{equation}
and the degree is the number of pieces containing an edge labeled
by some arbitrary but fixed $x\in V(\Gamma_{w})$.

Hence, the problem can be turned into a linear programming one. For
every possible type $\mathfrak{g}$ of a piece (namely, for every
subgraph of $\wh_{\eta_{w}}\left(o\right)$ satisfying the $n_{b}$
conditions and having at least two vertices), we denote by $n_{\mathfrak{g}}$
the number of pieces of type $\mathfrak{g}$. These numbers should
satisfy the inequalities $n_{\mathfrak{g}}\ge0$ and the following
equations: For every generator $a$ of $\F$ (equivalently, for every
edge of $\Omega$), and every subset $\{x_{1},\ldots,x_{k}\}\subseteq V(\Gamma)$
of the vertices of $\Gamma$ such that $x_{1},\ldots,x_{k}$ are incident
to an outgoing $a$-edge, the total number of pieces with an $\{x_{1},\ldots,x_{k}\}$-incoming-$a$-vertex
is equal to the total number of pieces with an $\{a.x_{1},\ldots,a.x_{k}\}$-outgoing-$a$-vertex.

Since the equations are homogeneous, we may look for rational solutions
and then multiply by a common denominator. We may thus normalize the
degree by adding an equation saying that the degree of the diagram
is one, and simply minimize the value of \eqref{eq:expression for Euler char}
among all solutions of the above linear equalities and inequalities.
Since all of these were defined over $\mathbb{Q}$, if the solution
space is not empty, there must be a minimizing solution over $\mathbb{Q}$,
which then translates, multiplying by a common denominator, to a solution
over $\mathbb{Z}$ and a valid, extremal diagram. 
\end{proof}
\begin{example}
\label{exa:computing sp3}We illustrate the procedure depicted in
the proof of Theorem \ref{thm:rationality of spm}. Consider the word
$w=a^{3}ba^{-1}b^{-1}$ and assume that $m=3$. The possible pieces
in this case are depicted in Figure \ref{fig:pieces - example with m=00003D3}.
Denote by $n_{i}$ the number of pieces of type $i$ for $i=1,\ldots,4$.
The linear programming problem consists of four inequalities ($n_{i}\ge0$
for all $i$) and six equations, as follows. First, the total number
of pieces with a $\{v_{0},v_{1},v_{2}\}$-outgoing-$a$ vertex is
equal to the total number of pieces with a $\{v_{1},v_{2},v_{3}\}$-incoming-$a$
vertex, which gives the trivial equation $n_{1}=n_{1}$. Likewise,
comparing pieces with $\left\{ v_{0},v_{5}\right\} $-outgoing-$a$
and $\{v_{1},v_{4}\}$-incoming-$a$ gives $n_{4}=n_{2}$, and the
additional two equations corresponding to $a$-edges are $n_{2}=n_{3}$
and $n_{3}=n_{4}$. There is a single equation for $b$-edges: $n_{1}+n_{4}=n_{2}+n_{3}+n_{4}$.
Finally, constructing the degree equation based on the decoration
$v_{1}$, we get $n_{1}+n_{2}=1$. Under these constraints, we ought
to minimize the value $-\chi=\frac{1}{2}(n_{1}+n_{2}+n_{3})+n_{4}$.
In this case, the equations admit a single solution: $n_{1}=\frac{2}{3},n_{2}=n_{3}=n_{4}=\frac{1}{3}$,
giving $-\chi=1$. Hence $\sp^{(3)}(w)=1$. 
\end{example}

The following proposition shows that in the definition of $\spm$,
we may restrict to diagrams where $b$ is \emph{algebraic} without
affecting the outcome. 
\begin{prop}
\label{prop:spf with algebraic morphisms}Let $1\ne w\in\F$ and $1\ne m\in\mathbb{Z}_{\ge0}$.
Then 
\begin{equation}
\spm\left(w\right)=\inf\left\{ \frac{-\chi\left(\Gamma\right)}{\deg\left(\rho\right)}\,\middle|\,\begin{gathered}\xymatrix{P\ar@{->>}[d]^{\rho~~}\ar[r]^{b} & \Gamma\ar[d]^{f}\\
\Gamma_{w}\ar[r]^{\eta_{w}} & \Omega
}
\end{gathered}
~\begin{gathered}\mathrm{s.t.}~\Gamma~\mathrm{is~a~core~graph~and}~f~\mathrm{an~immersion,}\\
\rho~\mathrm{is~a~finite\text{-}degree~covering~map,}\\
b~\mathrm{is~algebraic,the~diagram~efficient,}\\
\text{and}\,\,\forall e\in E\left(\Gamma\right),~n_{b}\left(e\right)\in m\mathbb{Z}.
\end{gathered}
\right\} .\label{eq:spm with algebraic}
\end{equation}
\end{prop}

As in Definition \ref{def:spm}, $P$ inherits the orientation from
$\Gamma_{w}$. Note that we may add in \eqref{eq:spm with algebraic}
the additional condition that $b$ is not an isomorphism on any connected
component of $\Gamma$ (as in Definition \ref{def:spi}) without affecting
the resulting diagram: the condition that $n_{b}(e)\in m\mathbb{Z}$
anyway forbids such components when $m\ne1$. However, adding this
condition yields a unified definition which works for $m=1$ as well.
The proof we give here is combinatorial in style. In \cite[Cor.~4.10]{PSh25}
we prove a generalization of Proposition \ref{prop:spf with algebraic morphisms}
for a generalization of $\spm$, and the argument is different and
much more algebraic in style.
\begin{proof}
Consider a diagram with $b\colon P\to\Gamma$ as in \eqref{eq:spm with core graphs}.
It is enough to show that we can obtain a diagram as in \eqref{eq:spm with algebraic}
with the quotient $\frac{-\chi\left(\Gamma\right)}{\deg\left(\rho\right)}$
at least as small. First, we may assume that $b$ is surjective, for
otherwise, replacing $\Gamma$ with the image $b(P)$, does not increase
$-\chi(\Gamma)$ and does not change $\deg\rho$. Now let $P\stackrel{b_{\mathrm{alg}}}{\to}\Sigma\stackrel{b_{\mathrm{free}}}{\to}\Gamma$
be the algebraic-free decomposition of $b$ (see $\S\S$\ref{subsec:Algerbaic-morphisms}),
and consider the diagram 
\begin{equation}
\begin{gathered}\xymatrix{P\ar@{->>}[d]^{\rho~~}\ar[r]^{b_{\mathrm{alg}}} & \Sigma\ar[d]^{f\circ b_{\mathrm{free}}}\\
\Gamma_{w}\ar[r]^{\eta_{w}} & \Omega
}
\end{gathered}
\label{eq:diagram with algebraic part of b}
\end{equation}
In this algebraic-free decomposition $\Sigma$ is a core graph, $b_{\text{free}}$
is an immersion (hence so is $f\circ b_{\mathrm{free}}$), and the
efficiency of $b_{\text{alg}}$ follows immediately from that of $b$.
It remains to show that the condition about $n_{b}(e)$ is satisfied
in \eqref{eq:diagram with algebraic part of b} and that $-\chi(\Sigma)\le-\chi(\Gamma)$.
The latter follows from \cite[Lem.~5.3]{hanany2020word}. By \cite[\S5]{hanany2020word},
$b_{\text{free}}\colon\Sigma\to\Gamma$ can be obtained by a sequence
of steps 
\[
\Sigma=\Sigma_{0}\stackrel{c_{1}}{\longrightarrow}\Sigma_{1}\stackrel{c_{2}}{\longrightarrow}\Sigma_{2}\stackrel{c_{3}}{\longrightarrow}\ldots\stackrel{c_{k}}{\longrightarrow}\Sigma_{k}=\Gamma,
\]
where for each $j=1,\ldots,k$, $c_{j}$ is obtained by either $\left(i\right)$
merging together two distinct vertices of $\Sigma_{j-1}$, or $\left(ii\right)$
folding (à la Stallings) two edges with exactly one joint endpoint
in $\Sigma_{j-1}$ (a Stallings fold which is a homotopy equivalence).
Consider the sequence of morphisms $b_{j}\colon P\to\Sigma_{j}$,
$j=0,\ldots,k$ defined by $b_{j}=c_{j}\circ\ldots\circ c_{1}\circ b_{\text{alg}}$
(so $b_{0}=b_{\text{alg}}$ and $b_{k}=b$). We claim that $b_{j}$
satisfies the $n_{b_{j}}(e)$-condition if and only if $b_{j-1}$
does. Indeed, if $c_{j}$ is a step of the first kind, this is obvious:
there is a full bijection between the edges of $\Sigma_{j}$ and these
of $\Sigma_{j-1}$ with the same numbers $n_{b_{j}}(e)=n_{b_{j-1}}(e)$.
\begin{figure}
\includegraphics[scale=0.5]{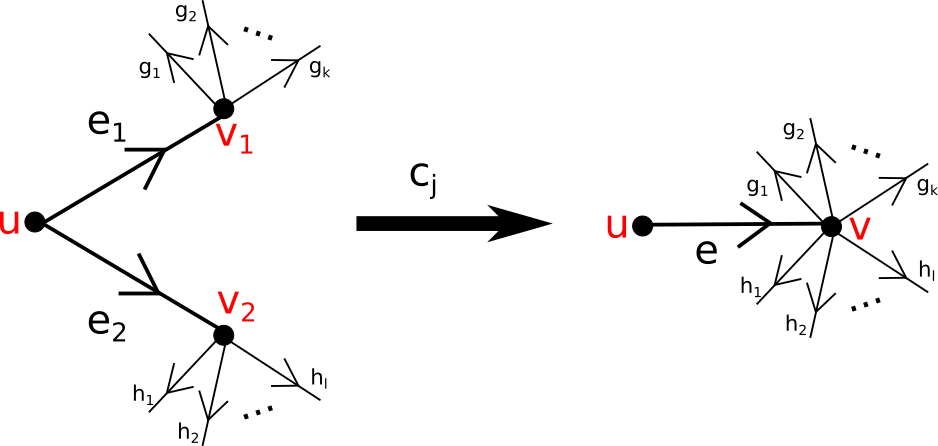}
\centering{}\caption{\label{fig:folding for proof of spm}A homotopy-equivalent folding
step $c_{j}\colon\Sigma_{j-1}\to\Sigma_{j}$}
\end{figure}

Now assume that $c_{j}$ is a step of the second kind, namely, a homotopy-equivalent
fold where we merge $e_{1},e_{2}\in E(\Sigma_{j-1})$ with $u:=\iota(e_{1})=\iota(e_{2})$
and $v_{1}:=\iota(\overline{e_{1}})\ne v_{2}:=\iota(\overline{e_{2}})$
into $e\in E(\Sigma_{j})$ with $\iota(e)=u$ and $v:=\iota(\overline{e})$
(see Figure \ref{fig:folding for proof of spm}). Assume first that
$n_{b_{j-1}}(e)\in m\mathbb{Z}$ for all $e\in E(\Sigma_{j-1})$.
For every $e'\in E(\Sigma_{j})$ with $e'\ne e,\overline{e}$, $n_{b_{j}}(e')=n_{b_{j-1}}(e')\in m\mathbb{Z}$.
As for $e$ itself, $n_{b_{j}}(e)=n_{b_{j-1}}(e_{1})+n_{b_{j-1}}(e_{2})\in m\mathbb{Z}$.
Conversely, assume the condition holds for $b_{j}$. Again, it is
enough to show the condition holds for $e_{1},e_{2}\in E(\Sigma_{j-1})$.
As $v_{1}\ne v_{2}$, at least one of them is different than $u$,
say $v_{1}$. Because $P$ is a collection of cycles, for every vertex
of $\Sigma_{j-1}$, and in particular $v_{1}$, we have 
\[
\sum_{e\colon\iota(e)=v_{1}}n_{b_{j-1}}(e)=0.
\]
Using the notation from Figure \ref{fig:folding for proof of spm},
$n_{b_{j-1}}(e_{1})=-\sum_{j=1}^{k}n_{b_{j-1}}(g_{j})\in m\mathbb{Z}$
(we used that $v_{1}\ne u,v_{2}$ so none of $g_{1},\ldots,g_{k}$
is $e_{2}$ or $\overline{e_{2}}$). As a consequence, $n_{b_{j-1}}(e_{2})=n_{b_{j}}(e)-n_{b_{j-1}}(e_{2})\in m\mathbb{Z}$.
\end{proof}
The following proposition lists additional properties of $\spm$.
It uses our convention that $\sp^{\left(1\right)}$ is a synonym of
$\sp$. We want to have all the ``abstract'' properties we know
of $\spm$ listed in this $\S\S$\ref{subsec:spm}, but one of the
arguments is deferred to $\S\S$\ref{subsec:spm-main-proof}, after
the proof of Theorem \ref{thm:result on wreath products and modulo-m spi}
on which it relies.
\begin{prop}
\label{prop:properties of spm}Let $m,k\in\mathbb{Z}_{\ge0}$ and
let $w\in\F$. Then, 
\begin{enumerate}
\item \label{enu:spm AutF-invariant}$\spm$ is $\mathrm{Aut}\F$-invariant.
\item \label{enu:m|k then spm <=00003D spk}If $m\mid k$ then $\spm\left(w\right)\le\sp^{\left(k\right)}\left(w\right)$.
\item \label{enu:sp0<=00003D2scl}$\sp^{\left(0\right)}\left(w\right)\le2\scl\left(w\right)$.
\item \label{enu:sp<=00003Dspm<=00003Dsp0<=00003D2scl}In particular, $\sp\left(w\right)\le\spm\left(w\right)\le\sp^{\left(0\right)}\left(w\right)\le2\scl\left(w\right)$. 
\item \label{enu:sp2<=00003Dssql}$\sp^{\left(2\right)}\left(w\right)\le\ssql\left(w\right)$.
\item \label{enu:spm gap }$\spm$ admits a gap in $\left(0,1\right)$:
its image avoids this interval.
\item \label{enu:spm outside =00005BF,F=00005D}Let $w\notin[\F,\F]$. Then
$\sp^{(0)}(w)=\infty$. If $m\ne0$, and when $w$ is written in some
basis $S$ of $\F$ every $s\in S^{\pm1}$ appears less than $m$
times, then $\spm(w)=\infty$. 
\item \label{enu:spm inside =00005BF,F=00005D}Let $w\in[\F,\F]$. Then
$\sp^{(m)}(w)\le\rk\F-1$ for all.
\end{enumerate}
\end{prop}

\begin{proof}
\eqref{enu:spm AutF-invariant} The proof of Claim \ref{claim:sp is Aut invariant}
applies here too (ignoring algebraicity), where the additional property
that $n_{b}\left(e\right)\in m\mathbb{Z}$ for all $e\in E\left(\Gamma\right)$
is equivalent to the group-theoretic property that $u_{1}w^{n_{1}}u_{1}^{-1}\cdots u_{m}w^{n_{m}}u_{m}^{-1}\in K_{m}\left(H\right)$. 

\eqref{enu:m|k then spm <=00003D spk} When $m=1$ and $k\ne1$, this
follows from the fact that every diagram in \eqref{eq:spm with algebraic}
is valid also under the conditions in \eqref{eq:def of spi with core graphs}.
If $m\ne1$, this is obvious from Definition \ref{def:spm}, as $m\mid k$
and $n_{b}\left(e\right)\in k\mathbb{Z}$ imply $n_{b}\left(e\right)\in m\mathbb{Z}$. 

\eqref{enu:sp0<=00003D2scl} This is immediate from the fact that
every valid diagram in Corollary \ref{cor:scl in free groups - def with WH graphs}
is also valid for $\sp^{\left(0\right)}$ in Definition \ref{def:spm}.

\eqref{enu:sp<=00003Dspm<=00003Dsp0<=00003D2scl} This follows from
Items \ref{enu:m|k then spm <=00003D spk} and \ref{enu:sp0<=00003D2scl}.

\eqref{enu:sp2<=00003Dssql} This is immediate as every valid diagram
in Corollary \ref{cor:ssql - def with WH graphs} is valid for $\sp^{\left(2\right)}$
in Definition \ref{def:spm}. 

\eqref{enu:spm gap } That every $w\ne1$ which is not a proper power
satisfies $\spm(w)\ge1$ follows from Item \ref{enu:sp<=00003Dspm<=00003Dsp0<=00003D2scl}
and Theorem \ref{thm:HW-LW}. As $\spm(1_{\F})=-1$, it remains to
show that $\spm$ of proper power avoids the interval $\left(0,1\right)$.
This is the content of Theorem \ref{thm:spm avoids (0,1)}.

\eqref{enu:spm outside =00005BF,F=00005D} Assume first that $m=0$.
Assume towards contradiction that there exists a valid diagram as
in \eqref{eq:spm with core graphs}. So $n_{b}(e)=0$ for all $e\in E(\Gamma)$.
Let $S$ be the basis of $\F$ whose elements are in bijection with
the edges of $\Omega$. For every $s\in S$ we have 
\[
\#_{s}(w)-\#_{s^{-1}}(w)=\frac{1}{\deg\rho}\left(\#_{s}(P)-\#_{s^{-1}}(P)\right)=\frac{1}{\deg\rho}\sum_{e\in E(\Gamma)\colon f(e)=s}n_{b}(e)=0,
\]
where $\#_{s}(w)$ ($\#_{s^{-1}}(w)$, respectively) is the number
of times $s$ ($s^{-1}$, respectively), appears in $w$, and $\#_{s}(P)$
($\#_{s^{-1}}(P)$, respectively) is the number of times $s$ appears
in $P$ with positive (negative, respectively) orientation. But this
contradicts our assumption that $w\notin[\F,\F]$.

Now assume that $m\ne0$ and every $s\in S^{\pm1}$ appears less than
$m$ times in $w$ for some basis $S$ of $\F$. By Item \ref{enu:spm AutF-invariant},
we may assume without loss of generality that the edges of $\Omega$
correspond to the elements of $S$. Assume again that there exists
a valid diagram as in \eqref{eq:spm with core graphs}. Let $e\in E(\Gamma)$.
Efficiency guarantees that for disjoint edges $\varepsilon_{1},\varepsilon_{2}\in b^{-1}(e)\subseteq E(P)$,
we have $\rho(\varepsilon_{1})\ne\rho(\varepsilon_{2})$. The condition
about the basis elements yields, therefore, that there are at most
$m-1$ different edges in $b^{-1}(e)$ with a given orientation. Hence
$|n_{b}(e)|\le m-1$, and as $n_{b}(e)\in m\mathbb{Z}$, we conclude
that $n_{b}(e)=0$. This is true for all $e\in E(\Gamma)$ and we
finish as in the case of $m=0$.

\eqref{enu:spm inside =00005BF,F=00005D} If $w\in[\F,\F]$, then
the diagram 
\begin{equation}
\xymatrix{\Gamma_{w}\ar@{->>}[d]^{\id~~}\ar[r]^{b} & \Omega\ar[d]^{\id}\\
\Gamma_{w}\ar[r]^{\eta_{w}} & \Omega
}
\label{eq:trivial diagram}
\end{equation}
satisfies all the conditions in \eqref{eq:spm with core graphs} for
all $1\ne m\in\mathbb{Z}_{\ge0}$ (in particular, $n_{b}(e)=0$ for
every $e\in E(\Omega)$), so $\spm(w)\le\frac{-\chi(\Omega)}{1}=\rk\F-1$.
For $m=1$ this follows from Corollary \ref{cor:ineq of spi, spi tilde and pi}
and the fact that $\pi(w)\le\rk\F-1$ when $w$ is imprimitive. 
\end{proof}
\begin{example}
\label{exa:spm of aBcbbaCac}Consider the word $w=ab^{-1}cb^{2}ac^{-1}ac$
from Table \ref{tab:values on selected words} and Example \ref{exa:ssql of aBcbbaCac}.
By Proposition \ref{prop:properties of spm}\eqref{enu:spm outside =00005BF,F=00005D},
$\spm(w)=\infty$ for $m\ge4$ and for $m=0$. Recall that the ordinary
Whitehead graph of $w$ (namely, $\wh_{\eta_{w}}\left(o\right)$ where
$\eta_{w}\colon\Gamma_{w}\to\Omega$) is the full bipartite graph
$K_{3,3}$ with no double edges, and that in every efficient morphism
$b\colon\Gamma_{w^{\sigma}}\to\Gamma$ and every vertex $v$ of $\Gamma$,
the Whitehead graph $\wh_{b}\left(v\right)$ is a subgraph of this
$K_{3,3}$. These subgraphs must satisfy, locally, the condition that
$n_{b}\left(e\right)\in m\mathbb{Z}$ -- see the proof of Theorem
\ref{thm:rationality of spm}. When $m=3$, a simple check yields
that the few valid pieces (à la proof of Theorem \ref{thm:rationality of spm})
cannot be glued together so $\sp^{(3)}\left(w\right)=\infty$. When
$m=2$, the valid pieces are exactly those we have for $\ssql$, so
$\sp^{\left(2\right)}\left(w\right)=\ssql\left(w\right)=\frac{9}{4}$.
Finally, for $m=1$, as there are only three instances of every letter
in $w$, in every valid diagram as in \ref{eq:def of spi with WH}
$f$ is an immersion. So $\wh_{b}(v)$ is a subgraph of $\wh_{\eta_{w}}\left(o\right)$.
We can proceed with analysis similar to the one from Example \ref{exa:ssql of aBcbbaCac},
noticing that 
\[
\frac{-\chi(\Gamma)}{\deg\rho}=\frac{|w|}{2}\cdot\frac{\sum_{v\in V(\Gamma)}(\deg v-2)}{\sum_{v\in V(\Gamma)}|E(\wh_{b}(v))|}.
\]
The full graph $\wh_{\eta_{w}}\left(o\right)$ gives the smallest
possible ratio $\frac{\deg\left(v\right)-2}{\left|E\left(\wh_{\eta_{w}}\left(o\right)\right)\right|}=\frac{4}{9}$,
so \eqref{eq:trivial diagram} is an extremal diagram, yielding $\sp(w)=\frac{9}{2}\cdot\frac{4}{9}=2$.
\end{example}

\subsection{The stable characters of $C_{m}\wr S_{\bullet}$\label{subsec:stable characters of CmSn}}

Recall that Theorem \ref{thm:result on wreath products and modulo-m spi},
which we prove in $\S\S$\ref{subsec:spm-main-proof}, discusses word
measures on $C_{m}\wr S_{N}$ for $1\ne m\in\mathbb{Z}_{\ge0}$ and
$C_{m}=\left\{ z\in\mathbb{S}^{1}\,\middle|\,z^{m}=1\right\} $. Fix
such $m$ and throughout this section and denote $G_{N}\defi C_{m}\wr S_{N}$.
The elements of $G_{N}$ are $\left(v,\sigma\right)\in C_{m}^{\,N}\times S_{N}$,
$v=\left(v_{1},\ldots,v_{N}\right)$ with $\left(v,\sigma\right)\cdot\left(u,\tau\right)=\left(\left(v_{1}u_{\sigma^{-1}\left(1\right)},\ldots,v_{N}u_{\sigma^{-1}\left(N\right)}\right),\sigma\tau\right)$.
Equivalently, the elements of $G_{N}$ are $N\times N$ monomial matrices
(also known as generalized permutation matrices) with entries in $C_{m}\sqcup\left\{ 0\right\} $.
The element $\left(v,\sigma\right)$ corresponds to the monomial matrix
$A$ with $A_{i,\sigma^{-1}\left(i\right)}=v_{i}$, and the multiplication
corresponds to the ordinary matrix multiplication. The conjugacy class
of $\left(v,\sigma\right)\in G_{N}$ is determined by the numbers
$(a_{\ell,x})_{\ell\in\left[N\right],x\in C_{m}}$, where $a_{\ell,x}$
is the number of $\ell$-cycles in $\sigma$ with the product of the
corresponding $\ell$ $v_{i}$'s giving $x$ \cite[\S\S\S2.3.3]{ceccherini2014representation}.

The irreducible characters of $G_{N}$ are parameterized by functions
$\arrm\colon\irr(C_{m})\to{\cal P}$ from $\irr\left(C_{m}\right)$
to the set of all partitions, with $|\arrm|\defi\sum_{\phi}|\arrm(\phi)|=N$
(see \cite[\S2.6]{ceccherini2014representation}\footnote{Even though the reference \cite[\S2]{ceccherini2014representation}
concerns only wreath products $G\wr S_{N}$ where $G$ is a \emph{finite}
group, which is not the case for $C_{0}=S^{1}$, most of its content
applies as-is to arbitrary compact groups. In \cite{hora2008limits}
the theory is explicitly developed for general compact groups.}). We let $\supp(\arrm)\defi\{\phi\in\irr(G)\,|\,\arrm(\phi)\ne\emptyset\}$
denote the support of $\arrm$. In particular, as $|\arrm|<\infty$,
its support is finite (even in the case $m=0$). We denote the irreducible
character corresponding to $\arrm$ by $\chi^{\arrm}$. In the following
proposition, an \emph{ordered set partition} is a set partition (of
some set) with a prescribed order on the blocks of the partition.
\begin{prop}
\label{prop:irred chars of C_m wr S_N}
\begin{enumerate}
\item \label{enu:formula when support of size 1}Let $\mu\vdash N$ and
$\phi\in\irr\left(C_{m}\right)$. Denote by $\chi^{\phi,\mu}$ the
irreducible character of $G_{N}$ corresponding to the function assigning
$\mu$ to $\phi$ and $\emptyset$ to any other $\phi'\in\irr(C_{m})$.
Then $\chi^{\phi,\mu}$ is given by
\begin{equation}
\chi^{\phi,\mu}\left(v,\sigma\right)=\chi^{\mu}\left(\sigma\right)\prod_{j=1}^{N}\phi\left(v_{j}\right),\label{eq:formula when arrm supported on a single irr-char}
\end{equation}
for every $(v,\sigma)\in G_{N}$, where $\chi^{\mu}\in\irr\left(S_{N}\right)$.
\item \label{enu:formula for character as induction}Let $\arrm\colon\irr(C_{m})\to{\cal P}$
with $|\arrm|=N$ and $\supp(\arrm)=\{\phi_{1},\ldots,\phi_{k}\}$.
Denote $d_{i}=|\arrm\left(\phi_{i}\right)|$, so that $d_{i}\ge1$
and $\sum d_{i}=N$. The irreducible character corresponding to $\arrm$
is given by 
\[
\chi^{\arrm}=\mathrm{Ind_{G_{d_{1}}\times\ldots\times G_{d_{k}}}^{G_{N}}\left(\chi^{\phi_{1},\arrm\left(\phi_{1}\right)}\boxtimes\ldots\boxtimes\chi^{\phi_{k},\arrm\left(\phi_{k}\right)}\right).}
\]
\item \label{enu:formula for value of induced character}Keeping the notation
from Item \ref{enu:formula for character as induction}, and fixing
$\left(v,\sigma\right)\in G_{N}$, let ${\cal B}_{\sigma}$ denote
the set of ordered set-partitions $\overline{B}=(B_{1},\ldots,B_{k})$
of $[N]$ into $k$ subsets $B_{1}\sqcup\ldots\sqcup B_{k}=[N]$ with
$\left|B_{i}\right|=d_{i}$ which are invariant under $\sigma\in S_{N}$,
namely, such that $\sigma(B_{i})=B_{i}$ for all $i$. Then
\begin{equation}
\chi^{\arrm}\left(v,\sigma\right)=\sum_{\overline{B}\in{\cal B}_{\sigma}}\prod_{i=1}^{k}\chi^{\arrm(\phi_{i})}\left(\sigma|_{B_{i}}\right)\prod_{x\in B_{i}}\phi_{i}(v_{x}).\label{eq:formula for induced char}
\end{equation}
\end{enumerate}
\end{prop}

\begin{proof}
All items follow quite easily from \cite{ceccherini2014representation}
and \cite[\S1.2]{hora2008limits}. We elaborate further in the proof
of a more general version of this theorem in the companion paper \cite[Prop.~2.13]{PSh25}:
that version applies to the wreath product $G\wr S_{N}$ with $G$
an arbitrary compact group. The formulas in Proposition \ref{prop:irred chars of C_m wr S_N}
are slightly simpler than in the general case thanks to $\phi$ being
a linear character. 
\end{proof}
\emph{Stable} irreducible characters of $C_{m}\wr S_{\bullet}$ are
parameterized by partition-valued functions $\arrm\colon\irr(C_{m})\to{\cal P}$
with $|\arrm|<\infty$ \cite{sam2019representations}. For every $N\ge|\arrm|+\arrm(\mathrm{triv)}_{1}$,
we define $\arrm[N]$ as the function which agrees with $\arrm$ on
all non-trivial representations and is $\arrm(\mathrm{triv)}$$\left[N'\right]$
where $N'=N-\sum_{\phi\ne\triv}|\arrm(\phi)|$ on the trivial representation
(the notation $\arrm(\mathrm{triv)}$$\left[N\right]$ is defined
in $\S\S$\ref{subsec:Evidence-towards-Conjecture}). The stable irreducible
character corresponding to $\arrm$ is 
\[
\chi^{\arrm[\bullet]}\defi\left\{ \chi^{\arrm\left[N\right]}\right\} {}_{N\ge\left|\arrm\right|+\arrm\left(\mathrm{triv}\right)_{1}}.
\]
 
\begin{prop}
\label{prop:dim of chi^arrm}Let $\arrm\colon\irr(C_{m})\to{\cal P}$
with $\left|\arrm\right|<\infty$. Then
\[
\dim\left(\chi^{\arrm\left[N\right]}\right)=\Theta\left(N^{\left|\arrm\right|}\right).
\]
\end{prop}

\begin{proof}
This follows quite easily from Proposition \ref{prop:irred chars of C_m wr S_N}.
We elaborate in \cite[Prop.~2.14]{PSh25}, where we state a more general
version of this fact.
\end{proof}

\subsection{Proof of Theorem \ref{thm:result on wreath products and modulo-m spi}
\label{subsec:spm-main-proof}}

\subsubsection{Preliminaries}

We keep the notation $G_{N}=C_{m}\wr S_{N}$ for some fixed $1\ne m\in\mathbb{Z}_{\ge0}$.
In the upcoming proof of Theorem \ref{thm:result on wreath products and modulo-m spi}
we use some additional concepts and lemmas. We start with random $G$-labelings
of a graph.
\begin{defn}
\label{def:G-lableing of a graph}Let $G$ be a compact group and
$\Gamma$ a finite Serre graph\footnote{Serre graphs were discussed in $\S\S$\ref{subsec:Whitehead-graphs}.}.
A\textbf{ $G$-labeling} of $\Gamma$ is a map $\gamma\colon E\left(\Gamma\right)\to G$
satisfying $\gamma\left(\overline{e}\right)=\gamma\left(e\right)^{-1}$
for all $e\in E\left[\Gamma\right]$. A \textbf{Haar-random $G$-labeling}
of $\Gamma$ is a random $G$-labeling such that $\gamma\left(e\right)$
is Haar-random for every $e\in E\left(\Gamma\right)$, and the $G$-labels
on all edges are independent except for the constraint that $\gamma\left(\overline{e}\right)=\gamma\left(e\right)^{-1}$. 
\end{defn}

If $\Gamma$ is connected, we may identify $\pi_{1}(\Gamma)\cong\F_{r}$
for some $r$, where $\F_{r}$ is a rank-$r$ free group. The set
$\Hom(\pi_{1}(\Gamma),G)$ has a natural probability measure coming
from the Haar measure on $G^{r}$ (the direct product of $r$ copies
of $G$) through the isomorphism of sets $\Hom(\pi_{1}(\Gamma),G)\cong\Hom(\F_{r},G)\cong G^{r}$.
A Haar-random $G$-labeling can be used to encode a Haar-random homomorphism
$\pi_{1}(\Gamma)\to G$:
\begin{lem}
\label{lem:labeling to homomorphism to G}Let $\Gamma$ be a finite,
connected Serre graph with a fixed vertex $v\in V(\Gamma)$ and let
$G$ be a compact group. Every $G$-labeling of $\Gamma$ gives rise
to a well-defined homomorphism $\pi_{1}(\Gamma,v)\to G$ where for
any closed path $p$ at $v$ we map $[p]$ to the product in $G$
of the labels of the edges $p$ traverses. A Haar-random $G$-labeling
of $\Gamma$ gives rise to a Haar-random homomorphism $\pi_{1}(\Gamma,v)\to G$.
\end{lem}

\begin{proof}
This is standard. See, for example, \cite[\S\S2.3]{hall2018ramanujan}.
\end{proof}
In the special case where $G$ is a symmetric group, labelings have
a natural geometric interpretation in the form of covering spaces.
This concept appeared countless times in the literature. We give a
technical description, but see \cite[P.~9244]{hanany2020word} for
a more intuitive one.
\begin{defn}
\label{def:random-cover}Let $\Gamma$ be a finite Serre graph, $N\in\mathbb{Z}_{\ge1}$
and $\gamma$ an $S_{N}$-labeling of $\Gamma$. Let $V$ and $E$
be the vertex and edge sets of $\Gamma$, respectively. The $N$-cover
of $\Gamma$ corresponding to $\gamma$ has vertex set $V\times[N]$
and edge set $E\times[N]$ with $\overline{(e,j)}=(\overline{e},\gamma_{e}(j))$
and $\iota((e,j))=(\iota(e),j)$. A \textbf{random $N$-cover of $\Gamma$}
is the $N$-cover corresponding to a random $S_{N}$-labeling of $\Gamma$. 
\end{defn}

It is immediate to check that this construction yields, indeed, a
topological $N$-cover of $\Gamma$. In particular, a random $N$-cover
of $\Gamma$ corresponds to a random homomorphism $\pi_{1}(\Gamma)\to S_{N}$.
One may use random covers also in the case of $G\wr S_{N}$: 
\begin{lem}
\label{lem:A-Haar-random-homomorphism by labeling of a random cover}Let
$\Gamma$ be a finite, connected Serre graph and $G$ a compact group.
A Haar-random homomorphism $\pi_{1}(\Gamma)\to G_{N}$ can be realized
by taking a random $N$-cover $\hat{\Gamma}$ of $\Gamma$, together
with a Haar-random $G$-labeling of $\hat{\Gamma}$.
\end{lem}

\begin{proof}
An $N$-cover $\hat{\Gamma}$ of $\Gamma$ corresponds to an $S_{N}$-labeling
of $\Gamma$, and a $G$-labeling of $\hat{\Gamma}$ consists of assigning
a (legal) label from $G$ to every edge in $E(\hat{\Gamma})=E(\Gamma)\times[N]$.
Overall, for every $e\in E[\Gamma]$ we chose some $(v,\sigma)\in G^{N}\times S_{N}$
(which, as a set, is precisely $G\wr S_{N}$). It is easy to check
that this gives rise to a valid $G\wr S_{N}$-labeling of $\Gamma$.
\end{proof}
We also use the following result about random $d$-covers of graphs:
\begin{lem}
\label{lem:introducing L}Let $\eta\colon\Gamma\to\Delta$ be an immersion
of core graphs and $\widehat{\Delta_{N}}$ a random $N$-cover of
$\Delta$. Denote by $L_{\eta}(N)$ the expected number of \emph{injective}
lifts of $\eta$ into $\widehat{\Delta_{N}}$. 
\[
\xymatrix{ & \widehat{\Delta_{N}}\ar[d]^{\rho}\\
\Gamma\ar[r]_{\eta}\ar@{^{(}-->}[ur]^{\#} & \Delta
}
\]
Then 
\[
L_{\eta}(N)=N^{\chi(\Gamma)}\left(1+O\left(N^{-1}\right)\right).
\]
\end{lem}

\begin{proof}
This is immediate from \cite[Prop.~6.8]{hanany2020word}.
\end{proof}

\subsubsection{The proof}

Recall that for $1\ne m\in\mathbb{Z}_{\ge0}$ we let $\phi_{m}\in\irr(C_{m})$
denote the defining character of $C_{m}$: for $m=0$ this is the
embedding $C_{0}\cong S^{1}\hookrightarrow\mathbb{C}^{*}$, and for
$m\ge2$ this is the character mapping the cyclic generator of $C_{m}$
to $e^{2\pi i/m}\in\mathbb{C}$. We denote by ${\cal I}_{m}$ the
set of stable irreducible characters of $C_{m}\wr S_{\bullet}$ corresponding
to $\arrm\colon\irr(C_{m})\to{\cal P}$ with $\supp(\arrm)=\{\phi_{m}\}$.
Theorem \ref{thm:result on wreath products and modulo-m spi} states
that for any $w\in\F$, 
\[
\spm(w)=\inf_{\chi\in{\cal I}_{m}}\beta(w,\chi)
\]
and that the infimum is attained. 
\begin{rem}
We may always define $\beta(w,\chi)$ as in \eqref{eq:formal-def-of-beta}.
In fact, for $m\ne0$, \cite[Thm~1.6]{shomroni2023wreathII} yields
that for every stable character $\chi=\{\chi_{N}\}_{N\ge N_{0}}$
of $C_{m}\wr S_{\bullet}$, $\mathbb{E}_{w}[\chi_{N}]$ coincides
with some $f_{w,\chi}\in\mathbb{C}(N)$ for every large enough $N$,
and in light of Proposition \ref{prop:dim of chi^arrm} we get that
$\beta(w,\chi)\in\mathbb{Q}\sqcup\{\infty\}$ as in \eqref{eq:beta for S_bullet}.
As for $m=0$, the following proof of Theorem \ref{thm:result on wreath products and modulo-m spi}
yields that the same conclusion is true at least for $\chi\in{\cal I}_{0}$. 
\end{rem}

\begin{proof}[Proof of Theorem \ref{thm:result on wreath products and modulo-m spi}]
We first show that every $\arrm$ with $\supp(\arrm)=\{\phi_{m}\}$
satisfies $\spm(w)\le\beta(w,\chi^{\arrm[\bullet]})$. Denote $\mu=\text{\ensuremath{\arrm\left(\phi_{m}\right)}}$
and $d=\left|\arrm\right|=\left|\mu\right|\ge1$. Imitating the notation
from Proposition \ref{prop:irred chars of C_m wr S_N} (with $\phi_{1}=\phi_{m}$
and $\phi_{2}=\triv$), the set ${\cal B}_{\sigma}$ is simply $\{B\in\binom{[N]}{d}|\sigma(B)=B\}$.
As $\chi^{\mathrm{triv},\left(N-d\right)}$ is the trivial character
on $G_{N-d}$, we get for $\left(v,\sigma\right)\in G_{N}$ 
\begin{equation}
\chi^{\arrm\left[N\right]}\left(v,\sigma\right)=\sum_{B\in\binom{[N]}{d}\colon\sigma(B)=B}\chi^{\mu}\left(\sigma|_{B}\right)\prod_{x\in B}\phi_{m}(v_{x}).\label{eq:formula for chi^arrm=00005BN=00005D  when supported on phi_m}
\end{equation}
Now, assume that $\left(v,\sigma\right)$ is a $w$-random element
of $G_{N}$. Word measures are conjugation-invariant, so for any given
$B\in\binom{[N]}{d}$, the probability that $\sigma(B)=B$ is equal
to the probability that $\sigma([d])=[d]$, and when using \eqref{eq:formula for chi^arrm=00005BN=00005D  when supported on phi_m}
to estimate $\mathbb{E}_{w}[\chi^{\arrm[N]}]$, the total contribution
of $B$ is the same as that of $[d]$. We conclude that

\begin{eqnarray*}
\mathbb{E}_{w}\left[\chi^{\arrm[N]}\right] & = & \binom{N}{d}\sum_{\alpha\in S_{d}}\text{Prob}_{w}\left(\sigma|_{[d]}=\alpha\right)\chi^{\mu}(\alpha)\mathbb{E}_{w}\left[\prod_{x\in[d]}\phi_{m}(v_{x})\,\big|\,\sigma|_{[d]}=\alpha\right].
\end{eqnarray*}
A $w$-random element of a group $G$, defined in $\S$\ref{sec:Introduction},
can be equivalently defined as the image of $w$ through a Haar-random
homomorphism $\F\to G$. By Lemma \ref{lem:A-Haar-random-homomorphism by labeling of a random cover},
we may realize a Haar-random homomorphism $\F\to G_{N}$ by a Haar-random
$C_{m}$-labeling $\gamma$ of a uniformly random $N$-cover $\hat{\Omega}$
of the bouquet $\Omega$. We denote by $o$ the unique vertex of $\Omega$
(as above), and by $o_{j}=(o,j)$, $j=1,\ldots,N$, the $N$ vertices
of $\hat{\Omega}$. 

For any $\alpha\in S_{d}$, label by $x_{1},\ldots,x_{d}$ the $d$
vertices of $\Gamma_{w^{\alpha}}$ representing the starting points
of $w$, so that there is a path reading $w$ from $x_{i}$ to $x_{j}$
if and only if $\alpha(i)=j$.\footnote{There is an immaterial issue here: we need to assume without loss
of generality that $w$ is cyclically reduced: otherwise, the starting
point of $w$ is not part of the Stallings core graphs $\Gamma_{w}$
or $\Gamma_{w^{\alpha}}$.} We have that $\sigma|_{\left[d\right]}=\alpha$ if and only if the
injection $x_{i}\mapsto o_{i}$ for $i=1,\ldots,d$, extends to a
morphism $\eta\colon\Gamma_{w^{\alpha}}\to\hat{\Omega}$ such that
\[
\xymatrix{\Gamma_{w^{\alpha}}\ar[d]\ar[r]^{\eta} & \hat{\Omega\ar[d]^{\rho}}\\
\Gamma_{w}\ar[r]^{\eta_{w}} & \Omega
}
\]
commutes. Now assume this is the case. The expectation of $\prod_{x\in[d]}\phi_{m}(v_{x})$
is then 
\[
\mathbb{E}_{\gamma}\left[\prod_{\left\{ e,\overline{e}\right\} \subseteq E\left(\hat{\Omega}\right)}\gamma(e)^{n_{\eta}(e)}\right]=\prod_{\left\{ e,\overline{e}\right\} \subseteq E(\hat{\Omega})}\mathbb{E}_{\gamma}\left[\gamma(e)^{n_{\eta}(e)}\right],
\]
and the latter product is either $1$ if $n_{\eta}(e)\in m\mathbb{Z}$
for all edges $e\in E(\hat{\Omega})$ or vanishes otherwise. 

By symmetry, there is no difference in these numbers between the injection
$x_{i}\mapsto o_{i}$ and any other injection $\left\{ x_{1},\ldots,x_{d}\right\} \hookrightarrow\left\{ o_{1},\ldots,o_{N}\right\} $.
The number of such injections is $\frac{N!}{\left(N-d\right)!}$.
A morphism $\eta\colon\Gamma_{w^{\alpha}}\to\hat{\Omega}$ is injective
on $\left\{ x_{1},\ldots,x_{d}\right\} $ if and only if it is efficient.
Every efficient morphism $\eta$ factors uniquely as $\Gamma_{w^{\alpha}}\stackrel{\eta_{1}}{\twoheadrightarrow}\Sigma\stackrel{\overline{\eta_{2}}}{\hookrightarrow}\hat{\Omega}$
with $\eta_{1}$ efficient and surjective, and $\overline{\eta_{2}}$
injective. Given such $\eta_{1}$, let $\eta_{2}=\rho\circ\overline{\eta_{2}}\colon\Sigma\to\Omega$.
Note that $n_{\eta}(e)\in m\mathbb{Z}$ for all edges $e\in E(\hat{\Omega})$
if and only if $n_{\eta_{1}}(e)\in m\mathbb{Z}$ for all edges $e\in E(\Sigma)$.
Using Lemma \ref{lem:introducing L}, we obtain 

\begin{eqnarray}
\mathbb{E}_{w}\left[\chi^{\arrm[N]}\right] & = & \binom{N}{d}\cdot\sum_{\substack{\alpha\in S_{d}}
}\chi^{\mu}\left(\alpha\right)\frac{\left(N-d\right)!}{N!}\cdot\sum_{\substack{\Gamma_{w^{\alpha}}\stackrel{\eta_{1}}{\twoheadrightarrow}\Sigma\stackrel{\eta_{2}}{\to}\Omega\colon\\
\eta_{1}~\mathrm{efficient~and~surjective}\\
n_{\eta_{1}}(e)\in m\mathbb{Z}\,\,\forall e\in E(\Sigma)
}
}L_{\eta_{2}}\left(N\right)\nonumber \\
 & = & \frac{1}{d!}\sum_{\alpha\in S_{d}}\chi^{\mu}\left(\alpha\right)\sum_{\substack{\Gamma_{w^{\alpha}}\stackrel{\eta_{1}}{\twoheadrightarrow}\Sigma\stackrel{\eta_{2}}{\to}\Omega\colon\\
\eta_{1}~\mathrm{efficient~and~surjective}\\
n_{\eta_{1}}(e)\in m\mathbb{Z}\,\,\forall e\in E(\Sigma)
}
}N^{\chi\left(\Sigma\right)}\left(1+O\left(N^{-1}\right)\right).\label{eq:step in proof for spf}
\end{eqnarray}
Now, in every summand in \eqref{eq:step in proof for spf}, the diagram 

\[
\xymatrix{\Gamma_{w^{\alpha}}\ar[d]^{\rho}\ar[r]^{\eta_{1}} & \Sigma\ar[d]^{\eta_{2}}\\
\Gamma_{w}\ar[r]^{\eta_{w}} & \Omega
}
\]
satisfies the conditions in Definition \ref{def:spm}, so $\spm\left(w\right)\le\frac{-\chi\left(\Sigma\right)}{\deg\left(\rho\right)}=\frac{-\chi\left(\Sigma\right)}{d}$,
namely, $\chi\left(\Sigma\right)\le-d\cdot\spm\left(w\right)$ and
by Proposition \ref{prop:dim of chi^arrm},
\[
\mathbb{E}_{w}\left[\chi^{\arrm[N]}\right]=O\left(\left(N^{d}\right)^{-\spm(w)}\right)=O\left(\left(\dim\chi^{\arrm[N]}\right)^{-\spm(w)}\right).
\]
It remains to show that the latter bound is tight for certain stable
characters $\chi\in{\cal I}_{m}$. By Theorem \ref{thm:rationality of spm},
there exist extremal diagrams as in \eqref{eq:spm with core graphs}
attaining the value of $\spm\left(w\right)$. Let $d$ be the degree
of $\rho$ in some extremal diagram for $w$ and consider $\arrm$
supported on $\phi_{m}$ with $\arrm\left(\phi_{m}\right)=\left(d\right)$
(the partition with one block of size $d$). As $\chi^{\left(d\right)}$
is the trivial character of $S_{d}$, in this case \eqref{eq:step in proof for spf}
becomes
\[
\mathbb{E}_{w}\left[\chi^{\arrm[N]}\right]=\frac{1}{d!}\sum_{\alpha\in S_{d}}\sum_{\substack{\Gamma_{w^{\alpha}}\stackrel{\eta_{1}}{\twoheadrightarrow}\Sigma\stackrel{\eta_{2}}{\to}\Omega\colon\\
\eta_{1}~\mathrm{efficient~and~surjective}\\
n_{\eta_{1}}(e)\in m\mathbb{Z}\,\,\forall e\in E(\Sigma)
}
}N^{\chi\left(\Sigma\right)}\left(1+O\left(N^{-1}\right)\right).
\]
The existence of an extremal diagram of degree $d$ guarantees that
the coefficient of $N^{-d\cdot\spm\left(w\right)}$ is strictly positive.
\end{proof}
We can now conclude Corollary \ref{cor:sp(m) profinite} that $\spm$
is profinite for every $1\ne m\in\mathbb{Z}_{\ge0}$.
\begin{proof}[Proof of Corollary \ref{cor:sp(m) profinite}]
 For $m\ne0$, the groups $C_{m}\wr S_{N}$ are finite, so the claim
is immediate from Theorem \ref{thm:result on wreath products and modulo-m spi}:
if $w_{1}$ and $w_{2}$ induce the same measure on $C_{m}\wr S_{N}$
for all $N$, then clearly $\spm\left(w_{1}\right)=\spm\left(w_{2}\right)$.
For $m=0$, the argument in the proof of Proposition \ref{prop:properties of spm}\eqref{enu:spm outside =00005BF,F=00005D}
shows that if $m>|w|$ then $\spm(w)=\sp^{(0)}(w)$. Hence, if $m>|w_{1}|,|w_{2}|$
and both words induce the same measure on $C_{m}\wr S_{N}$ for every
$N$, then $\sp^{\left(0\right)}\left(w_{1}\right)=\spm\left(w_{1}\right)=\spm\left(w_{2}\right)=\sp^{\left(0\right)}\left(w_{2}\right)$.
\end{proof}
We end this subsection with a result similar in nature to Theorem
\ref{thm:result on wreath products and modulo-m spi}. It applies
to more general stable irreducible characters, but gives a weaker
conclusion.
\begin{thm}
\label{thm:sp lower bound for stables supported outside triv}Let
$\arrm\colon\irr(C_{m})\to{\cal P}$ satisfy $|\arrm|<\infty$ and
$\arrm(\triv)=\emptyset$. Then for all $w\in\F$,
\[
\sp(w)\le\beta\left(w,\chi^{\arrm[\bullet]}\right).
\]
\end{thm}

We defer the proof of Theorem \ref{thm:sp lower bound for stables supported outside triv}
to \cite{PSh25}, where we prove (in Theorem 1.5) a more general version
applying to $G\wr S_{\bullet}$ for arbitrary compact group $G$.
We do rely on Theorem \ref{thm:sp lower bound for stables supported outside triv}
in the proof of Theorem \ref{thm:spm avoids (0,1)} below.

\subsection{$\protect\spm$ admits a gap in $\left(0,1\right)$\label{subsec:gap of spm} }

Throughout this paper, we introduced stable invariants of words, recalled
or established many of their properties, and recalled, established
or conjectured their relations to word measures on groups, often relying
on the ``intrinsic'' properties of the invariants to establish these
connections (for example, we used the rationality Theorem \ref{thm:rationality of spm}
in the proof of Theorem \ref{thm:result on wreath products and modulo-m spi}).
In this final part of $\S$\ref{sec:spm + G wr S}, we exploit the
inverse direction and use the results on word measures to establish
an ``intrinsic'' property of $\spm$: that this stable invariant,
too, admits a gap in $\left(0,1\right)$. To the best of our knowledge,
this is the first result that goes in this direction.

Recall that if $w$ is a non-power, then $\spm(w)\ge\sp(w)\ge1$ by
Proposition \ref{prop:properties of spm} and Theorem \ref{thm:HW-LW},
and that $\spm(1_{\F})=-1$. It remains to show that $\spm$ avoids
$\left(0,1\right)$ also when evaluated on proper powers. 
\begin{thm}
\label{thm:spm avoids (0,1)}Let $1\ne w=u^{t}\in\F$ be a proper
power, with $u$ a non-power and $t\ge2$. Then $\spm(w)\notin(0,1)$.
\end{thm}

In the proof, we rely on a representation-theoretic result. Let $\chi$
be a character of the group $G$. For every $t\in\mathbb{Z}_{\ge1}$
denote by $\chi^{(t)}\colon G\to\mathbb{C}$ the class function defined
by $\chi^{(t)}(g)=\chi(g^{t})$. If $\chi_{\bullet}$ is a stable
character of $C_{m}\wr S_{\bullet}$, then $\chi_{\bullet}^{(t)}$
has a stable decomposition. Namely, 
\begin{equation}
\chi_{\bullet}^{(t)}=\sum_{\xi_{\bullet}}\alpha_{\xi}\xi_{\bullet}\label{eq:decomposition of a stable under t power}
\end{equation}
where the sum is over the stable irreducible characters of $C_{m}\wr S_{\bullet}$,
$\alpha_{\xi}\in\mathbb{C}$,\footnote{In fact, $\chi^{(t)}$ is always a virtual character of $G$ \cite[Problem 4.7]{isaacs1994character},
so $\alpha_{\xi}\in\mathbb{Z}$.} and all but finitely many of the $\alpha_{\xi}$'s vanish.\footnote{This follows from \cite[Cor.~2.30]{shomroni2023wreathII} for the
groups $C_{m}\wr S_{\bullet}$, but should be true for any stable
character $\chi_{\bullet}$ of an arbitrary sequence of compact groups
$G_{\bullet}$ with stable representation theory. In fact, Nir Gadish
explained to us the proof of this fact in private communication.} The equality \eqref{eq:decomposition of a stable under t power}
means that for every large enough $N$ there is an actual equality
when both sides are specialized to $C_{m}\wr S_{N}$.
\begin{prop}
\label{prop:on the decomposition of arrm}Let $\arrm\colon\irr(C_{m})\to{\cal P}$
\textup{satisfy $\supp(\arrm)=\{\phi_{m}\}$ and let $t$ be a positive
integer. If $m\ne0$, assume further that $t<m$. Consider the decomposition
\begin{equation}
(\chi^{\arrm[\bullet]})^{(t)}=\sum_{\arrn}\alpha_{\arrn}\chi^{\arrn[\bullet]}.\label{eq:decomp of power-t-of-stable in CmS_bullet}
\end{equation}
Then $\alpha_{\arrn}\ne0$ only if} $\arrn(\triv)=\emptyset$ and
$|\arrn|\ge|\arrm|$.
\end{prop}

We prove Proposition \ref{prop:on the decomposition of arrm} after
we show how it is used to prove Theorem \ref{thm:spm avoids (0,1)}.
\begin{proof}[Proof of Theorem \ref{thm:spm avoids (0,1)} assuming Proposition
\ref{prop:on the decomposition of arrm}]
 Let $1\ne w=u^{t}\in\F$ with $u$ a non-power and $t\ge2$. First,
if $m\ne0$ and $t\ge m$, Example \ref{exa:spm} shows that $\spm(w)=0$.
So we may assume that either $m=0$ or $t<m$ and show that in these
cases $\spm(w)\ge1$. By Theorem \ref{thm:result on wreath products and modulo-m spi},
there is some $\arrm$ supported on $\phi_{m}$ with $\spm(w)=\beta(w,\chi^{\arrm[\bullet]})$.
Denote $\mu=\arrm(\phi_{m})$ and $d=|\mu|=|\arrm|$. So, $\mathbb{E}_{w}[\chi^{\arrm[N]}]=\Theta(N^{-d\cdot\spm(w)})$. 

On the other hand, using the decomposition of $\chi^{\arrm[\bullet]}$
as in \eqref{eq:decomp of power-t-of-stable in CmS_bullet},
\begin{eqnarray*}
\mathbb{E}_{w}\left[\chi^{\arrm\left[N\right]}\right] & = & \mathbb{E}_{u}\left[\left(\chi^{\arrm\left[N\right]}\right)^{(t)}\right]=\sum_{\arrn}\alpha_{\arrn}\mathbb{E}_{u}\left[\chi^{\arrn[N]}\right],
\end{eqnarray*}
and by Proposition \ref{prop:on the decomposition of arrm}, in the
latter summation the non-vanishing summands satisfy $\arrn(\triv)=\emptyset$
and $\left|\arrn\right|\ge d$. For these, by Theorem \ref{thm:sp lower bound for stables supported outside triv},
\[
\mathbb{E}_{u}\left[\chi^{\arrn[N]}\right]=O\left(N^{-\left|\arrn\right|\cdot\sp(u)}\right)=O\left(N^{-d}\right),
\]
where in the last equality we used Theorem \ref{thm:HW-LW} that $\sp\ge1$
for non-powers. We conclude that $-d\cdot\spm(w)\le-d$, namely, $\spm(w)\ge1$.
\end{proof}
Finally, we prove Proposition \ref{prop:on the decomposition of arrm}.
\begin{proof}[Proof of Proposition \ref{prop:on the decomposition of arrm}]
 Denote $\mu=\arrm(\phi_{m})$ and $d=|\mu|=|\arrm|$. Using the
formula \eqref{eq:formula for chi^arrm=00005BN=00005D  when supported on phi_m},
we have 
\begin{eqnarray}
\left(\chi^{\arrm[N]}\right)^{(t)}\left(v,\sigma\right) & = & \chi^{\arrm[N]}\left(\left(v,\sigma\right)^{t}\right)\nonumber \\
 & = & \sum_{B\in\binom{[N]}{d}\colon\sigma^{t}(B)=B}\chi^{\mu}\left(\sigma^{t}|_{B}\right)\prod_{x\in B}\phi_{m}\left(v_{x}\right)\phi_{m}\left(v_{\sigma^{-1}(x)}\right)\cdots\phi_{m}\left(v_{\sigma^{-(t-1)}(x)}\right).\label{eq:formula for chi^=00005Carrm^(t)}
\end{eqnarray}
Let $\arrn\colon\irr(C_{m})\to{\cal P}$ satisfy $|\arrn|<\infty$
and imitate the notation from Proposition \ref{prop:irred chars of C_m wr S_N}
for $\arrn[N]$, only with $\psi_{1},\ldots,\psi_{k}$ (instead of
$\phi_{1},\ldots,\phi_{k}$), with $\psi_{k}=\triv$, and with $c_{1}=|\arrn(\psi_{1})|,\ldots,c_{k-1}=|\arrn(\psi_{k-1})|,c_{k}=N-c_{1}-\ldots-c_{k-1}$.
Then, by the same proposition,
\begin{equation}
\chi^{\arrn[N]}(v,\sigma)=\sum_{\overline{B}\in{\cal B}_{\sigma}}\left[\prod_{i=1}^{k-1}\chi^{\arrn\left(\psi_{i}\right)}\left(\sigma|_{B_{i}}\right)\prod_{x\in B_{i}}\psi_{i}\left(v_{x}\right)\right]\cdot\chi^{\arrn[N](\triv)}\left(\sigma|_{B_{k}}\right).\label{eq:formula for chi^arrn=00005BN=00005D}
\end{equation}
The coefficient $\alpha_{\arrn}$ is equal to the inner product $\left\langle \left(\chi^{\arrm[N]}\right)^{(t)},\chi^{\arrn[N]}\right\rangle _{G_{N}}$
estimated at any large enough $N$ (the value is constant for every
large enough $N$ -- see $\S\S$\ref{subsec:Stable-representation-theory}).
This inner product is equal to 
\begin{eqnarray}
 &  & \frac{1}{\left|G_{N}\right|}\sum_{\left(v,\sigma\right)\in G_{N}}\left(\chi^{\arrm[N]}\right)^{(t)}\left(v,\sigma\right)\cdot\overline{\chi^{\arrn[N]}\left(v,\sigma\right)}\nonumber \\
 & = & \sum_{B\in\binom{[N]}{d},\overline{B}\in\binom{[N]}{c_{1}\,\ldots\,c_{k}}}\frac{1}{\left|G_{N}\right|}\sum_{\substack{\sigma\in S_{N}\colon\\
\sigma(B)=B,\overline{B}\in{\cal B}_{\sigma}
}
}\sum_{v\in C_{m}^{\,N}}\chi^{\mu}\left(\sigma^{t}|_{B}\right)\prod_{x\in B}\phi_{m}\left(v_{x}\right)\cdots\phi_{m}\left(v_{\sigma^{-(t-1)}(x)}\right)\cdot\label{eq:sum over P and overline-P}\\
 &  & \,\,\,\,\,\,\,\,\,\,\,\,\,\,\,\,\,\,\,\,\,\,\,\,\,\,\,\,\,\,\,\,\,\,\,\,\,\,\,\,\,\,\,\,\,\,\,\,\,\,\,\,\,\,\,\,\,\,\,\,\,\,\,\,\,\,\,\,\,\,\,\,\,\,\,\,\,\,\,\,\,\,\cdot\overline{\left[\prod_{i=1}^{k-1}\chi^{\arrn\left(\psi_{i}\right)}\left(\sigma|_{B_{i}}\right)\prod_{x\in B_{i}}\psi_{i}\left(v_{x}\right)\right]\cdot\chi^{\arrn[N](\triv)}\left(\sigma|_{B_{k}}\right)}\nonumber 
\end{eqnarray}
We ought to show that this inner product vanishes if $|\arrn|<d$
or if $\arrn(\triv)\ne\emptyset$.

Assume first that $|\arrn|<d$. Fix some $B\in\binom{[N]}{d}$, $\overline{B}\in\binom{[N]}{c_{1}\,\ldots\,c_{k}}$
and $\sigma\in S_{N}$ with $\sigma(B)=B$ and $\overline{B}\in{\cal B}_{\sigma}$,
and consider the corresponding summand in \eqref{eq:sum over P and overline-P}.
As $c_{1}+\ldots+c_{k-1}\le|\arrn|<d$, there is some $x\in B\cap B_{k}$.
Then $\phi_{m}(v_{x})$ appears in the summand exactly $t$ times,
and $v_{x}$ does not appear with the $\psi_{i}$'s. As $m=0$ or
$t<m$, we have $\sum_{y\in C_{m}}\phi_{m}(y^{t})=0$, and the entire
summand vanishes.

Now assume that $\arrn(\triv)\ne\emptyset$, and fix some $B\in\binom{[N]}{d}$,
$\overline{B}\in\binom{[N]}{c_{1}\,\ldots\,c_{k}}$. If $B\cap B_{k}\ne\emptyset$,
we finish as in the previous case. Otherwise, $B_{k}\subseteq[N]\setminus B$
and the sum over $\sigma\in S_{N}$ with $\sigma(B)=B$ and $\overline{B}\in{\cal B}_{\sigma}$
can be split as a sum over $\sigma_{1}\in S_{B_{k}}$ (with no additional
conditions) and $\sigma_{2}\in S_{[N]\setminus B_{k}}$ (such that
$B$ and $B_{1},\ldots,B_{k-1}$ are invariant under $\sigma_{2}$).
The values $\{v_{x}\}_{x\in B_{k}}$ do not take part in the corresponding
summand in \eqref{eq:sum over P and overline-P}, and the only role
of $\sigma_{1}$ is in $\chi^{\arrn[N](\triv)}(\sigma_{1})$. But
as $\arrn(\triv)\ne\emptyset$, $\chi^{\arrn[N](\triv)}$ is a non-trivial
irreducible character of $S_{B_{k}}$ and 
\[
\sum_{\sigma_{1}\in S_{B_{k}}}\chi^{\arrn[N](\triv)}(\sigma_{1})=0.
\]
Again, the summand corresponding to $B$ and $\overline{B}$ vanishes.
\end{proof}

\section{Open questions\label{sec:Open-questions}}

In this final section (excluding the appendix), we list some of the
many intriguing open questions that this work gives rise to. 

\paragraph{Stable invariants and word measures}
\begin{itemize}
\item Many conjectures along this paper revolve around the role of the different
stable invariants in word measures of groups, and whether they can
be recovered using these measures: $\sp$ and $S_{\bullet}$ (Conjectures
\ref{conj:spi and stable characters of S_N} and \ref{conj:C^alg from cycle and non-efficient}),
$\ssql$ and $\U(\bullet)$, $\O(\bullet)$ and $\Sp(\bullet)$ (Conjectures
\ref{conj:ssql and its role for U,O and Sp}, \ref{conj:arbitrary stable chars of U}
and \ref{conj:stable chars of O}), and $\spq$ and $\gl_{\bullet}(q$)
(Conjecture \ref{conj:spq and stable characters}). There is also
Conjecture \ref{conj:universal spectral gap } about a universal gap
for the numbers $\beta(w,\chi)$ in $(0,1$). 
\item Related to these conjectures is the question of whether $\cl$, $\sqlh$,
$\ssql$, $\pi_{q}$, $\spq$ and $\spf$ (see \cite[Def.~4.7]{PSh25})
of a word $w$ are determined by the $w$-measures on compact groups,
and whether they are, moreover, profinite invariants (the latter question
is also open for $\scl$).
\end{itemize}

\paragraph*{Equalities and Inequalities between the different invariants}

Consider first the stable invariants $2\cdot\scl$, $\sp$, $\ssql$
and $\spm$ for $1\ne m\in\mathbb{Z}_{\ge0}$. They are all distinct.
Already the examples in Table \ref{tab:values on selected words}
show that every two of these stable invariants are distinct, except
for the pair $\sp^{\left(2\right)}$ and $\ssql$, and the pair $\sp^{\left(0\right)}$
and $2\cdot\scl$, which agree on all the examples. As for $\sp^{\left(0\right)}$
and $2\scl$, by \cite[Thm.~4.1]{calegari2013random}, a generic word
of (even) length $\ell$ in $\left[\F,\F\right]$ satisfies $2\scl\left(w\right)\approx\frac{\ell}{\log\left(\ell\right)}$.
However, every word in $\left[\F,\F\right]$ satisfies $\sp^{\left(0\right)}\left(w\right)\le\rk\F-1$
by Proposition \ref{prop:properties of spm}\eqref{enu:spm inside =00005BF,F=00005D}.
As for $\sp^{\left(2\right)}$ and $\ssql$, the lower bound in the
proof of \cite[Thm.~4.1]{calegari2013random} works almost intact
for random words and $\ssql$, showing that a generic word of length
$\ell$ in $\F$ satisfies that $\ssql\left(w\right)\ge c\frac{\ell}{\log\left(\ell\right)}$
for some constant $c$. However, if $\ell$ is even, then with probability
of about $2^{-\rk\F}$ (in particular, bounded away from zero), $w\in K_{2}\left(\F\right)$
(see \eqref{eq:pi modulo m}), and then $\sp^{\left(2\right)}\left(w\right)\le\rk\F-1$.
However, 
\begin{itemize}
\item Are the invariants $\spk$ for different fields $K$ different from
one another? We suspect they are all equal and, moreover, equal to
$\sp$ (compare with \cite[Conj.~1.9]{West}). 
\item Is it true that $\sp=\tilde{\sp}=\pi-1$ (Conjectures \ref{conj:spi and spi tilde}
and \ref{conj:spi=00003Dpi-1})? Is it true that $\spk=\pi_{K}-1$
for all $K$ (see Appendix \ref{sec:spq and GL_N(q)})? 
\item Does omitting the requirement that $f$ be $\pi_{1}$-injective in
Definition \ref{def:ssql} of $\ssql$ alter the resulting invariant?
Is it generally true that $\ssql\le\sqlh-1$? (See Remark \ref{rem:def of ssql}\eqref{enu:ssql without pi_1-inj issue}.)
\end{itemize}

\paragraph{Additional properties of the invariants}
\begin{itemize}
\item Do the invariants $\spk$ admit, too, a gap in $\left(0,1\right)$?
Are they rational? Computable?
\item In the companion paper \cite[\S4]{PSh25} we generalize $\spm$ in
the form of an invariant $\spf$ defined for an arbitrary non-trivial,
irreducible character of a compact group. Which of the properties
of $\spm$ does $\spf$ share? Is it rational? computable? Does it
admit a gap?
\item We know that $\sp$, $\ssql$, $\sp^{(2)}$ and $\spk$ all obtain
the value $\infty$ precisely for primitive words, and that $\scl$
and $\sp^{(0)}$ obtain the value $\infty$ precisely for words outside
$[\F,\F]$. When $m\ge3$, when does $\spm(w)=\infty$? Are there
reasons to obtain $\infty$ beyond the one recorded in Proposition
\ref{prop:properties of spm}\eqref{enu:spm outside =00005BF,F=00005D}?
\item We have a good understanding of the distribution of some invariants
evaluated on random words of a given length $\ell$, as $\ell\to\infty$.
This is true for $\scl$ \cite{calegari2013random} and for the ordinary
primitivity rank $\pi$ \cite{puder2015expansion,kapovich2022primitivity}.
What is the distribution of the other invariants -- $\sp,\ssql,\spm,\spk$,
$\spf$ and the non-stable counterparts?
\end{itemize}

\paragraph{Extensions of the theory}
\begin{itemize}
\item \textbf{The invariants in non-free groups}: Many of the invariants
defined in this paper -- both in the combinatorial/topological/algebraic
side and in the random-matrix side -- make sense also for elements
of finitely generated but non-free groups, and it seems that the theory
may be extended further to such groups, and especially to surface
groups or to free products of finite groups. For some seeds of this
theory consult \cite{magee2020asymptotic,magee2021surface-unitary1}
and especially \cite[\S7]{puder2024local}.
\item \textbf{The invariants for multiple words: }For some potential applications
(for example, for finite presentations of groups with more than one
relator) it makes sense to define the invariants not only for a single
word at a time but, more generally, for a tuple of words. There exist
some definitions in the literature -- Calegari defines $\scl$ of
a tuple of elements \cite[Def.~2.6.1]{calegari2009scl}, and Wilton
defined $\sp$ and $\ssql$ for general $2$-complexes -- but it
is not clear whether these are the right definitions for dealing with
word measures on groups.
\item \textbf{The supremum invariants:} In many of the definitions of stable
invariants along this paper, one may consider not only the infimum
over the numbers $\frac{-\chi(\Gamma)}{\deg\rho}$, but also the supremum.
In fact, Wilton considers these suprema in the case of $\sp$ and
$\ssql$ and shows they are equal \cite[Thm.~A]{wilton2024rational}.
They are also equal, therefore, to the supremum of the intermediate
$\sp^{(2)}$. Do the supremum invariants of $\spm$ ($m\ne2)$ or
$\spq$ or in the restrictive definition \eqref{eq:scl in free groups with WH}
of $\scl$ make sense? Do the different supremum invariants have any
random-matrix significance?
\end{itemize}

\section*{Appendix}

\begin{appendices}

\section{Stable $K$-primitivity rank and stable characters of $\protect\gl_{\bullet}\left(q\right)$\label{sec:spq and GL_N(q)}\protect \\
By Danielle Ernst-West, Doron Puder, Matan Seidel and Yotam Shomroni}

Fix an arbitrary field $K$. In this appendix we define the stable
$K$-primitivity rank $\spk$ of a word, and describe some of its
basic properties. When $K$ is finite of (prime power) size $q$,
we also denote this invariant $\spq$. The $K$-primitivity rank is
a creature different from $\scl$, $\sp$, $\ssql$ and $\spm$: it
is more algebraic in nature and less topological\textbackslash combinatorial.
We also know less about it than we currently know about the other
invariants. For example, we do not know whether it is computable and
always rational, and whether it admits a gap in $\left(0,1\right)$.
However, it is very much part, and an important one, in the story
of stable invariants and their role in word measures on groups, so
we choose to include this short appendix in the current paper. We
hope to write a separate paper dedicated to the stable $K$-primitivity
rank once we know more about it. 

When $K=\mathbb{F}_{q}$ is the finite field of size $q$, the $K$-primitivity
rank was born out of the study of word measures on $\gl_{\bullet}\left(q\right)=\gl_{\bullet}(\mathbb{F}_{q})$
and is defined in the world of free groups algebras. In analogy with
the Nielsen-Schreier theorem for free groups, Cohn and Lewin \cite{cohn1964free,lewin1969free}
proved that for any field $K$, every submodule of a free module of
$K\left[\F\right]$ is free, with a well-defined rank. In particular,
every one-sided ideal of $K\left[\F\right]$ is a free $K\left[\F\right]$-module. 

In \cite{West}, three of us show how word measures on $\gl_{\bullet}\left(q\right)$
are related to right ideals inside $\mathbb{F}_{q}\left[\F\right]$.
We defined there the $q$-primitivity rank of $w\in\F$, which can
be generalized to the $K$-primitivity rank for an arbitrary field
$K$ as
\[
\pi_{K}\left(w\right)=\min\left\{ \rk I\,\middle|\,~\begin{gathered}I\lvertneqq K\left[\F\right]~\mathrm{a~proper~right~ideal,}\\
w-1\in I~\mathrm{~not~a~primitive~element}
\end{gathered}
\right\} ,
\]
(an element of an ideal $I\le K[\F]$ is \emph{primitive} if it belongs
to some free basis), and showed that when $K=\mathbb{F}_{q}$, this
invariant is related to the expected number of fixed vectors in the
action of a $w$-random element of $\gl_{N}\left(q\right)$ on $\mathbb{F}_{q}^{~N}$.\footnote{When $K$ is infinite, $\gl_{N}(K)$ is not compact, admits no natural
and canonical probability measure and therefore also no natural probability
word measure.}

Here we suggest the following stable version of $\pi_{K}\left(w\right)$
which is defined in analogy with $\sp$. We first define a $w$-module
for a word $1\ne w\in\F$, which should be thought of as the analog
of a cover of $\Gamma_{w}$.
\begin{defn}
\label{def:w-module}A \textbf{$w$-module} (over a field $K$) is
a right $K[\F]$-module of the form
\[
W=W_{0}\otimes_{K\left[\langle w\rangle\right]}K[\F],
\]
where $W_{0}$ is a non-zero, finite dimensional right $K[\langle w\rangle]$-module.
The degree of $W$ is defined as 
\[
\deg(W)\defi\dim(W_{0})\in\mathbb{Z}_{\ge1}.
\]
A $w$-submodule of $W$ is $W'=W_{0}'\otimes_{K[\langle w\rangle]}K[\F]$,
where $W_{0}'\le W_{0}$ is a $K[\langle w\rangle]$-submodule. A
quotient of a $w$-module is $W/W'\cong(W_{0}/W_{0}')\otimes_{K[\langle w\rangle]}K[\F]$
where $W'$ is a $w$-submodule as above. 
\end{defn}

In the following definition, an \emph{algebraic $K[\F]$-module} is
a module $M$ which does not split as a direct sum $M=M_{1}\oplus M_{2}$
with $M_{2}$ a non-zero \emph{free} module. Equivalently, $M$ does
not admit a non-zero quotient which is free\footnote{An algebraic module is also called \emph{bound }in \cite[\S5.1]{cohn2006free}.}. 

\begin{defn}
\label{def:spq}Let $K$ be a field and $1\ne w\in\F$. The \textbf{stable
$K$-primitivity rank} of $w$ is 
\[
\spk\left(w\right)\defi\inf\left\{ \frac{-\chi\left(Q\right)}{\deg\left(W\right)}\,\middle|\,\begin{gathered}W~\mathrm{is~a~}w\textnormal{-}\mathrm{module,}\\
\phi\colon W\twoheadrightarrow Q~\mathrm{a~surjective\,homomorphism\,of\,right}\,K[\F]\textnormal{-}\mathrm{modules,}\\
\ker(\phi)\mathrm{\,is~an~algebraic~}K[\F]\textnormal{-}\mathrm{submodule,}\\
\phi~\mathrm{does~not~factor~through~a}\mathrm{~proper~quotient~}W\stackrel{\not\cong}{\twoheadrightarrow}\overline{W}\,\,\mathrm{of}\,\,w\textnormal{-}\mathrm{modules},\\
\mathrm{and~no~quotient~}W\twoheadrightarrow\overline{W}\,\,\mathrm{of}\,\,w\textnormal{-}\mathrm{modules}\,\mathrm{factors\,\,through}\,\,\phi.
\end{gathered}
\right\} .
\]
When $K$ is the finite field $\mathbb{F}_{q}$, we also denote this
invariant by $\spq$. When $w=1$ we set $\spk(1)=-1$.
\end{defn}

If $Q$ is a quotient of a $w$-module with an algebraic kernel, $Q$
is guaranteed to be finitely presented, namely, it is the quotient
$Q\cong P_{1}/P_{2}$ of two free $K[\F]$-modules of finite rank.
The Euler characteristic of $Q$ is then defined to be 
\[
\chi(Q)\defi\rk(P_{1})-\rk(P_{2}),
\]
and this number does not depend on the particular presentation (e.g.,
by Schanuel's lemma). The condition that $\phi$ does not factor through
any non-identity quotient $W\twoheadrightarrow\overline{W}$, namely,
that we do not allow 
\[
\xymatrix{W\ar@{->>}[rr]^{\phi}\ar@{->>}[rd]_{\not\not\cong} &  & Q\\
 & \overline{W}\ar@{->>}[ru]
}
\]
is the analog of the efficiency condition in Definition \ref{def:spi}
of $\sp$. The condition that no quotient $W\twoheadrightarrow\overline{W}$
of $w$-modules factors through $\phi$, namely, that we do not allow
\[
\xymatrix{W\ar@{->>}[rr]\ar@{->>}[rd]_{\phi} &  & \overline{W}\\
 & Q\ar@{->>}[ru]
}
\]
is the analog for the condition from Definition \ref{def:spi} that
$b$ is not as isomorphism of any connected component of $\Gamma$.
\begin{rem}
This last condition in Definition \ref{def:spq}, which did not appear
in this correct form in a previous version of this paper, is necessary:
without it, the invariant vanishes for any non-primitive $w\ne1$.
Indeed, let $w\ne1$ be non-primitive, and let $W_{1}=W_{0}\otimes_{K[\langle w\rangle]}K[\F]$
be the degree-1 $w$-module where $W_{0}$ is the one-dimensional
trivial $K[\langle w\rangle]$-module. Let $Q_{1}$ be the trivial
one-dimensional $K[\F]$-module. Then $\phi_{1}\colon W_{1}\to Q_{1}$
satisfies the conditions in the definition (and $-\chi(Q_{1})=\rk(\F)-1$,
showing that $\spk(w)\le\rk(\F)-1$). For every $d\in\mathbb{Z}_{\ge1}$,
let $W_{d}$ be any degree-$d$ $w$-module such that $W_{1}\le W_{d}$
is a $w$-submodule. One can check that $\phi_{d}\colon W_{d}\twoheadrightarrow Q_{d}=\nicefrac{W_{d}}{\ker(\phi_{1})}$
satisfies almost all the conditions of Definition \ref{def:spq},
and that $\chi(Q_{d})=\chi(Q_{1})$, so allowing such modules in the
definition would make $\spk$ vanish for $w$. The only condition
that $\phi_{d}\colon W_{d}\twoheadrightarrow Q_{d}$ does not satisfy
is the new condition that no quotient of $w$-modules factors through
$\phi$: in the example, the quotient of $w$-modules $W_{d}\twoheadrightarrow\nicefrac{W_{d}}{W_{1}}$
factors through $\phi_{d}$ by
\[
W_{d}\stackrel{\phi_{d}}{\twoheadrightarrow}\nicefrac{W_{d}}{\ker\left(\phi_{1}\right)}\twoheadrightarrow\nicefrac{W_{d}}{W_{1}}.
\]
\end{rem}

As mentioned above, we know less about $\spk$ than we do about $\sp$
or other stable invariants discussed above. We do know, however, that
it is $\mathrm{Aut}\F$-invariant. We also know that $\spk\left(w\right)=0$
if $w$ is a proper power, $\spk\left(w\right)=\infty$ for every
primitive word $w$, and that $\spk\left(w\right)\in[0,\rk(\F)-1]$
for any other word $w\ne1$.  

We need to stress that we do not know whether the stable $K$-primitivity
rank is truly a new invariant. In fact, we have reasons to believe
that $\spk$, for an arbitrary field $K$, is identical to $\sp$.
For example, in \cite[Conj.~1.9]{West}, we conjecture that $\pi_{q}=\pi$
for all $q$. \medskip{}

The stable representation theory of $\gl_{\bullet}\left(q\right)$
is well understood -- see, e.g., \cite{gan2018representation,putman2017representation,ERNSTWEST2025}.
We believe that the general conjectural picture defined throughout
this paper is true also in the case of $\gl_{\bullet}\left(q\right)$
with $\spq$ the suitable stable invariant. The fact that $\beta(w,\chi)$
is rational for every $w\in\F$ and non-trivial stable character of
$\gl_{\bullet}(q)$ follows from \cite[Cor.~3.6]{ERNSTWEST2025}.
\begin{conjecture}
\label{conj:spq and stable characters}Fix a prime power $q$, and
let ${\cal I}_{q}$ denote the set of all stable irreducible characters
of $\gl_{\bullet}(q)$. Then for every $w\in\F$,
\[
\spq(w)=\inf_{\triv\ne\chi\in{\cal I}_{q}}\beta(w,\chi),
\]
and, moreover, the infimum is attained.
\end{conjecture}

\end{appendices}

\bibliographystyle{alpha}
\bibliography{stable-invariants-bib}

\noindent Doron Puder, School of Mathematical Sciences, Tel Aviv University,
Tel Aviv, 6997801, Israel\\
\texttt{doronpuder@gmail.com}~\\

\noindent Yotam Shomroni, School of Mathematical Sciences, Tel Aviv
University, Tel Aviv, 6997801, Israel\\
\texttt{yotam.shomroni@gmail.com}~\\

\noindent Danielle Ernst-West, School of Mathematical Sciences, Tel
Aviv University, Tel Aviv, 6997801, Israel\\
\texttt{daniellewest@mail.tau.ac.il}~\\

\noindent Matan Seidel, School of Mathematical Sciences, Tel Aviv
University, Tel Aviv, 6997801, Israel\\
\texttt{matanseidel@gmail.com}
\end{document}